\documentclass[hidelinks,onefignum,onetabnum]{siamart251216}
\usepackage{aligned-overset}
\usepackage{amsmath}
\usepackage{amsfonts}
\usepackage{amssymb}
\usepackage{booktabs}
\usepackage{cancel}
\usepackage{diagbox}
\usepackage[hmargin={26mm,26mm},vmargin={30mm,35mm}]{geometry}
\usepackage{hyperref}
\usepackage{multirow}
\usepackage[position=b]{subcaption}
\usepackage[table]{xcolor}

\usepackage{loglogslopetriangle}
\usepackage{pgfplots}

\theoremstyle{remark}
\newtheorem{remark}[theorem]{Remark}

\theoremstyle{definition}
\newtheorem{assumption}{Assumption}

\newcommand{\st}{\;:\;}

\newcommand{\Real}{\mathbb{R}}

\newcommand{\truth}[1]{\langle #1\rangle}

\newcommand{\stokes}{\mathrm{s}}
\newcommand{\darcy}{\mathrm{d}}
\newcommand{\sstokes}{\mathrm{S}}
\newcommand{\ddarcy}{\mathrm{D}}

\newcommand{\bnd}{\mathrm{b}}

\newcommand{\compl}{\mathrm{c}}

\newcommand{\Th}{\mathcal{T}_h}
\newcommand{\Thd}{\mathcal{T}_h^\darcy}
\newcommand{\Ths}{\mathcal{T}_h^\stokes}
\newcommand{\Fh}{\mathcal{F}_h}
\newcommand{\Fhb}{\mathcal{F}_h^\bnd}
\newcommand{\FT}{\mathcal{F}_T}

\newcommand{\Poly}[1]{\mathcal{P}^{#1}}
\newcommand{\cGoly}[1]{\mathcal{G}^{\compl,#1}}
\newcommand{\Nedelec}[1]{\mathcal{N}^{#1}}

\newcommand{\VT}{\underline{V}_T^k}
\newcommand{\Vh}{\underline{V}_h^k}

\newcommand{\VhZ}{\underline{V}_{h,0}^k}
\newcommand{\QhZ}{Q_{h,0}^k}

\newcommand{\Ih}{\underline{I}_h^k}
\newcommand{\IT}{\underline{I}_T^k}

\newcommand{\lproj}[2]{\pi^{#1}_{\mathcal{P},#2}}
\newcommand{\nproj}[2]{\pi^{#1}_{\mathcal{N},#2}}

\newcommand{\GT}{G_T^k}
\newcommand{\PST}{P_{\stokes,T}^{k+1}}
\newcommand{\PDT}{P_{\darcy,T}^k}
\newcommand{\PDh}{P_{\darcy,h}^k}
\newcommand{\PDTz}{P_{\darcy,T}^0}
\newcommand{\DT}{D_T^k}
\newcommand{\DTz}{D_T^0}

\newcommand{\norm}[2]{\|#2\|_{#1}}
\newcommand{\seminorm}[2]{|#2|_{#1}}

\newcommand{\tnorm}[2]{|\kern-0.25ex|\kern-0.25ex|#2|\kern-0.25ex|\kern-0.25ex|_{#1}}

\newcommand{\term}{\mathfrak{T}}

\newcommand{\huline}[1]{\widehat{\underline{#1}}}

\newcommand{\Err}{\mathcal{E}}

\newcommand{\CfT}{C_{\mathrm{f},T}(u)}
\newcommand{\CfOmg}{C_{\mathrm{f},\Omega}}
\newcommand{\kappaT}{\kappa_T(u)}

\newcommand{\EnergyError}[1][q_r]{\alpha_\mu \norm{\mu,r,h}{\underline{e}_h}^{#1} + \norm{\nu,h}{\underline{e}_h}^2}

\title{A Hybrid High-Order method for the power-law Brinkman problem with robust error estimates in all regimes\thanks{%
    \funding{%
      The work of D. Casta\~{n}\'on Quiroz was partially supported by UNAM-PAPIIT grant IA101926. Funded by the European Union (ERC Synergy, NEMESIS, project number 101115663).
      Views and opinions expressed are however those of the authors only and do not necessarily reflect those of the European Union or the European Research Council Executive Agency. Neither the European Union nor the granting authority can be held responsible for them.
    }%
  }
}

\headers{A HHO method for the power-law Brinkman problem}{D. Casta\~{n}\'on Quiroz, D. A. Di Pietro, J. Droniou, and M. Salah}

\author{%
  D. Casta\~{n}\'on Quiroz\thanks{Instituto de Investigaciones en Matemáticas Aplicadas y en Sistemas, Universidad Nacional Autónoma de México, Circuito Escolar s/n, Ciudad Universitaria C.P. 04510 Cd. Mx. (México) (\email{daniel.castanon@iimas.unam.mx})}%
  \and
  D. A. Di Pietro\thanks{IMAG, Univ Montpellier, CNRS, Montpellier 34090, France (\email{daniele.di-pietro@umontpellier.fr})}%
  \and
  J. Droniou\thanks{IMAG, Univ Montpellier, CNRS, Montpellier 34090, France and School of Mathematics, Monash University, Melbourne, Australia (\email{jerome.droniou@cnrs.fr})}%
  \and
  M. Salah\thanks{IMAG, Univ Montpellier, CNRS, Montpellier 34090, France (\email{marwa.salah@umontpellier.fr})}
}

\ifpdf
\hypersetup{
  pdftitle={A Hybrid High-Order method for the power-law Brinkman problem with robust error estimates in all regimes},
  pdfauthor={C. Casta\~{n}\'on Quiroz, D. A. Di Pietro, J. Droniou, and M. Salah}
}
\fi

\begin{document}

\maketitle

\begin{abstract}
  In this work we propose and analyze a new Hybrid High-Order method for the Brinkman problem for fluids with power-law viscosity.
  The proposed method supports general meshes and arbitrary approximation orders and is robust in all regimes, from pure (power-law) Stokes to pure Darcy.
  Robustness is reflected by error estimates that distinguish the contributions from Stokes- and Darcy-dominated elements as identified by an appropriate dimensionless number, and that additionally account for pre-asymptotic orders of convergence.
  Theoretical results are illustrated by a complete panel of numerical experiments.
\end{abstract}

\begin{keywords}
  $p$-Brinkman, Darcy, $p$-Stokes, Hybrid High-Order methods, polytopal methods
\end{keywords}

\begin{MSCcodes}
  65N30, 65N08, 76S05, 76D07
\end{MSCcodes}

%------------------------------------------------------------------------------%

\section{Introduction}

The Brinkman problem arises in the modeling of fluid flows through a porous matrix where obstacles are present; see \cite{Allaire:91,Allaire:91*1} and also \cite{Allaire:91*2} for a review.
Mathematically, it translates into a partial differential equation which is a superposition of the Stokes and Darcy systems. %% , the latter being recovered as a singular limit.
The extension of the Brinkman model to power-law non-Newtonian fluids has been recently considered in~\cite{Anguiano.Suarez-Grau:21}.
The purpose of the present paper is to design an arbitrary-order numerical method  robust in all regimes for a version of the latter model simplified, but yet conceived to retain all the major mathematical difficulties.
We are, in particular, concerned with the robustness of the method in all regimes, and particularly the ability of the method to handle the singular limit corresponding to the Darcy problem and with the (pre-asymptotic) orders of convergence in the Darcy-dominated regime; see Remark~\ref{rem:regimes} on this subject.

Already in the standard Newtonian case, the construction of robust numerical approximations of the Brinkman problem is not straightforward; see, e.g.,~\cite{Mardal.Tai.ea:02,Burman.Hansbo:05,Burman.Hansbo:07,Juntunen.Stenberg:10,Konno.Stenberg:11,Anaya.Gatica.ea:15,Evans.Hughes:13,Caceres.Gatica.ea:17,Vacca:18,Araya.Harder.ea:17,Zhao.Chung.ea:20,Mora.Reales.ea:22} and also \cite{Gatica.Munar.ea:18} for a nonlinear version in the Hilbertian framework.
In \cite{Botti.Di-Pietro.ea:18}, the authors introduced a Hybrid High-Order (HHO) method on matching simplicial meshes and derived what appears to be the first error estimate accounting for the local regime through a (dimensionless) coefficient of friction.
Thanks to the presence of cutoff factors, this estimate is valid in the entire spectrum of regimes, ranging from the pure Stokes problem all the way to the pure Darcy problem.
The previous work was later extended in~\cite{Di-Pietro.Droniou:23} to cover general polyhedral meshes.
The key idea to obtain error estimates valid in all regimes was, in this case, to use different approximations of the velocity depending on the local value of the coefficient of friction: an HHO approximation was used in Stokes-dominates elements, while a Discrete de Rham (DDR) approximation was considered for Darcy-dominated elements.
Unfortunately, this strategy does not carry out to the non-Newtonian case, since the local coefficient of friction depends on the velocity itself and is therefore no longer computable.
In this work, we present a new method that removes this difficulty by showing that a DDR approximation can be used across all regimes for both the Darcy and source terms without negatively impacting on the accuracy of the method; the HHO approximation remains the key to handle the Stokes terms.
Inspired by the techniques of~\cite{Beirao-da-Veiga.Di-Pietro.ea:24}, we prove error estimates that are again robust across the entire spectrum of local regimes and that deliver orders of convergence that coincide with those found in \cite{Di-Pietro.Droniou:23} for the Newtonian Brinkman problem, with those of~\cite{Botti.Castanon-Quiroz.ea:21} in the (power-law) Stokes limit, and with those of~\cite{Di-Pietro.Ern:15} in the Darcy limit.
For the non-linear case corresponding to a power-law exponent different from $2$, the contributions to the error from the Stokes and Darcy terms do not converge at the same rate.
Our estimates account for this fact by predicting pre-asymptotic convergence rates corresponding to the Darcy-regime.
This is an important and new element since, in practical situations, one may not be able to use meshes fine enough for asymptotic convergence rates corresponding to the Stokes-dominated regime to be observed.
We notice, for the sake of completeness, that the convergence rates for the Stokes-dominated regime could in all likelihood be improved adapting the techniques of~\cite{Di-Pietro.Droniou.ea:21}.
This would, however, add another layer of complexity in arguments that are already very technical, which is why we have preferred to stick to (relatively) simpler estimates inspired by~\cite{Di-Pietro.Droniou:17*1}.
We refer to Remark~\ref{rem:conv.rates} for a comparison with other methods and techniques.%

  The scheme proposed in this work inherits all the key features of HHO and DDR methods including, in particular, the possibility to handle general meshes made of polytopal elements, as well as possibility to increase the approximation order to leverage the regularity of the solution when this is available or to perform non-conforming mesh adaptation in the presence of singularities.
  In the practical implementation, one can additionally benefit from static condensation, which applies to both velocity and pressure (see, e.g.,~\cite[Section~6]{Di-Pietro.Ern.ea:16}) and can therefore dramatically reduce the size of the linear systems to be solved at each non-linear iteration. %

The rest of this work is organized as follows.
In Section~\ref{sec:setting} we recall the setting (continuous problem, meshes, polynomial spaces).
The discrete problem and the main convergence result are stated in Section~\ref{sec:discrete.problem}.
Section~\ref{sec:error.analysis} contains the proof of the convergence result.
Finally, a set of numerical results supporting the theory is provided in Section~\ref{sec:numerical.results}.

%------------------------------------------------------------------------------%

\section{Setting}\label{sec:setting}

\subsection{Continuous setting}\label{sec:setting:continuous}

Let $d \in \{2,3\}$, a real exponent $r > 1$,
and a function $\mu : \Omega \to \Real$ such that $\mu \ge \underline{\mu} > 0$ in $\Omega$ (with $\underline{\mu}$ constant)
be given. For $\mathbb{T} \in \left\{ \Real^d, \Real^{d\times d}\right\}$, let $\sigma : \mathbb{T} \to \mathbb{T}$ be such that
\begin{equation}\label{eq:sigma}
  \sigma(\tau) = \mu |\tau|^{r-2} \tau \qquad \forall \tau \in \mathbb{T}.
\end{equation}
Given a non-negative function $\nu : \Omega \to \lbrack 0,\infty)$
as well as forcing terms $f : \Omega \to \Real^d$ and $g : \Omega \to \Real$ satisfying $\int_\Omega g = 0$, we consider the following problem:
Find $u : \Omega \to \Real^d$ and $p : \Omega \to \Real$ such that
\begin{subequations}\label{eq:strong}
  \begin{alignat}{2}\label{eq:strong:momentum}
    - \nabla \cdot \sigma(\nabla u) + \nu u + \nabla p &= f
    &\qquad& \text{in $\Omega$},
    \\ \label{eq:strong:mass}
    \nabla \cdot u &= g
    &\qquad& \text{in $\Omega$},
    \\ \label{eq:strong:bc}
    u &= 0
    &\qquad& \text{on $\partial\Omega$},
    \\ \label{eq:strong:closure}
    \int_\Omega p &= 0.
  \end{alignat}
\end{subequations}
For the sake of simplicity, throughout the rest of this work we assume that $\mu$ and $\nu$ are piecewise constant on a polytopal partition $P_\Omega$ of the domain. The weak formulation of \eqref{eq:strong} is: Find $(u,p)\in W^{1,r}_0(\Omega)^d\times L^{r'}(\Omega)$ such that $\int_\Omega p=0$ and
\begin{equation}\label{eq:weak}
  \begin{alignedat}{2}
    \int_\Omega \sigma(\nabla u):\nabla v
    + \int_\Omega \nu u\cdot v
    -  \int_\Omega p \, (\nabla\cdot v) &= \int_\Omega f\cdot v
    &\qquad&\forall v\in W^{1,r}_0(\Omega)^d,
    \\
    \int_\Omega q \, (\nabla \cdot u) &= \int_\Omega g \, q
    &\qquad& \forall q\in L^{r'}(\Omega).
  \end{alignedat}
\end{equation}

\begin{remark}[Choice of model problem]
  Equation~\ref{eq:strong} has to be intended as a simplified version of physically relevant problems that retains all the major mathematical difficulties.
  A more physically accurate version of the Stokes term would be obtained replacing $\sigma( \nabla u)$ with $\mu | \varepsilon(u) |^{r-2} \varepsilon(u)$ with $\varepsilon(u) \coloneq \frac12\left(\nabla u + \nabla u^\top\right)$ denoting the strain tensor.
  The extension of the present method to this case is classical; see~\cite{Botti.Di-Pietro.ea:17} and also \cite{Botti.Di-Pietro.ea:19*3} and \cite[Chapter~7]{Di-Pietro.Droniou:20} for further insight.
\end{remark}

\subsection{Discrete setting}

We consider a polytopal mesh $(\Th,\Fh)$ matching the geometrical requirements detailed in \cite[Definition 1.4]{Di-Pietro.Droniou:20}, with $\Th$ set of elements and $\Fh$ set of faces.
We additionally assume that $\Th$ is compatible with $P_\Omega$, meaning that, for each $T\in\Th$, there exists $\omega\in P_\Omega$ such that $T\subset\omega$.
We then set $\mu_T \coloneq \mu_{|T}$ and $\nu_T\coloneq \nu_{|T}$ for all $T\in\Th$, noticing that these constant values are uniquely defined in each element.
For any $Y \in \Th \cup \Fh$, we denote by $h_Y$ its diameter, so that $h=\max_{T\in\Th}h_T>0$.
For every mesh element $T\in\Th$, we denote by $\FT$ the subset of $\Fh$ containing the faces that lie on the boundary $\partial T$ of $T$ and, for all $F \in \FT$, we denote by $n_{TF}$ the unit vector normal to $F$ and pointing out of $T$.
We also define the vector field $n_{\partial T} : \partial T \to \Real^d$ such that $(n_{\partial T})_{|F} \coloneq n_{TF}$ for all $F \in \FT$.
Boundary faces lying on $\partial\Omega$  are collected in the set $\Fhb$.

From this point on, $a \lesssim b$ means $a\le Cb$ with $C$ only depending on $\Omega$, the mesh regularity parameter of \cite[Definition 1.9]{Di-Pietro.Droniou:20}, and, possibly, polynomial degrees or Lebesgue exponents involved in the expressions.
We stress that this means, in particular, that $C$ is independent of the mesh size $h$, of $\mu$, of $\nu$, and of the exact solution.
The notation $a\simeq b$ means that $a\lesssim b$ and $b\lesssim a$.

Given $Y\in\Th\cup\Fh$ and an integer $m\ge 0$, we denote by $\Poly{m}(Y)$ the space spanned by the restriction to $Y$ of $d$-variate polynomials of total degree $\le m$ and we conventionally set $\Poly{-1}(Y) \coloneq \{ 0 \}$.
The $L^2$-orthogonal projector on $\Poly{m}(Y)$ (or $\Poly{m}(Y)^d$ or $\Poly{m}(Y)^{d\times d}$) is denoted by $\lproj{m}{Y}$.
For all $T \in \Th$, we denote by $\Poly{k}(\partial T)$ the space of broken polynomials of total degree $\le k$ on $\partial T$ and by $\lproj{k}{\partial T}$ the corresponding $L^2$-orthogonal projector.
The set of broken polynomials of total degree $\le m$ on the mesh is denoted by $\Poly{m}(\Th)$, and its $L^2$-orthogonal projector by $\lproj{m}{h}$.
For any $T \in \Th$, we will additionally need the Nédélec space $\Nedelec{m}(T)\subset \Poly{m}(T)^d$ such that
\[
\begin{gathered}
\Nedelec{m}(T)
\coloneq
\nabla \Poly{m}(T) \oplus \cGoly{m}(T)
\\
\text{
  with
  $\cGoly{m}(T) \coloneq
  \begin{cases}
    (x - x_T)^\perp \, \Poly{m-1}(T) & \text{if $d = 2$},
    \\
    (x - x_T) \times \Poly{m-1}(T) & \text{if $d = 3$},
  \end{cases}
  $
}
\end{gathered}
\]
where $y^\perp$ denotes the rotation by $-\frac{\pi}{2}$ of a vector $y \in \Real^2$,
while $x_T$ is a point inside $T$.
The $L^2$-orthogonal projector on $\Nedelec{m}(T)$ will be denoted by $\nproj{m}{T}$.
By definition, $\Nedelec{0}(T) = \{ 0 \}$ and, since $\Poly{m-1}(T)^d \subset \Nedelec{m}(T)$ for all $m\ge 0$, we have
\begin{equation}\label{eq:nproj.circ.lproj}
  \lproj{m-1}{T}\circ \nproj{m}{T} = \lproj{m-1}{T}.
\end{equation}

Let us recall that the conjugate exponent of $q\in [1,\infty]$ is $q'\in [1,\infty]$ such that
\begin{equation}\label{eq:r'}
  \frac1q+ \frac1{q'} = 1.
\end{equation}

\begin{lemma}[Bound on the $L^q$-norm of a polynomial]
  Let $Y \in \Th \cup \Fh$.
  For any exponent $q \in [1,+\infty]$ and any polynomial $\varphi \in \Poly{m}(Y)$, it holds
  \begin{equation}\label{eq:norm.Lq.poly}
    \norm{L^q(Y)}{\varphi} \lesssim \sup_{\psi \in \Poly{m}(Y) \setminus \{0\}} \frac{\int_Y \varphi~\psi}{\norm{L^{q'}(Y)}{\psi}}.
  \end{equation}
\end{lemma}

\begin{proof}
  It classically holds
  \[
  \norm{L^q(Y)}{\varphi} = \sup_{\psi \in L^{q'}(\Omega) \setminus \{0\}} \frac{\int_Y \varphi~\psi}{\norm{L^{q'}(Y)}{\psi}}.
  \]
  For $q<\infty$ this relation is a direct consequence of the fact that $L^{q'}(Y)$ is the dual space of $L^q(Y)$ while, for $q=\infty$, the equality is achieved by considering the sequence $(\psi_n)_n=(\chi_n\mathrm{sign}(\varphi))$ where $\chi_n$ is the characteristic functions of the set $\left\{x\in Y\,:\,|\varphi(x)|\ge \norm{L^\infty(Y)}{\varphi}-\frac1n\right\}$.

  For $\psi\in L^{q'}(Y) \setminus \{0\}$, we can write, denoting by $\$$ the supremum in the right-hand side of~\eqref{eq:norm.Lq.poly},
  \[
  \int_Y \varphi~\psi = \int_Y \varphi~\lproj{m}{Y}\psi
  \le \$ ~ \norm{L^{q'}(Y)}{\lproj{m}{Y}\psi}
  \lesssim \$ ~ \norm{L^{q'}(Y)}{\psi},
  \]
  where the first equality follows from the definition of $\lproj{m}{Y}$ and $\varphi\in\Poly{m}(T)$, while the last bound is a consequence of the $L^{q'}$-boundedness of the $L^2$-orthogonal projector proved, e.g., in~\cite[Lemma~1.44]{Di-Pietro.Droniou:20}.
  Dividing both sides by $\norm{L^{q'}(Y)}{\psi}$ and passing to the supremum over $\psi$ concludes the proof.
\end{proof}

%------------------------------------------------------------------------------%

\section{Discrete problem and main results}\label{sec:discrete.problem}

This section contains the definition of the scheme and the statement of the main results of the analysis.

\subsection{Discrete spaces}

Let a polynomial degree $k \ge 0$ be fixed.
We define the following variation of the HHO space for the velocity, in which the element unknowns are taken in $\Nedelec{k}(T)$ instead of $\Poly{k}(T)^d$:
\[
\begin{aligned}
  \Vh \coloneq \big\{
  &\underline{v}_h = ( (v_T)_{T \in \Th}, (v_F)_{F \in \Fh}) \st
  \\
  &\text{$v_T \in \Nedelec{k}(T)$ for all $T \in \Th$,
    $v_F \in \Poly{k}(F)^d$ for all $F \in \Fh$}
  \big\},
\end{aligned}
\]
as well as its subspace with strongly enforced boundary conditions:
\[
\VhZ \coloneq \left\{
\underline{v}_h \in \Vh \st \text{%
  $v_F = 0$ for all $F \in \Fhb$
}
\right\}.
\]
The interpolator $\Ih : W^{1,1}(\Omega)^d \to \Vh$ is such that, for all $v \in W^{1,1}(\Omega)^d$,
\begin{equation}\label{eq:def.Ih}
  \Ih v
  \coloneq
  ( (\nproj{k}{T} v)_{T \in \Th}, (\lproj{k}{F} v)_{F \in \Fh} ).
\end{equation}
The restrictions of $\Vh$, of $\underline{v}_h \in \Vh$, and of the interpolator $\Ih$ to a mesh element $T \in \Th$ are denoted replacing the subscript $h$ with $T$.
Given $\underline{v}_T = (v_T, (v_F)_{F \in \FT}) \in \VT$, we let $v_{\partial T} \in \Poly{k}(\partial T)^d$ be such that $(v_{\partial T})_{|F} \coloneq v_F$ for all $F \in \FT$.

\begin{remark}[Element component for $k=0$]
  For $k=0$, since $\Nedelec{0}(T) = \{0\}$, there are no element values attached to $\underline{v}_h\in\Vh$. In this case, as usual in HHO methods \cite[Section 5.1.1]{Di-Pietro.Droniou:20}, we let
  \[
  v_T \coloneq \frac{1}{|T|}\sum_{F\in\FT} \varpi_{TF} |F|v_F,
  \]
  with non-negative weights $(\varpi_{TF})_{F\in\FT}$ such that $\sum_{F\in\FT}|F|\varpi_{TF}=|T|$.
  The standard choice of weights on elements that are star-shaped with respect to their center of mass $\overline{x}_T$ is obtained setting, for all $F \in \FT$, $\varpi_{TF} \coloneq \frac{|F| d_{TF}}{d |T|}$ with $d_{TF}$ denoting the distance of $\overline{x}_T$ from the plane containing $F$.
\end{remark}

For the pressure, we define the following space of fully discontinuous polynomials with zero-average over $\Omega$:
\[
\QhZ \coloneq \left\{
q_h \in L^{r'}(\Omega) \st \text{%
  $q_T \coloneq (q_h)_{|T} \in \Poly{k}(T)$ for all $T \in \Th$
and $\int_\Omega q_h=0$}
\right\}.
\]

\subsection{Operator reconstructions and discrete functions}

\subsubsection{Stokes term}

For all $T \in \Th$, denote by $\GT : \VT \to \Poly{k}(T)^{d\times d}$ and $\PST : \VT \to \Poly{k+1}(T)^d$ the usual HHO gradient and Stokes velocity reconstructions such that, for all $\underline{v}_T \in \VT$ and, respectively, all $\tau \in \Poly{k}(T)^{d\times d}$, and all $w \in \Poly{k+1}(T)^d$,
\[
\begin{gathered}
  \int_T \GT \underline{v}_T : \tau
  = -\int_T v_T \cdot (\nabla \cdot \tau)
  + \int_{\partial T} v_{\partial T} \cdot \tau n_{TF},
  \\
  \int_T \nabla \PST \underline{v}_T : \nabla w
  + \beta_T \int_T \PST \underline{v}_T {\cdot} \int_T w
  = \int_T \GT \underline{v}_T : \nabla w
  + \beta_T \int_T v_T {\cdot} \int_T w,
\end{gathered}
\]
where the purpose of the factor $\beta_T \simeq \frac1{h_T^2 |T|}$ is to ensure that all the terms have the same scaling.
$\PST$ is called the \emph{Stokes velocity reconstruction} because it is used to discretise the (nonlinear) Stokes component in the model.

By~\eqref{eq:nproj.circ.lproj}, we have the usual commutation property $\GT \circ \IT = \lproj{k}{T} \circ \nabla$,
while $\PST \circ \IT$ coincides with the elliptic projector onto $\Poly{k+1}(T)^d$ (or the modified projector of \cite[Definition 5.4]{Di-Pietro.Droniou:20} if $k=0$).
We define the stabilization function $s_{\mu,h} : \Vh \times \Vh \to \Real$ such that, for all $(\underline{w}_h, \underline{v}_h) \in \Vh \times \Vh$,
\begin{equation}\label{eq:jh.jT}
  \begin{aligned}
    s_{\mu,h}(\underline{w}_h, \underline{v}_h)
    &\coloneq \sum_{T \in\Th} s_{\mu,T}(\underline{w}_T, \underline{v}_T),
    \\
    s_{\mu,T}(\underline{w}_T, \underline{v}_T)
    &\coloneq
    h_T^{1-r} \int_{\partial T} \sigma(\Delta_{\partial T}^k \underline{w}_T) \cdot \Delta_{\partial T}^k \underline{v}_T,
  \end{aligned}
\end{equation}
where
\begin{equation}\label{eq:Delta}
  \Delta_{\partial T}^k \underline{v}_T
  \coloneq
  \lproj{k}{\partial T}(\PST \underline{v}_T - v_{\partial T})
  - \nproj{k}{T}(\PST \underline{v}_T - v_T).
\end{equation}
Reasoning as in \cite[Eq.~(6.48)]{Di-Pietro.Droniou:20}, with $p$ replaced by $q$ and $h_F$ by $h_T$ in the definition of the seminorm, it can be proved that, for all $q > 1$ and all $w \in W^{m+2,q}(T)^d$ with $m \in \{0,\ldots,k\}$,
\begin{equation}\label{eq:Delta.k:consistency}
  \norm{L^q(\partial T)^d}{\Delta_{\partial T}^k \IT w}
  \lesssim h_T^{m + 2 - \frac1q} \seminorm{W^{m+2,q}(T)^d}{w}.
\end{equation}
Moreover, for all $w \in W^{1,\infty}(T)^d$, using triangle inequalities followed by the $L^\infty$-boundedness of $\lproj{k}{\partial T}$ and $\nproj{k}{T}$ and the $L^\infty$-approximation properties of the elliptic projector (see \cite[Lemma 1.44 and Theorem~1.48]{Di-Pietro.Droniou:20}), it holds
\begin{equation}\label{eq:Delta.k:boundedness}
  \norm{L^\infty(\partial T)^d}{\Delta_{\partial T}^k \IT w}
  \lesssim h_T \norm{L^\infty(T)^{d\times d}}{\nabla w}.
\end{equation}

To discretise the Stokes term, we define the function $a_{\mu,h} : \Vh \times \Vh \to \Real$ such that, for all $(\underline{w}_h, \underline{v}_h) \in \Vh \times \Vh$,
\begin{equation}\label{eq:a.mu.h.T}
  \begin{gathered}
    a_{\mu,h}(\underline{w}_h, \underline{v}_h)
    \coloneq \sum_{T \in \Th} a_{\mu,T}(\underline{w}_T, \underline{v}_T),
    \\
    a_{\mu,T}(\underline{w}_T, \underline{v}_T)
    \coloneq \int_T \sigma(\GT \underline{w}_T) \cdot \GT \underline{v}_T
    + s_{\mu,T}(\underline{w}_T, \underline{v}_T).
  \end{gathered}
\end{equation}

\begin{remark}[Stabilisation term]
  As in any polytopal method, the role of the stabilisation term $s_{\mu,T}$ is to bound the degrees of freedom that are not bounded by the consistent part (volumetric integral) of the function $a_{\mu,T}$.
  As per \cite[Lemma 2.11]{Di-Pietro.Droniou:20}, to ensure that $s_{\mu,T}$ does not impact the polynomial consistency of $a_{\mu,T}$, we make it dependent only on $\IT\PST\underline{w}_T-\underline{w}_T$.
  The particular form we have chosen in \eqref{eq:Delta} corresponds in the linear case to the original HHO stabilisation \cite[Proposition 2.13]{Di-Pietro.Droniou:20}. In the nonlinear case, using $\sigma$ in the definition of $s_{\mu,T}$ ensures that this stabilisation inherits appropriate monotonicity and scaling properties, and was already done for the $p$-Laplacian in \cite{Di-Pietro.Droniou:17,Di-Pietro.Droniou:20}.
\end{remark}

\subsubsection{Darcy term}

Denote by $\DT \coloneq \operatorname{tr} \GT$ the usual HHO divergence, which satisfies, for all $\underline{v}_T \in \VT$,
\begin{equation}\label{eq:DT}
  \int_T \DT \underline{v}_T\, q
  = \int_T (\nabla \cdot v_T) \, q
  + \int_{\partial T} (v_{\partial T} - v_T) \cdot n_{\partial T}\, q
  \qquad \forall q \in \Poly{k}(T).
\end{equation}
Since $\GT\circ\IT=\lproj{k}{T}\circ\nabla$, we have the commutation property
\begin{equation}\label{eq:DT.proj.div}
  \DT\IT w = \lproj{k}{T} (\nabla\cdot w) \qquad \forall w\in W^{1,1}(T)^d.
\end{equation}
Inspired by the Discrete de Rham method~\cite{Di-Pietro.Droniou:23*1} (see also~\cite{Di-Pietro.Ern:17,Di-Pietro.Droniou:23}), we define the \emph{Darcy velocity reconstruction} $\PDT:\VT\to\Poly{k}(T)^d$ such that, for all $\underline{v}_T \in \VT$ and all $(\varphi,w) \in \Poly{k+1}(T) \times \cGoly{k}(T)$,
\begin{equation}\label{eq:PDT}
  \int_T \PDT \underline{v}_T \cdot (\nabla \varphi + w)
  = -\int_T \DT \underline{v}_T\, \varphi
  + \int_{\partial T} (v_{\partial T} \cdot n_{\partial T})\, \varphi
  + \int_T v_T \cdot w.
\end{equation}
Applying this definition with $(\varphi,w)\in\Poly{k}(T)\times\cGoly{k-1}(T)$, using \eqref{eq:DT} with $q=\varphi$ and the right-hand side integrated by parts, and recalling that $\Poly{k-1}(T)^d=\nabla \Poly{k}(T)\oplus \cGoly{k-1}(T)$, it can be checked that
\begin{equation}\label{eq:pi.k-1.PDT=vT}
  \int_T(\PDT\underline{v}_T-v_T) \cdot \Psi = 0
  \qquad \forall (\underline{v}_T, \Psi) \in \VT \times \Poly{k-1}(T)^d.
\end{equation}

For the Darcy term, we define the bilinear form $a_{\nu,h} : \Vh \times \Vh \to \Real$ such that, for all $(\underline{w}_h, \underline{v}_h) \in \Vh \times \Vh$,
\begin{equation}\label{eq:adh}
  a_{\nu,h}(\underline{w}_h, \underline{v}_h)
  \coloneq \sum_{T \in \Th} \nu_T a_{\darcy,T}(\underline{w}_T, \underline{v}_T),
\end{equation}
with, for all $T \in \Th$,
\begin{equation}\label{eq:adT}
  \begin{gathered}
    a_{\darcy,T}(\underline{w}_T, \underline{v}_T)
    \coloneq \int_T \PDT \underline{w}_T \cdot \PDT \underline{v}_T
    + s_{\darcy,T}(\underline{w}_T, \underline{v}_T),
    \\
    \begin{aligned}[t]
      s_{\darcy,T}(\underline{w}_T, \underline{v}_T)
      &\coloneq \lambda_T \int_T \nproj{k}{T} (\PDT \underline{w}_T - w_T) \cdot \nproj{k}{T} (\PDT \underline{v}_T - v_T)
      \\
      &\quad
        + h_T \sum_{F \in \FT \setminus \Fhb} \int_F (\PDT \underline{w}_T - w_F) \cdot (\PDT \underline{v}_T - v_F)
    \end{aligned}
  \end{gathered}
\end{equation}
where the role of $\lambda_T \coloneq \frac{h_T^d}{|T|} \operatorname{card}(\FT) \simeq 1$ is to balance out the two contributions \cite[Section 4]{Di-Pietro.Droniou:23}.

\subsubsection{Velocity-pressure coupling term}

The pressure-velocity coupling is realized by the bilinear form $b_h : \Vh \times \QhZ \to \Real$ such that, for all $(\underline{v}_h, q_h) \in \Vh \times \QhZ$,
\begin{equation}\label{eq:bh}
  b_h(\underline{v}_h, q_h)
  \coloneq \sum_{T \in \Th} b_T(\underline{v}_T, q_T),\qquad
  b_T(\underline{v}_T, q_T) \coloneq -\int_T \DT \underline{v}_T\,q_T.
\end{equation}
By \eqref{eq:DT.proj.div}, for all $w \in W^{1,1}(\Omega)^d$,
\begin{equation}\label{eq:fortin}
  b_h(\Ih w, q_h) = -\int_\Omega (\nabla \cdot w) \, q_h
  \qquad \forall q_h \in \Poly{k}(\Th).
\end{equation}

\subsection{Discrete problem}

The discrete problem reads: Find $(\underline{u}_h, p_h) \in \VhZ \times \QhZ$ such that
\begin{subequations}\label{eq:discrete}
  \begin{alignat}{2}\label{eq:discrete:momentum}
    a_{\mu,h}(\underline{u}_h, \underline{v}_h)
    + a_{\nu,h}(\underline{u}_h, \underline{v}_h)
    + b_h(\underline{v}_h, p_h)
    &= \int_\Omega f \cdot \PDh \underline{v}_h
    &\qquad& \forall \underline{v}_h \in \VhZ,
    \\ \label{eq:discrete:mass}
    -b_h(\underline{u}_h, q_h) &= \int_\Omega g \, q_h
    &\qquad& \forall q_h \in \QhZ,
  \end{alignat}
\end{subequations}
where $\PDh \underline{v}_h \in \Poly{k}(\Th)^d$ is such that $(\PDh \underline{v}_h)_{|T} \coloneq \PDT \underline{v}_T$ for all $T \in \Th$.

\begin{remark}[A variation with element velocities in full polynomial spaces]
  A variation of the method with element velocities in $\Poly{k}(T)^d$ instead of $\Nedelec{k}(T)$ is obtained replacing $\nproj{k}{T}$ by $\lproj{k}{T}$ in the definitions~\eqref{eq:Delta} of $\Delta_{\partial T}^k$ and~\eqref{eq:adT} of $s_{\darcy,T}$.
While the discrete spaces for this variation are the same as for the method of~\cite{Di-Pietro.Droniou:23}, the present scheme for $r=2$ differs from the one in the reference since the velocity reconstruction $\PDT$ is used in both the Darcy and forcing terms irrespectively of the regime.
\end{remark}

\subsection{Norms and local coefficient of friction}

\subsubsection{Discrete $L^2$-norms}

We define the Darcy norm such that, for all $\underline{v}_h \in \Vh$,
\begin{equation}\label{eq:norm.darcy}
  \norm{\nu,h}{\underline{v}_h}
  \coloneq \left(
  \sum_{T \in \Th} \nu_T \norm{\darcy,T}{\underline{v}_T}^2
  \right)^{\frac12},\qquad
  \norm{\darcy,T}{\underline{v}_T}
  \coloneq a_{\darcy,T}(\underline{v}_T,\underline{v}_T)^{\frac12}.
\end{equation}
Moreover, for all $T \in \Th$, we define on $\VT$ the discrete $L^2$-norm $\norm{0,2,T}{{\cdot}}$ such that
\begin{equation}\label{eq:norm.0.2.T}
  \norm{0,2,T}{\underline{v}_T} \coloneq \left(
  \norm{L^2(T)^d}{v_T}^2
  + h_T \norm{L^2(\partial T)^d}{v_{\partial T} - v_T}^2
  \right)^{\frac12}\qquad
  \forall \underline{v}_T \in \VT.
\end{equation}

\begin{proposition}[Discrete $L^2$- and Darcy-norms]
  For all $T \in \Th$ and all $\underline{v}_T = (v_T, (v_F)_{F \in \FT}) \in \VT$ such that $v_F = 0$ for all $F \in \FT \cap \Fhb$, it holds
  \begin{equation}\label{eq:L2.norm.darcy}
    \norm{0,2,T}{\underline{v}_T} \lesssim \norm{\darcy,T}{\underline{v}_T}.
  \end{equation}
\end{proposition}

\begin{proof}
  Let $\underline{v}_T \in \VT$. Using the definition \eqref{eq:norm.0.2.T} of the discrete $L^2$-norm together with triangle and discrete trace inequalities, we have
  \begin{equation}\label{eq:equiv.norms.01}
    \norm{0,2,T}{\underline{v}_T}
    \lesssim \norm{L^2(T)^d}{v_T}
    + \left(
    h_T \sum_{F \in \FT \setminus \Fhb} \norm{L^2(F)^d}{v_F}^2.
    \right)^{\frac12}
    \eqcolon \term_1 + \term_2.      
  \end{equation}
  Notice that boundary faces have been removed from the sum using $v_F = 0$ for all $F \in \FT \cap \Fhb$.
  For the boundary term, we insert $\pm \PDT \underline{v}_T$ and use a triangle inequality together with a discrete trace inequality:
  \begin{equation}\label{eq:equiv.norms.1}
    \begin{aligned}
      \term_2
      &\lesssim\left(
      h_T \sum_{F \in \FT \setminus \Fhb} \norm{L^2(F)^d}{\PDT \underline{v}_T - v_F}^2
      \right)^{\frac12}
      + \norm{L^2(T)^d}{\PDT \underline{v}_T}
      \\
      \overset{\eqref{eq:norm.darcy},\,\eqref{eq:adT}}&\lesssim
      \norm{\darcy,T}{\underline{v}_T}.  
    \end{aligned}
  \end{equation}
  For the volumetric term, we insert $\pm \nproj{k}{T} \PDT \underline{v}_T$, then use a triangle inequality to get
  \begin{align}
    \term_1
    &\le \norm{L^2(T)^d}{ \nproj{k}{T} \PDT \underline{v}_T}
    + \norm{L^2(T)^d}{\nproj{k}{T} \PDT \underline{v}_T - v_T}\nonumber\\
    &\le \norm{L^2(T)^d}{ \PDT \underline{v}_T}
    + \norm{L^2(T)^d}{\nproj{k}{T} (\PDT \underline{v}_T - v_T)}
    \overset{\eqref{eq:norm.darcy},\,\eqref{eq:adT}}\lesssim
    \norm{\darcy,T}{\underline{v}_T},
    \label{eq:equiv.norms.001}
  \end{align}
  where the second inequality is justified by the $L^2$-contraction property of $\nproj{k}{T}$ and the fact that $v_T\in \Nedelec{k}(T)$.
  The desired result follows by plugging \eqref{eq:equiv.norms.1} and \eqref{eq:equiv.norms.001} into \eqref{eq:equiv.norms.01}.
\end{proof}

\subsubsection{Discrete $W^{1,q}$-norms}

For all $q \in (1,+\infty)$, we define on $\Vh$ the $W^{1,q}$-like seminorm such that, for all $\underline{v}_h \in \Vh$,
\begin{equation}\label{eq:norm.1.q.h}
  \begin{gathered}
    \norm{1,q,h}{\underline{v}_h}
    \coloneq \left(
    \sum_{T \in \Th} \norm{1,q,T}{\underline{v}_T}^q
    \right)^{\frac1q},
    \\
    \norm{1,q,T}{\underline{v}_T}
    \coloneq \left(
    \norm{L^q(T)^{d\times d}}{\nabla v_T}^q
    +   h_T^{1-q} \norm{L^q(\partial T)^d}{v_{\partial T} - v_T}^q
    \right)^{\frac1q}.
  \end{gathered}
\end{equation}
Discrete inverse and trace inequalities show that
\begin{equation}\label{eq:norm.1.2.T<=hT-1.norm.0.2.T}
  \norm{1,2,T}{\underline{v}_T}
  \lesssim h_T^{-1} \norm{0,2,T}{\underline{v}_T}
  \qquad \forall \underline{v}_T \in \VT.
\end{equation}
The following inequality is a consequence of the local norm equivalence of \cite[Lemma~6.11]{Di-Pietro.Droniou:20}: For all $T\in\Th$,
  \begin{equation}\label{eq:norm.equiv.GT.1qT}
    \left(\norm{L^q(T)^{d\times d}}{\GT\underline{v}_T}^q+h_T^{1-q}\norm{L^q(\partial T)^d}{\Delta_{\partial T}\underline{v}_T}^q\right)^{\frac1q}
    \lesssim \norm{1,q,T}{\underline{v}_T}
    \qquad \forall \underline{v}_T \in \VT.
  \end{equation}
We will also need, in what follows, the following $\mu$-weighted version of the above $W^{1,q}$-like seminorm:
\[
\norm{\mu,q,h}{\underline{v}_h}
\coloneq \left(
\sum_{T \in \Th} \mu_T \norm{1,q,T}{\underline{v}_T}^q
\right)^{\frac1q}
\qquad \forall \underline{v}_h \in \Vh.
\]

\begin{proposition}[Local Poincaré--Wirtinger inequality for the Darcy velocity]
  For all $q \in (1,+\infty)$, all $T \in \Th$, and all $\underline{v}_T \in \VT$, it holds
  \begin{equation}\label{eq:poincare.PDT}
    \norm{L^q(T)^d}{\PDT \underline{v}_T - v_T}
    \lesssim h_T \norm{1,q,T}{\underline{v}_T}.
  \end{equation}
\end{proposition}

\begin{proof}
  Let $\Psi\in\Poly{k}(T)^d$ and let $(\varphi,w) \in \Poly{k+1}(T) \times \cGoly{k}(T)$ such that $\Psi = \nabla \varphi + w$.
  We have
  \[
  \int_T \DT\underline{v}_T \, \varphi
  = \int_T \DT\underline{v}_T\,\lproj{k}{T}\varphi\overset{\eqref{eq:DT}}=
  \int_T (\nabla \cdot v_T) \, %% \cancel{\lproj{k}{T}}
  \varphi + \int_{\partial T} (v_{\partial T} - v_T) \cdot n_{\partial T}\, \lproj{k}{T}\varphi,
  \]
  where the introduction of the projector in the first equality is justified by its definition since $\DT\underline{v}_T\in\Poly{k}(T)$, and its cancellation in the first term in the right-hand side by the fact that $\nabla\cdot v_T\in\Poly{k}(T)$.
  Plugging this relation into the definition~\eqref{eq:PDT} of $\PDT$ and subtracting from the resulting expression $\int_T v_T \cdot (\nabla \varphi + w) = - \int_T (\nabla \cdot v_T) \, \varphi + \int_{\partial T} (v_T \cdot n_{\partial T}) \, \varphi + \int_T v_T \cdot w$, we obtain,
  \[
  \int_T (\PDT \underline{v}_T - v_T) \cdot \Psi
  = \int_{\partial T} (v_{\partial T} - v_T )\cdot n_{\partial T} \, (\varphi - \lproj{k}{T} \varphi).
  \]
  Using H\"older inequalities together with the fact that $1 \overset{\eqref{eq:r'}}= h_T^{\frac{1-q}{q}} h_T^{\frac1{q'}}$, we go on writing
  \[
  \begin{aligned}
    \int_T (\PDT \underline{v}_T - v_T) \cdot \Psi
    &\lesssim
    h_T^{\frac{1-q}{q}} \norm{L^q(\partial T)^d}{v_{\partial T} - v_T}
    \, h_T^{\frac1{q'}} \norm{L^{q'}(\partial T)}{\varphi  - \lproj{k}{T} \varphi}
    \\
    \overset{\eqref{eq:norm.1.q.h}}&\lesssim h_T
    \norm{1,q,T}{\underline{v}_T} \, \norm{L^{q'}(T)^d}{\nabla \varphi},
  \end{aligned}
  \]
  where, in the conclusion, we have additionally used the trace approximation properties of $\lproj{k}{T}$, see \cite[Theorem 1.45]{Di-Pietro.Droniou:20}.
  Reasoning as in the proof of~\cite[Lemma~8]{Di-Pietro.Droniou.ea:24}, we have that $\norm{L^{q'}(T)^d}{\nabla \varphi} \lesssim \norm{L^{q'}(T)^d}{\Psi}$.
  By~\eqref{eq:norm.Lq.poly}, dividing by $\norm{L^{q'}(T)^d}{\Psi}$ and passing to the supremum the conclusion follows.
\end{proof}

\subsubsection{Coefficient of friction}

We define the following dimensionless number, which can be interpreted as a coefficient of friction:
\begin{equation}\label{eq:CfT}
  \text{
    $\CfT \coloneq \frac{\nu_T h_T^2}{\kappaT}$
    where $\kappaT \coloneq \mu_T \norm{L^\infty(T)}{|\nabla u|^{r-2}}$.
  }
\end{equation}
If $\kappaT = 0$, we conventionally set $\CfT \coloneq +\infty$.
Notice that, if $1 < r < 2$, $\CfT = +\infty$ is only possible if $\mu_T = 0$, case considered in Remark~\ref{rem:darcy.limit} below.
Consistently, we also set $\CfT=0$ whenever $|\nabla u|^{r-2}\not\in L^\infty(T)$ and $\mu_T>0$.

Based on the above dimensionless number, we can partition the sets of elements into Darcy- and Stokes-dominated writing $\Th = \Thd(u) \sqcup \Ths(u)$ with
\[
\Thd(u) \coloneq
\left\{
T \in \Th \st \CfT \ge 1
\right\},\qquad
\Ths(u) \coloneq \left\{
T \in \Th \st \CfT < 1
\right\}.
\]
For each component of the error, we prove consistency estimates in which we identify the convergence rate of the contribution of each mesh element depending on the local regime determined by the value of $\CfT$.

\begin{remark}[Asymptotic and pre-asymptotic regimes]\label{rem:regimes}
  Notice that, as long as $\kappaT > 0$ for all $T \in \Th$, we asymptotically have $\CfT < 1$ for all $T \in \Th$, i.e., the problem is Stokes-dominated.
  In these circumstances, the Darcy regime can therefore only be observed pre-asymptotically for meshes coarse enough.
\end{remark}

\subsection{Regime-dependent velocity error estimate}

In this section we state the main result of this work, namely a regime-dependent estimate of the error on the velocity.
A key element in its proof are the monotonicity (with values $q_r$ depending on whether the power-law is strongly monotone or not, cf. \eqref{eq:def.qr} below) and coercivity properties detailed below, which are used to estimate the error starting from an error equation derived in the spirit of \cite{Di-Pietro.Droniou:18}.

\begin{lemma}[Monotonicity of $a_{\mu,h}$]\label{lem:ah:monotonicity}
  Denote by $(\underline{u}_h,p_h)\in\VhZ\times\QhZ$ the solution of \eqref{eq:discrete}
  and by $(u,p)\in W^{1,r}_0(\Omega)^d\times L^{r'}(\Omega)$ the solution of \eqref{eq:weak}.
  Let, for the sake of brevity, $\widehat{\underline{u}}_h \coloneq \Ih u$, and set
  \[
  \underline{e}_h \coloneq \underline{u}_h - \huline{u}_h.
  \]
  Then, it holds
  \begin{equation}\label{eq:ah:monotonicity}
    \alpha_\mu \norm{\mu,r,h}{\underline{e}_h}^{q_r}
    \lesssim a_{\mu,h}(\underline{u}_h, \underline{e}_h)
    -  a_{\mu,h}(\huline{u}_h, \underline{e}_h),
  \end{equation}
  where, recalling that $\underline{\mu} > 0$ is a real number such that $\mu \ge \underline{\mu}$ in $\Omega$,
  \begin{equation}\label{eq:def.qr}
  \alpha_\mu \coloneq \begin{cases}
    \underline{\mu}^{\frac{2-r}{r(r-1)}}
    & \text{if $1 < r < 2$,}
    \\
    1
    & \text{if $r \ge 2$},
  \end{cases}
  \qquad
  q_r \coloneq \begin{cases}
    2 & \text{if $1 < r < 2$,}
    \\
    r & \text{if $r \ge 2$}.
  \end{cases}
  \end{equation}
\end{lemma}

\begin{proof}
  See Section~\ref{sec:ah:monotonicity}.
\end{proof}

\begin{corollary}[Coercivity for the error]
  There exists a real number $C_{\rm stab} > 0$ independent of the mesh size, of $\mu$, and of $\nu$ such that
  \begin{equation}\label{eq:stability}
    C_{\rm stab} \left(
    \alpha_\mu \norm{\mu,r,h}{\underline{e}_h}^{q_r}
    + \norm{\nu,h}{\underline{e}_h}^2
    \right)
    \le a_{\mu,h}(\underline{u}_h, \underline{e}_h)
    -  a_{\mu,h}(\huline{u}_h, \underline{e}_h)
    + a_{\nu,h}(\underline{e}_h, \underline{e}_h).
  \end{equation}
\end{corollary}
To establish convergence rates, as usual we require regularity assumptions on the solution $(u,p)$ that are higher than the minimal ones provided by the model. These assumptions are detailed here.

\begin{assumption}[Regularity of $(u,p)$]\label{ass:regularity}
  Let $(u,p) \in W^{1,r}_0(\Omega)^d \times L^{r'}_0(\Omega)$ denote the solution to \eqref{eq:weak}.
  We assume $p \in W^{1,1}(\Omega)$ and that there exists $m \in \{0, \ldots, k\}$ such that
  \begin{itemize}
  \item $u_{|T} \in W^{m+2,r}(T)^d \cap W^{m+1,r'}(T)^d$,
    $\sigma(\nabla u)_{|T} \in W^{m+1,r'}(T)^{d\times d}$,
    and $p_{|T} \in W^{m+1,r'}(T)$ if $T \in \Ths(u)$;
  \item $u_{|T} \in H^{m+2}(T)^d \cap W^{m+2,r}(T)^d$,
    $\sigma(\nabla u)_{|T} \in H^{m+1}(T)^{d\times d}$,
    and $p_{|T} \in H^{m+2}(T)$ if $T \in \Thd(u)$.
  \end{itemize}
\end{assumption}

The error estimate is expressed in terms of the consistency error contributions listed hereafter, in which $\truth{P}\in\{0,1\}$ denotes the truth value of the logical proposition $P$.

\begin{itemize}
\item Contributions stemming from the Stokes term:
  \begin{equation}\label{eq:mathfrak.S}
    \begin{aligned}
      \mathfrak{S}_{\stokes,1}
      &\coloneq
      \sum_{T\in\Ths(u)} \Bigg(
      \begin{aligned}[t]
        &h_T^{r'(m+1)}\mu_T^{-\frac{1}{r-1}}\seminorm{W^{m+1,r'}(T)^{d\times d}}{\sigma(\nabla u)}^{r'}
        \\
        &+ h_T^{r(m+1)}\mu_T\seminorm{W^{m+2,r}(T)^d}{u}^r
        \Bigg),
      \end{aligned}
      \\
      \mathfrak{S}_{\stokes,2}
      &\coloneq
      \truth{r \ge 2} \sum_{T\in\Ths(u)} h_T^{r(m+1)} \mu_T^{-1} \kappaT^{r} \seminorm{W^{m+2,r}(T)^{d}}{u}^{r},
      \\
      \mathfrak{S}_{\darcy,1}
      &\coloneq
      \sum_{T\in\Thd(u)} h_T^{2(m+1)} \CfT^{-1} \kappaT^{-1}\seminorm{H^{m+1}(T)^{d\times d}}{\sigma(\nabla u)}^{2}
      \\
      &+ \truth{1 < r < 2} \hspace{-1.2ex} \sum_{T\in\Thd(u)} \hspace{-1.2ex}
      h_T^{r(m+1)} \CfT^{-\frac{r'}{2}} \mu_T^{r'} \kappaT^{-\frac{r'}{2}} \seminorm{W^{m+2,r}(T)^{d}}{u}^{r},
      \\
      \mathfrak{S}_{\darcy,2}
      &\coloneq
      \sum_{T\in\Thd(u)} h_T^{2(m+1)}\kappaT^{-1} \CfT^{-1} \seminorm{H^{m+1}(T)^{d\times d}}{\sigma(\nabla u)}^{2}
      \\
      &\quad
      + \truth{r \ge 2} \sum_{T\in\Thd(u)} h_T^{2(m+1)} \CfT^{-1} \kappaT \seminorm{H^{m+2}(T)^{d}}{u}^{2};
    \end{aligned}
  \end{equation}
\item Contributions stemming from the Darcy term:
  \begin{equation}\label{eq:mathfrak.D}
    \begin{aligned}
      \mathfrak{D}_{\stokes,1}
      &\coloneq
      \sum_{T \in \Ths(u)}
      h_T^{r'\left(m+\frac{2}{r'}\right)} \CfT^{\frac{1}{r-1}} \nu_T \norm{L^\infty(T)}{|\nabla u|^{r-2}}^{\frac{1}{r-1}} \seminorm{W^{m+1,r'}(T)^d}{u}^{r'},
      \\
      \mathfrak{D}_{\stokes,2}
      &\coloneq
      \truth{k = 0} \sum_{T\in\Ths(u)}\nu_Th_T^2\seminorm{H^1(T)^d}{u}^2,
      \\
      \mathfrak{D}_{\darcy}
      &\coloneq
      \sum_{T \in \Thd(u)} h_T^{2(m+1)} \nu_T \seminorm{H^{m+1}(T)^d}{u}^2;
    \end{aligned}
  \end{equation}

\item Contributions stemming from the velocity-pressure coupling term:
  \begin{equation}\label{eq:mathfrak.C}
    \begin{aligned}
      \mathfrak{C}_{\stokes}
      &\coloneq
      \sum_{T \in \Ths(u)} h_T^{r'(m+1)} \mu_T^{-\frac{1}{r-1}} \seminorm{W^{m+1,r'}(T)}{p}^{r'},
      \\
      \mathfrak{C}_{\darcy}
      &\coloneq
      \sum_{T \in \Thd(u)} h_T^{2(m+1)} \nu_T^{-1} \seminorm{H^{m+2}(T)}{p}^2.
    \end{aligned}
  \end{equation}
\end{itemize}

\begin{theorem}[Regime-dependent velocity error estimate]\label{thm:err.estimate}
  Let Assumption~\ref{ass:regularity} as well as the hypotheses of Lemma~\ref{lem:ah:monotonicity} hold.
  Then,
  \begin{equation}\label{eq:err.estimate}
    \begin{aligned}
      \alpha_\mu \norm{\mu,r,h}{\underline{e}_h}^{q_r}
      + \norm{\nu,h}{\underline{e}_h}^2
      &\lesssim
      \alpha_{\mu}^{-\frac{q_r'}{q_r}} \left(
      \mathfrak{S}_{\stokes,1}^{\frac{q_r'}{r'}}
      + \mathfrak{S}_{\stokes,2}^{\frac{q_r'}r}
      + \mathfrak{D}_{\stokes,1}^{\frac{q_r'}{r'}}
      + \mathfrak{C}_{\stokes}^{\frac{q_r'}{r'}}
      \right)
      \\
      &\quad
      + \mathfrak{S}_{\darcy,1}^{\frac2{r'}}
      + \mathfrak{S}_{\darcy,2}
      + \mathfrak{D}_{\stokes,2}
      + \mathfrak{D}_{\darcy}
      + \mathfrak{C}_{\darcy}.
    \end{aligned}
  \end{equation}
\end{theorem}

\begin{proof}
  See Section~\ref{sec:error.analysis}.
\end{proof}

\begin{table}[h!]\centering
  \subcaptionbox{Contributions in the Stokes-dominated regime}[\textwidth]{%
    \begin{tabular}{ccc}
      \toprule
      Contribution & $1 < r < 2$ & $r \ge 2$ \\
      \midrule
      $\mathfrak{S}_{\stokes,1}^{\frac{q_r'}{r'}}$ & \cellcolor{black!10}{$h^{2(r-1)(m+1)}$} & \cellcolor{black!10}{$h^{r'(m+1)}$} \\
      $\mathfrak{S}_{\stokes,2}^{\frac{q_r'}{r}}$ & 0 & \cellcolor{black!10}{$h^{r'(m+1)}$} \\
      \midrule
      $\mathfrak{D}_{\stokes,1}^{\frac{q_r'}{r'}}$ & $h^{2\left(m+\frac{2}{r'}\right)}$ & $h^{r'\left(m+\frac{2}{r'}\right)}$ \\
      \midrule
      $\mathfrak{D}_{\stokes,2}$ ($k = 0$ only) & $h^2$ & $h^2$ \\
      \midrule
      $\mathfrak{C}_{\stokes}^{\frac{q_r'}{r'}}$ & $h^{2(m+1)}$ & \cellcolor{black!10}{$h^{r'(m+1)}$} \\
      \bottomrule
    \end{tabular}
  }
  \\
  \subcaptionbox{Contributions in the Darcy-dominated regime}[\textwidth]{%
    \begin{tabular}{ccc}
      \toprule
      Contribution & $1 < r < 2$ & $r \ge 2$ \\
      \midrule
      $\mathfrak{S}_{\darcy,1}^{\frac2{r'}}$ & \cellcolor{black!10}{$h^{2(r-1)(m+1)}$} & $h^{\frac{4}{r'} (m+1)}$ \\
      $\mathfrak{S}_{\darcy,2}$ & $h^{2(m+1)}$ & \cellcolor{black!10}{$h^{2(m+1)}$} \\
      \midrule
      $\mathfrak{D}_{\darcy}$ & $h^{2(m+1)}$ & \cellcolor{black!10}{$h^{2(m+1)}$} \\
      \midrule
      $\mathfrak{C}_{\darcy}$ & $ h^{2(m+1)}$ & \cellcolor{black!10}{$h^{2(m+1)}$} \\
      \bottomrule
    \end{tabular}
  }
  \caption{Convergence rate of the dominant terms for the error components defined by~\eqref{eq:mathfrak.S},~\eqref{eq:mathfrak.D}, and~\eqref{eq:mathfrak.C}. $q_r$ is defined in \eqref{eq:def.qr} and the dominant terms in each column are highlighted in gray.
    \label{tab:conv.rate.dominant}}
\end{table}

\begin{table}[h!]\centering
  \begin{tabular}{ccccc}
    \toprule
    \multirow{2}{*}{Global regime}
    & \multicolumn{2}{c}{$\norm{\mu,r,h}{\underline{e}_h}$}
    & \multicolumn{2}{c}{$\norm{\nu,h}{\underline{e}_h}$} \\
    {} & $1 < r < 2$ & $r \ge 2$ & $1 < r < 2$ & $r \ge 2$ \\
    \midrule
    Stokes-dominated
    & $h^{(r-1)(m+1)}$ & $h^{\frac{m+1}{r-1}}$ & $h^{(r-1)(m+1)}$ & $h^{\frac{r'}{2}(m+1)}$ \\
    Darcy-dominated
    & $h^{(r-1)(m+1)}$ & $h^{\frac2r (m+1)}$ & $h^{(r-1)(m+1)}$ & $h^{m+1}$ \\
    \bottomrule
  \end{tabular}
  \caption{Convergence rate of the error components in the left-hand side of~\eqref{eq:err.estimate} as a function of the regime and of the power-law exponent $r$.
    Notice that the approximation properties of the space are matched for $r = 2$ and also for the $\norm{\nu,h}{\cdot}$-norm of the error in the Darcy-dominated regime when $r \ge 2$.
    \label{tab:oc}}
\end{table}

\begin{remark}[Convergence rates]\label{rem:conv.rates}
  The convergence rates of the dominant terms in each contribution raised to the power with which it appears in~\eqref{eq:err.estimate} are summarized in Table~\ref{tab:conv.rate.dominant}.
  Table~\ref{tab:oc} shows how the estimate~\eqref{eq:err.estimate} simplifies for Stokes- or Darcy-dominated flows, i.e., flows for which either $\CfT < 1$ for all $T \in \Th$ or $\CfT \ge 1$ for all $T \in \Th$.
    As noticed in Remark~\ref{rem:regimes}, the Darcy regime is pre-asymptotic and can only be observed on sufficiently coarse meshes.
    
  In the pure Stokes case corresponding to $\nu \equiv 0$, we have $\Thd(u)=\emptyset$ and $\CfT=0$ for all $T\in\Th$, so all the terms except $\mathfrak{S}_{\stokes,1}$, $\mathfrak{S}_{\stokes,2}$ and $\mathfrak{C}_\stokes$ vanish, and we recover the estimates of \cite{Botti.Castanon-Quiroz.ea:21}.
  In the linear case corresponding to $r = 2$, the above estimates are additionally consistent with the ones proved in \cite{Di-Pietro.Droniou:23} for a different HHO scheme.
  For $r < 2$, it is most likely possible to improve the orders of convergence in the Stokes-dominated regime adapting the ideas of~\cite{Di-Pietro.Droniou.ea:21}.
  In this reference, which focuses on scalar Leray--Lions (scalar) problems, it is shown that the approximation properties of the HHO space in the discrete $W^{1,r}$-norm can be matched for non-degenerated fluxes.
  The numerical experiments in~\cite[Section~6.4]{Di-Pietro.Droniou.ea:21} (see also~\cite[Section~6.2]{Beirao-da-Veiga.Di-Pietro.ea:24} concerning a discontinuous Galerkin method) seem, on the other hand, to confirm that the original estimates of~\cite{Di-Pietro.Droniou:17*1} are recovered in practice in the degenerated case.
  To avoid further increasing the technicality of this paper, we have preferred to stick to estimates inspired by~\cite{Di-Pietro.Droniou:17*1} in the Stokes-dominated regime.

  An extensive comparison between these estimates and the ones obtained for conforming methods~\cite{Barrett.Liu:94} in the framework of scalar $p$-diffusion is provided in~\cite[Remark~3.3]{Di-Pietro.Droniou:17*1}.
  We also notice, for the sake of completeness, that improved estimates for lowest-order discontinuous Galerkin methods have been obtained using techniques based on the Scott--Zhang interpolator in~\cite{Diening.Koner.ea:14}; see~\cite{Belenki.Berselli.ea:12} for similar estimates for the $p$-Stokes problem and the recent work~\cite{Diening.Hirn.ea:25}, which additionally addresses pressure-roubstness.
  In the latter contribution, a key ingredient is a divergence-preserving smoothing operator.
  It is not clear whether the techniques used in the above works apply to the present setting, as they heavily rely on the use of standard meshes (e.g., for the construction of the interpolation and smoothing operators), whereas here we consider general polytopal meshes.
  %% They also specifically target the lowest-order case, where the main technical problem lies in the interpolation of functions in negative Sobolev spaces.
\end{remark}

\begin{remark}[Robustness in the Darcy limit]\label{rem:darcy.limit}
  Using arguments similar to those in \cite[Remark~13]{Botti.Di-Pietro.ea:18}, it can be checked that the
  proposed method is also applicable in the singular limit corresponding to the Darcy case ($\mu \equiv 0$ and $\min_{\Omega}\nu>0$), which is obtained dropping the first term in~\eqref{eq:strong:momentum} and modifying~\eqref{eq:strong:bc} by only enforcing $u \cdot n = 0$ (with $n$ denoting the outward normal vector to $\partial \Omega$).

  We note that the estimate \eqref{eq:err.estimate} is robust in this Darcy limit: since we have $\CfT=+\infty$ for all $T\in\Th$, the only remaining contributions to the error bound come from $\mathfrak{D}_\darcy$ and $\mathfrak{C}_\darcy$ (all the other terms correspond to sums over $\Ths(u)=\emptyset$ or are scaled by $\CfT^{-1}=0$). The resulting scheme has orders of convergence consistent with the ones found in~\cite{Di-Pietro.Ern:17} for the Mixed High-Order method, which corresponds to a variation with element velocities in $\nabla \Poly{k}(T)$ instead of $\Nedelec{k}(T)$.
\end{remark}

\begin{remark}[Pressure error estimate]
  Denoting by $\epsilon_h \coloneq p_h - \lproj{k}{h} p$ the error on the pressure, we introduce the following norm:
  \[
    \norm{\rm P}{\epsilon_h}\coloneq \sup_{\underline{v}_h \in \VhZ \setminus\{\underline{0}\}}\frac{b_h(\underline{v}_h,\epsilon_h)}{\norm{\mu,r,h}{\underline{v}_h}+\norm{\nu,h}{\underline{v}_h}}.
  \]
  Using the result of Theorem~\ref{thm:err.estimate}, we can derive the following estimate:
  \[
  \norm{\rm P}{\epsilon_h}\lesssim
  \mathfrak{S}_{\stokes,1}^{\frac1{r'}}
  + \mathfrak{S}_{\stokes,2}^{\frac1r}
  + \mathfrak{D}_{\stokes,1}^{\frac1{r'}}
  + \mathfrak{C}_{\stokes}^{\frac1{r'}}
  + \mathfrak{S}_{\darcy,1}^{\frac1{r'}}
  + \mathfrak{S}_{\darcy,2}^{\frac12}
  +\mathfrak{D}_{\stokes,2}^{\frac12}
  + \mathfrak{D}_{\darcy}^{\frac12}
  + \mathfrak{C}_{\darcy}^{\frac12}.
  \]
  The resulting convergence rate only depends on $r$, not on the global regime: it is $\mathcal O(h^{(r-1)(m+1)})$ if $1<r<2$, and $\mathcal O(h^{m+1})$ if $r\ge 2$.
\end{remark}

%------------------------------------------------------------------------------%

\section{Error analysis}\label{sec:error.analysis}

In this section, we prove Lemma~\ref{lem:ah:monotonicity} and establish regime-dependent consistency estimates for the contributions to the error stemming from the Stokes, Darcy, and velocity-pressure coupling terms that will be needed in the proof of Theorem~\ref{thm:err.estimate}.
We will often use the fact that $\sigma_{|T}$ is simply $\mu_T$ times the (tensorised) $r$-Laplace operator when we invoke results for the $r$-Laplace equation from \cite[Chapter~6]{Di-Pietro.Droniou:20}.

\subsection{Monotonicity of the Stokes bilinear form}\label{sec:ah:monotonicity}

\begin{proof}[Proof of Lemma~\ref{lem:ah:monotonicity}]
  If $r\ge 2$, the bound \eqref{eq:ah:monotonicity} directly follows by adding together the local equations \cite[Eqs.~(6.66) and (6.67)]{Di-Pietro.Droniou:20} (which are established in the case $\mu\equiv 1$), multiplying the resulting equation by $\mu_T$,
  using the norm equivalence of \cite[Lemma 6.11]{Di-Pietro.Droniou:20} (with minor adjustments to account for the element unknown in the Nédélec space), and summing over $T$.

  Let us therefore consider  the case $1<r<2$. Adapting the reasoning in \cite[Proof of Theorem 6.19, Part (ii.B)]{Di-Pietro.Droniou:20} to include the parameter $\mu$, we obtain
  \begin{equation}\label{eq:monotony.ash.r<2}
    \norm{\mu,r,h}{\underline{e}_h}^r
    \lesssim
    \left(
    a_{\mu,h}(\underline{u}_h, \underline{e}_h)
    - a_{\mu,h}(\Ih u, \underline{e}_h)
    \right)^{\frac{r}{2}}
    \left(
    \norm{\mu,r,h}{\underline{u}_h}^r+\norm{\mu,r,h}{\Ih u}^r
    \right)^{1-\frac{r}{2}}.
  \end{equation}
  We now bound the terms in the last bracket.
  Make $\underline{v}_h=\underline{u}_h$ in \eqref{eq:discrete:momentum}, $q_h=p_h$ in \eqref{eq:discrete:mass} and add up the two equations to get, since $a_{\nu,h}(\underline{u}_h,\underline{u}_h)\ge 0$,
  \[
  a_{\mu,h}(\underline{u}_h,\underline{u}_h)\le \norm{L^{r'}(\Omega)^d}{f}
  \norm{L^r(\Omega)^d}{\PDh\underline{u}_h}.
  \]
  For each $T\in\Th$, the definitions~\eqref{eq:a.mu.h.T} of $a_{\mu,T}$ and~\eqref{eq:sigma} of $\sigma$ together with the norm equivalence \cite[Lemma 6.11]{Di-Pietro.Droniou:20} as above yield
  \[
  a_{\mu,T}(\underline{u}_h,\underline{u}_h)
  = \mu_T\left(\norm{L^r(T)^{d\times d}}{\GT\underline{u}_T}^r
  + h_T^{1-r}\norm{L^r(\partial T)^d}{\Delta_{\partial T}\underline{u}_T}^r\right)\simeq \norm{\mu,r,T}{\underline{u}_h}^r.
  \]
  The local Poincar\'e--Wirtinger inequality \eqref{eq:poincare.PDT} with $r=q$, the global Poincar\'e inequality \cite[Proposition 5.4]{Di-Pietro.Droniou:17}, and the fact that $\underline{\mu}^{-1}\mu_T\ge 1$ for all $T\in\Th$ yield $\norm{L^r(\Omega)^d}{\PDh\underline{u}_h}\lesssim \underline{\mu}^{-\frac1r}\norm{\mu,r,T}{\underline{u}_h}$.
  Plugging this bound into the equations above and simplifying gives
  \begin{equation}\label{eq:bound.uh}
    \norm{\mu,r,T}{\underline{u}_h}^{r-1}\lesssim \underline{\mu}^{-\frac1r}.
  \end{equation}
  To estimate $\norm{\mu,r,h}{\Ih u}$, we start with the boundedness of the local interpolator \cite[Proposition 6.24]{Di-Pietro.Droniou:20}
  (with $p=r$), which we raise to the power $r$, multiply by $\mu_T$, and sum over $T$ to get
  \[
  \norm{\mu,r,h}{\Ih u}^r\lesssim \int_\Omega \mu|\nabla u|^r.
  \]
  Using $(u,p)$ as test functions in the weak formulation of \eqref{eq:strong} and applying a Poincar\'e inequality, we infer
  $\norm{\mu,r,h}{\Ih u}^{r-1}\lesssim \underline{\mu}^{-\frac1r}$, as in the discrete case.
  Plugging this equation together with \eqref{eq:bound.uh} into \eqref{eq:monotony.ash.r<2} and taking the power $\frac2r$ of the resulting inequality
  concludes the proof.
\end{proof}

\subsection{Consistency of the Stokes term}

\begin{lemma}[Regime-dependent consistency estimate for the Stokes term]\label{lem:Err.stokes.h}
  Let Assumption~\ref{ass:regularity} hold and define the consistency error linear form $\Err_{\stokes,h} : \VhZ \to \Real$ for $a_{\mu,h}$ such that, for all $\underline{v}_h \in \VhZ$,
  \begin{equation}\label{eq:Err.stokes.h}
    \Err_{\stokes,h}(\underline{v}_h)
    \coloneq
    -\int_\Omega \nabla \cdot \sigma(\nabla u) \cdot \PDh \underline{v}_h
    - a_{\mu,h}(\Ih u, \underline{v}_h).
  \end{equation}
  Then, recalling \eqref{eq:mathfrak.S}, it holds, for all $\underline{v}_h \in \VhZ$,
  \begin{equation}\label{eq:Err.stokes.h:estimate}
    \Err_{\stokes,h}(\underline{v}_h)
    \lesssim \left(
    \mathfrak{S}_{\stokes,1}^{\frac1{r'}}
    + \mathfrak{S}_{\stokes,2}^{\frac1r}
    \right) \norm{\mu,r,h}{\underline{v}_h}
    + \left(
    \mathfrak{S}_{\darcy,1}^{\frac1{r'}}
    + \mathfrak{S}_{\darcy,2}^{\frac12}
    \right) \norm{\nu,h}{\underline{v}_h}.
  \end{equation}
\end{lemma}

\begin{proof}
  We start by writing
  \[
  \begin{aligned}
    &-\int_T \nabla \cdot \sigma(\nabla u) \cdot \PDT \underline{v}_T
    \\
    &\begin{aligned}[t]
       &= \int_T \nabla \cdot \sigma(\nabla u) \cdot (v_T-\PDT \underline{v}_T)
       - \int_T \nabla \cdot \sigma(\nabla u) \cdot v_T
       \\
       \overset{\eqref{eq:pi.k-1.PDT=vT}}&=\int_T \left[
         \nabla \cdot \sigma(\nabla u)
         - \lproj{k-1}{T} (\nabla \cdot \sigma(\nabla u))
         \right] \cdot (v_T-\PDT \underline{v}_T) - \int_T \nabla \cdot \sigma(\nabla u) \cdot v_T.
     \end{aligned}
  \end{aligned}
  \]
  Substituting in \eqref{eq:Err.stokes.h} and using standard arguments (see, e.g., Point (i) in the proof of \cite[Theorem~6.19]{Di-Pietro.Droniou:20}), we get the following decomposition of the error:
  \begin{equation}\label{eq:Err.stokes.h:basic}
    \begin{aligned}
      \Err_{\stokes,h}(\underline{v}_h)
      &=
      \sum_{T \in \Th} \int_T \left[
        \nabla \cdot \sigma(\nabla u) - \lproj{k-1}{T} (\nabla \cdot \sigma(\nabla u))
        \right] \cdot (v_T - \PDT \underline{v}_T)
      \\
      &\quad
      +  \sum_{T \in \Th} \int_T \left[
        \sigma(\nabla u) - \sigma(\lproj{k}{T} \nabla u)
        \right] : \GT \underline{v}_T
      - s_{\mu,h}(\huline{u}_h, \underline{v}_h)
      \\
      &\quad
      + \sum_{T \in \Th} \int_{\partial T} \left[
        \sigma(\nabla u)_{|T} - \lproj{k}{T} \sigma(\nabla u)
        \right] n_{\partial T} \cdot (v_{\partial T} - v_T)
      \\
      &\eqcolon \term_1 + \term_2 + \term_3 + \term_4.
    \end{aligned}
  \end{equation}
  We write $\term_i = \sum_{T \in \Th} \term_i(T)$ and estimate the contributions for an element $T \in \Th$ in each local regime.%
  \medskip\\
  (i) \emph{$\CfT < 1$ (Stokes-dominated regime).}
  To estimate the first term, we use a H\"older inequality and write
  \[
  \begin{aligned}
    \term_1(T)
    &\le
    \norm{L^{r'}(T)^d}{\nabla \cdot \sigma(\nabla u)
      - \lproj{k-1}{T} (\nabla \cdot \sigma(\nabla u))} ~
    \norm{L^r(T)^d}{v_T - \PDT \underline{v}_T}
    \\
    \overset{\eqref{eq:poincare.PDT}}&\lesssim
    h_T^{m+1} \mu_T^{-\frac1r}\seminorm{W^{m+1,r'}(T)^{d\times d}}{\sigma(\nabla u)} ~
    \mu_T^{\frac1r}\norm{1,r,T}{\underline{v}_T},
  \end{aligned}
  \]
  where we have additionally used the approximation properties of $\lproj{k-1}{T}$ in the conclusion if $k\ge 1$ and we have simply noticed that $\lproj{-1}{T} (\nabla \cdot \sigma(\nabla u))=0$ if $k=0$ (in which case $m=0$).

  Let us now turn to the second term.
  We start with a H\"older inequality and \eqref{eq:norm.equiv.GT.1qT} with $q=r$:
  \begin{equation}\label{eq:est.T2.stokes}
    \term_2(T)
    \le \norm{L^{r'}(T)^{d\times d}}{\sigma(\nabla u)-\sigma(\lproj{k}{T}\nabla u)}
    ~ \norm{1,r,T}{\underline{v}_T}.
  \end{equation}
  If $1<r<2$, using $|\sigma(\nabla u)-\sigma(\lproj{k}{T}\nabla u)|\le \mu_T2^{2-r}|\nabla u-\lproj{k}{T}\nabla u|^{r-1}$ (see \cite[Eq.~(6.56)]{Di-Pietro.Droniou:20}), we have
  \[
  \term_2(T)
  \lesssim \mu_T\norm{L^{r}(T)^{d\times d}}{\nabla u-\lproj{k}{T}\nabla u}^{r-1}
  ~ \norm{1,r,T}{\underline{v}_T}.
  \]
  Invoking the approximation properties of $\lproj{k}{T}$, we arrive at
  \[
  \term_2(T)
  \lesssim h_T^{(r-1)(m+1)}\mu_T^{\frac{1}{r'}}\seminorm{W^{m+2,r}(T)^{d}}{u}^{r-1}
  ~ \mu_T^{\frac1r}\norm{1,r,T}{\underline{v}_T}
  \quad \text{if $1 < r < 2$.}
  \]
  If $r\ge 2$, we apply \cite[Eq.~(6.60)]{Di-Pietro.Droniou:20} to \eqref{eq:est.T2.stokes} and get
  \[
  \term_2(T)\le
  h_T^{m+1} \mu_T \seminorm{W^{m+2,r}(T)^{d}}{u} \norm{L^r(T)^{d\times d}}{\nabla u}^{r-2}
  ~ \norm{1,r,T}{\underline{v}_T}.
  \]
  We then notice, still using $r\ge 2$, that
  \[
  \norm{L^r(T)^{d\times d}}{\nabla u}^{r-2}\le |T|^{\frac{r-2}{r}} \norm{L^\infty(T)}{|\nabla u|^{r-2}}
  \lesssim |T|^{\frac{r-2}{r}}\mu_T^{-1}\kappaT.
  \]
  Additionally multiplying by $1 = \mu_T^{-\frac1r}\,\mu_T^{\frac1r}$, this leads to
  \begin{equation}\label{eq:Err.stokes:T2:stokes}
    \term_2(T)\le |T|^{\frac{r-2}{r}} ~
    h_T^{m+1} \mu_T^{-\frac1r}\kappaT\seminorm{W^{m+2,r}(T)^{d}}{u}
    ~ \mu_T^{\frac1r} \norm{1,r,T}{\underline{v}_T}
    \quad \text{if $r \ge 2$.}
  \end{equation}

  Concerning the remaining terms, known estimates for HHO methods for the $r$-Laplace equation (see \cite[Theorem~6.19]{Di-Pietro.Droniou:20}) together with \eqref{eq:norm.equiv.GT.1qT} with $q=r$ and $\mu_T=\mu_T^{\frac{1}{r'}}\mu_T^{\frac1r}$ give, for any element $T \in \Ths(u)$ in the Stokes-dominated regime,
  \[
  \begin{aligned}
    \term_3(T)
    &\lesssim
    h_T^{(r-1)(m+1)}  \mu_T^{\frac1{r'}} \seminorm{W^{m+2,r}(T)^d}{u}^{r-1}
    ~ \mu_T^{\frac1r} \norm{1,r,T}{\underline{v}_T},
    \\
    \term_4(T)
    &\lesssim h_T^{m+1} \mu_T^{-\frac1r} \seminorm{W^{m+1,r'}(T)^{d\times d}}{\sigma(\nabla u)}
    ~ \mu_T^{\frac1r}\norm{1,r,T}{\underline{v}_T}.
  \end{aligned}
  \]
  \\
  (ii) \emph{$1\le \CfT < +\infty$ (Darcy-dominated regime).}
  Let us now consider an element $T \in \Thd(u)$ in the Darcy-dominated regime.
  For the first term, we start with a Cauchy--Schwarz inequality to write
  \[
  \begin{aligned}
    \term_1(T)
    &\lesssim \norm{L^2(T)^d}{\nabla \cdot \sigma(\nabla u) - \lproj{k-1}{T}(\nabla \cdot \sigma(\nabla u))} \, \norm{L^2(T)^d}{\PDT \underline{v}_T -v_T}
    \\
    &\lesssim h_T^m \seminorm{H^{m+1}(T)^{d\times d}}{\sigma(\nabla u)} \, \norm{\darcy,T}{\underline{v}_T},
  \end{aligned}
  \]
  where we have used the approximation properties of $\lproj{k-1}{T}$ for the first factor if $k\ge 1$ or $\lproj{-1}{T}(\nabla \cdot \sigma(\nabla u))=0$ if $k = m = 0$,
  and~\eqref{eq:poincare.PDT} with $q = 2$ followed by~\eqref{eq:norm.1.2.T<=hT-1.norm.0.2.T} and~\eqref{eq:L2.norm.darcy} for the second factor.
  We then write
  \begin{equation}\label{eq:stokes.CfT>=1.law}
    1
    = \nu_T^{-\frac12}\nu_T^{\frac12}
    \overset{\eqref{eq:CfT}}= h_T\kappaT^{-\frac12}\CfT^{-\frac12}\nu_T^{\frac12},
  \end{equation}
  and deduce
  \[
  \term_1(T)
  \lesssim
  h_T^{m+1} \kappaT^{-\frac12} \CfT^{-\frac12}\seminorm{H^{m+1}(T)^{d\times d}}{\sigma(\nabla u)}
  ~ \nu_T^{\frac12} \norm{\darcy,T}{\underline{v}_T}.
  \]
  For the second and third term, recalling the definition \eqref{eq:jh.jT} of $s_{\mu,T}$ and using Cauchy--Schwarz inequalities, we get
  \begin{multline*}
    \term_2(T) + \term_3(T)
    \\
    \le
    h_T^{-1} \left(
    \norm{L^2(T)^{d\times d}}{\sigma(\nabla u) - \sigma(\lproj{k}{T}\nabla u)}
    + h_T^{\frac32-r} \norm{L^2(\partial T)^d}{\sigma(\Delta_{\partial T}^k \huline{u}_T)}
    \right)
    \\
    \times h_T \left(
    \norm{L^2(T)^{d\times d}}{\GT \underline{v}_T}
    + h_T^{-\frac12} \norm{L^2(\partial T)^d}{\Delta_{\partial T}^k \underline{v}_T}
    \right).
  \end{multline*}
  By \eqref{eq:norm.equiv.GT.1qT} with $q=2$, the last factor is $\lesssim h_T \norm{1,2,T}{\underline{v}_T} \overset{\eqref{eq:norm.1.2.T<=hT-1.norm.0.2.T}}\lesssim \norm{0,2,T}{\underline{v}_T} \overset{\eqref{eq:L2.norm.darcy}}\lesssim \norm{\darcy,T}{\underline{v}_T}$.
  Accounting for this fact and using \eqref{eq:stokes.CfT>=1.law}, we get
  \begin{multline}\label{eq:stokes:T1(T)+T2(T):basic}
    \term_2(T) + \term_3(T)
    \lesssim
    \kappaT^{-\frac12} \CfT^{-\frac12}
    \\
    \times
    \left(
    \norm{L^2(T)^{d\times d}}{\sigma(\nabla u) - \sigma(\lproj{k}{T}\nabla u)}
    + h_T^{\frac32-r} \norm{L^2(\partial T)^d}{\sigma(\Delta_{\partial T}^k \huline{u}_T)}
    \right)
    \nu_T^{\frac12} \norm{\darcy,T}{\underline{v}_T}.
  \end{multline}
  We next estimate the terms in parentheses.
  If $r \ge 2$, by \cite[Lemma~6.28]{Di-Pietro.Droniou:20} it holds
  \[
  \begin{aligned}
    |\sigma(\nabla u) - \sigma(\lproj{k}{T} \nabla u)|
    &\lesssim \mu_T \left(|\nabla u|^{r-2} + |\lproj{k}{T} \nabla u|^{r-2}\right) |\nabla u - \lproj{k}{T} \nabla u|
    \\
    &\lesssim \kappaT|\nabla u - \lproj{k}{T} \nabla u|,
  \end{aligned}
  \]
  where the conclusion follows by noticing that, on $T$,
  $|\nabla u|^{r-2} \le \norm{L^\infty(T)}{|\nabla u|^{r-2}}$ and
  \[
  \begin{aligned}
    |\lproj{k}{T} \nabla u|^{r-2}
    &\le \norm{L^\infty(T)}{|\lproj{k}{T} \nabla u|^{r-2}}
    = \norm{L^\infty(T)^{d\times d}}{\lproj{k}{T} \nabla u}^{r-2}
    \\
    &\lesssim \norm{L^\infty(T)^{d\times d}}{\nabla u}^{r-2}
    = \norm{L^\infty(T)}{|\nabla u|^{r-2}},
  \end{aligned}
  \]
  where we have used the fact that $x \mapsto x^{r-2}$ is monotonically increasing for $r \ge 2$ in the equalities and
  the $L^\infty(T)$-boundedness of $\lproj{k}{T}$ in the second inequality.
  We infer
  \[
  \begin{aligned}
    \norm{L^2(T)^{d\times d}}{\sigma(\nabla u) - \sigma(\lproj{k}{T}\nabla u)}
    &\lesssim
    \kappaT \norm{L^2(\Omega)^{d\times d}}{\nabla u - \lproj{k}{T} \nabla u}
    \\
    &\lesssim
    h_T^{m+1} \kappaT \seminorm{H^{m+2}(T)^{d}}{u},
  \end{aligned}
  \]
  where the conclusion follows from the approximation properties of $\lproj{k}{T}$.
  If $1 < r < 2$, on the other hand, \cite[Lemma~6.28]{Di-Pietro.Droniou:20} gives
  $|\sigma(\nabla u) - \sigma(\lproj{k}{T}\nabla u)| \lesssim \mu_T |\nabla u - \lproj{k}{T}\nabla u|^{r-1}$, and thus
  \[
  \begin{aligned}
    \norm{L^2(T)^{d\times d}}{\sigma(\nabla u) - \sigma(\lproj{k}{T}\nabla u)}
    &\le |T|^{\frac12-\frac{1}{r'}} ~
    \norm{L^{r'}(T)^{d\times d}}{\sigma(\nabla u) - \sigma(\lproj{k}{T}\nabla u)}\\
    &\lesssim |T|^{\frac12-\frac{1}{r'}} ~
    \mu_T \norm{L^r(T)^{d\times d}}{\nabla u - \lproj{k}{T}\nabla u}^{r-1}\\
    &\lesssim |T|^{\frac12-\frac{1}{r'}} ~
    h_T^{(r-1)(m+1)} \mu_T \seminorm{W^{m+2,r}(T)^d}{u}^{r-1},
  \end{aligned}
  \]
  where we have used, in the first line, a H\"older inequality with exponent $\frac{r'}{2} > 1$ and, in the conclusion, the approximation properties of $\lproj{k}{T}$.
  Gathering the above estimates, we thus have
  \begin{multline}\label{eq:stokes:est.sigma}
    \norm{L^2(T)^{d\times d}}{\sigma(\nabla u) - \sigma(\lproj{k}{T}\nabla u)}
    \\
    \lesssim
    \begin{cases}
      |T|^{\frac12-\frac{1}{r'}} ~ h_T^{(r-1)(m+1)} \mu_T \seminorm{W^{m+2,r}(T)^d}{u}^{r-1}
      & \text{if $1 < r < 2$},
      \\
      h_T^{m+1} \kappaT \seminorm{H^{m+2}(T)^d}{u}
      & \text{if $r \ge 2$}.
    \end{cases}
  \end{multline}

  The estimate of the second term in parentheses in~\eqref{eq:stokes:T1(T)+T2(T):basic} is similar.
  Consider first the case $r \ge 2$.
  Recalling the definition \eqref{eq:sigma} of $\sigma$ and using the fact that $|\Delta_{\partial T}^k \huline{u}_T|^{r-2} \le \norm{L^\infty(\partial T)^d}{\Delta_{\partial T}^k \huline{u}_T}^{r-2}$ on $\partial T$ since $r \ge 2$, we can write
  \[
  \begin{aligned}
    \norm{L^2(\partial T)^d}{\sigma(\Delta_{\partial T}^k \huline{u}_T)}
    &\le \mu_T \norm{L^\infty(\partial T)^d}{\Delta_{\partial T}^k \huline{u}_T}^{r-2}
    \norm{L^2(\partial T)^d}{\Delta_{\partial T}^k \huline{u}_T}
    \\
    \overset{\eqref{eq:Delta.k:consistency},\,\eqref{eq:Delta.k:boundedness}}&\lesssim
    h_T^{m+r-\frac12} \kappaT \seminorm{H^{m+2}(T)^d}{u},
  \end{aligned}
  \]
  where, for the second factor, we have additionally used the fact that $\norm{L^\infty(T)^{d\times d}}{\nabla u}^{r-2} = \norm{L^\infty(T)}{|\nabla u|^{r-2}}$ since $r\ge 2$ together with the definition~\eqref{eq:CfT} of $\kappaT$ in the last inequality.
  If, on the other hand, $1 < r < 2$,
  \begin{align*}
    \norm{L^2(\partial T)^d}{\sigma(\Delta_{\partial T}^k \huline{u}_T)}
    &\le |\partial T|^{\frac12-\frac{1}{r'}} ~ \norm{L^{r'}(\partial T)^d}{\sigma(\Delta_{\partial T}^k \huline{u}_T)}\\
    \overset{\eqref{eq:sigma}}&=
    |\partial T|^{\frac12-\frac{1}{r'}} ~ \mu_T \norm{L^{r}(\partial T)^{d}}{\Delta_{\partial T}^k \huline{u}_T}^{r-1}\\
    \overset{\eqref{eq:Delta.k:consistency}}&\lesssim
    |\partial T|^{\frac12-\frac{1}{r'}} ~ h_T^{(r-1)(m+2)-\frac{1}{r'}} \mu_T \seminorm{W^{m+2,r}(T)^d}{u}^{r-1}\\
    &\lesssim |T|^{\frac12-\frac{1}{r'}} ~ h_T^{(r-1)(m+2)-\frac{1}{2}} \mu_T \seminorm{W^{m+2,r}(T)^d}{u}^{r-1},
  \end{align*}
  where the conclusion follows by mesh regularity which allows us to write $|\partial T|^{\frac12-\frac{1}{r'}}\lesssim h_T^{\frac{1}{r'}-\frac12} | T|^{\frac12-\frac{1}{r'}}$.
  In conclusion,
  \begin{multline}\label{eq:stokes:est.Delta.k}
    h_T^{\frac32-r} \norm{L^2(\partial T)}{\sigma(\Delta_{\partial T}^k \huline{u}_T)}
    \\
    \lesssim
    \begin{cases}
      |T|^{\frac12-\frac{1}{r'}} ~ h_T^{(r-1)(m+1)}\mu_T \seminorm{W^{m+2,r}(T)^d}{u}^{r-1}
      & \text{if $1 < r < 2$},
      \\
      h_T^{m+1} \kappaT \seminorm{H^{m+2}(T)^d}{u}
      & \text{if $r \ge 2$}.
    \end{cases}
  \end{multline}
  Plugging \eqref{eq:stokes:est.sigma} and \eqref{eq:stokes:est.Delta.k} into \eqref{eq:stokes:T1(T)+T2(T):basic}, we thus arrive at
  \begin{multline}\label{eq:Err.stokes:T2+T3:darcy}
    \term_2(T) + \term_3(T)
    \\
    \lesssim
    \begin{cases}
      \begin{aligned}
        |T|^{\frac{r'-2}{2r'}} ~
        h_T^{(r-1)(m+1)} \CfT^{-\frac12} \mu_T \kappaT^{-\frac12} \seminorm{W^{m+2,r}(T)^d}{u}^{r-1}
        \\
        \times \nu_T^{\frac12} \norm{\darcy,T}{\underline{v}_T}
      \end{aligned}
      & \text{if $1 < r < 2$},
      \\
      h_T^{m+1} \CfT^{-\frac12}
      \kappaT^{\frac12} \seminorm{H^{m+2}(T)^d}{u}
      ~ \nu_T^{\frac12} \norm{\darcy,T}{\underline{v}_T}.
      & \text{if $r \ge 2$}.
    \end{cases}
  \end{multline}

  For the fourth term, we start with a $(2,\infty,2)$-H\"{o}lder inequality together with $\norm{L^\infty(\partial T)^d}{n_{\partial T}}\le 1$ to write
  \[
  \begin{aligned}
    \term_4(T)
    &\le
    \norm{L^2(\partial T)^{d\times d}}{\sigma(\nabla u) - \lproj{k}{T} \sigma(\nabla u)}
    ~ \norm{L^2(\partial T)^d}{v_{\partial T} - v_T}
    \\
    \overset{\eqref{eq:stokes.CfT>=1.law}}&=
    h_T^{\frac12} \kappaT^{-\frac12} \CfT^{-\frac12}
    \norm{L^2(\partial T)^{d\times d}}{\sigma(\nabla u) - \lproj{k}{T} \sigma(\nabla u)}
    \\
    &\quad \times
    \nu_T^{\frac12} h_T^{\frac12} \norm{L^2(\partial T)^d}{v_{\partial T} - v_T}
    \\
    &\lesssim
    h_T^{m+1} \kappaT^{-\frac12} \CfT^{-\frac12}
    \seminorm{H^{m+1}(T)^{d\times d}}{\sigma(\nabla u)}
    ~ \nu_T^{\frac12} \norm{\darcy,T}{\underline{v}_T},
  \end{aligned}
  \]
  where we have used the approximation properties of $\lproj{k}{T}$ and the definition \eqref{eq:norm.0.2.T} of the $\norm{0,2,T}{{\cdot}}$-norm followed by~\eqref{eq:L2.norm.darcy} to conclude.
  \medskip\\
  (iii) \emph{$\CfT = +\infty$ (pure Darcy).}
  The pure Darcy case corresponds to $\mu_T \norm{L^\infty(T)}{|\nabla u|^{r-2}} = 0$, i.e., $\mu_T = 0$ if $1 < r < 2$ and either $\mu_T = 0$ or $\norm{L^\infty(T)}{|\nabla u|^{r-2}} = 0$ for $r \ge 2$.
  Let us show that $\term_i(T) = 0$ for $i \in \{1,\ldots,4\}$ in this case.
  If $r = 2$, the condition $\mu_T \norm{L^\infty(T)}{|\nabla u|^{r-2}} = 0$ reduces to $\mu_T = 0$, as for $1 < r < 2$, which readily implies the conclusion since all the terms $\term_i(T)$ are multiplied by $\mu_T$.
  A similar conclusion can be drawn when $\mu_T = 0$ and $r > 2$.
  Let us now consider the case $r > 2$ and $\norm{L^\infty(T)}{|\nabla u|^{r-2}} = 0$.
  The latter condition gives $\nabla u \equiv 0$ in $T$, showing that $\term_1(T) = \term_2(T) = \term_4(T) = 0$.
  On the other hand $\nabla u \equiv 0$ also implies $\Delta_{\partial T}^k \huline{u}_T = 0$ since $u$ is constant in $T$, justifying $\term_3(T) = 0$.
  \medskip\\
  (iv) \emph{Conclusion.}
    To conclude, we split the sums over the elements in \eqref{eq:Err.stokes.h:basic} writing $\sum_{T \in \Th} \term_i(T) = \sum_{T \in \Ths(u)} \term_i(T) + \sum_{T \in \Thd(u)} \term_i(T)$, use the local bounds derived above, and invoke Cauchy--Schwarz and H\"{o}lder inequalities as needed.

    Specifically, for the sum on Stokes-dominated elements, we use $(r',r)$-H\"older inequalities for the terms $\term_1(T)$, $\term_3(T)$, $\term_4(T)$, and $\term_2(T)$ if $1 < r < 2$.
    If $r \ge 2$, we use a $\left(\frac{2}{r-2}, r, r\right)$-H\"{o}lder inequality to write
    \begin{multline*}
      \sum_{T \in \Ths(u)} \term_2(T)
      \\
      \overset{\eqref{eq:Err.stokes:T2:stokes}}\lesssim
      \left(
      \sum_{T \in \Ths(u)} \hspace{-1ex} |T|
      \right)^{\frac{r-2}{r}} \hspace{-1ex} \left(
      \sum_{T \in \Ths(u)} \hspace{-1ex} h_T^{r(m+1)} \mu_T^{-1} \kappaT^r \seminorm{W^{m+2,r}(T)^d}{u}^r
      \right)^{\frac1r}\norm{1,\mu,h}{\underline{v}_h},
    \end{multline*}
    further noticing that $\sum_{T \in \Ths(u)} |T| \le |\Omega| \lesssim 1$.

    Concerning the sum on Darcy-dominated elements, we apply Cauchy--Schwarz inequalities to $\term_1(T)$ and $\term_4(T)$ and also for $\term_2(T) + \term_3(T)$ for $r \ge 2$.
    When $1 < r < 2$, for the latter terms we instead use a $\left(\frac{2 r'}{r'-2}, r',2\right)$-H\"{o}lder inequality, further recalling that $r' \overset{\eqref{eq:r'}}= \frac{r}{r-1}$, along with $\sum_{T \in \Thd(u)} |T| \le |\Omega| \lesssim 1$ to write
    \begin{multline*}
    \sum_{T \in \Thd(u)} \left(
    \term_2(T) + \term_3(T)
    \right)
    \\
    \overset{\eqref{eq:Err.stokes:T2+T3:darcy}}\lesssim
    \left(
    \sum_{T \in \Thd(u)} h_T^{r(m+1)} \CfT^{-\frac{r'}2} \mu_T^{r'} \kappaT^{-\frac{r'}2}
    \seminorm{W^{m+2,r}(T)^d}{u}^r
    \right)^{\frac1{r'}}
    \norm{\nu,h}{\underline{v}_h}.
    \end{multline*}
\end{proof}

\subsection{Consistency of the Darcy term}

\begin{lemma}[Regime-dependent consistency estimate for the Darcy term]\label{lem:Err.darcy.h}
  Let Assumption~\ref{ass:regularity} hold and define the consistency error linear form $\Err_{\darcy,h} : \VhZ \to \Real$ for $a_{\nu,h}$ such that, for all $\underline{v}_h \in \VhZ$,
  \begin{equation}\label{eq:Err.darcy.h}
    \Err_{\darcy,h}(\underline{v}_h)
    \coloneq
    \int_\Omega \nu u\cdot\PDh\underline{v}_h
    - a_{\nu,h}(\Ih u, \underline{v}_h).
  \end{equation}
  Then, recalling~\eqref{eq:mathfrak.D}, it holds, for all $\underline{v}_h \in \Vh$, further assuming $\DTz \underline{v}_T = 0$ for all $T \in \Th$ if $k = 0$,
  \begin{equation}\label{eq:Err.darcy.h:estimate}
  \Err_{\darcy,h}(\underline{v}_h)
  \lesssim \mathfrak{D}_{\stokes,1}^{\frac1{r'}} \norm{\mu,r,h}{\underline{v}_h}
  + \left(
  \mathfrak{D}_{\stokes,2}^{\frac12} + \mathfrak{D}_{\darcy}^{\frac12}
  \right)\norm{\nu,h}{\underline{v}_h}.
  \end{equation}
\end{lemma}

\begin{proof}
  Expanding $a_{\nu,h}$ according to its definition \eqref{eq:adh}, we can decompose the consistency error as follows:
  \[
  \Err_{\darcy,h}(\underline{v}_h)
  = \sum_{T \in \Th} \nu_T \int_T (u - \PDT \huline{u}_T) \cdot \PDT \underline{v}_T
  + \sum_{T \in \Th} \nu_T s_{\darcy,T}(\huline{u}_T, \underline{v}_T)
  \coloneq \term_1 + \term_2.
  \]
  As in the proof of Lemma \ref{lem:Err.stokes.h}, we denote by $\term_i(T)$ the contribution of an element $T$ to $\term_i$ and we bound each
  of these contributions according to the local regime.
  Proceeding as in the proof of \cite[Lemma~13]{Di-Pietro.Droniou:23} (see, in particular, Eqs.~(54) and~(56) therein), we get
  \begin{equation}\label{eq:Edh:est.darcy}
    \term_1(T) + \term_2(T)
    \lesssim h_T^{m+1} \nu_T^{\frac12} \seminorm{H^{m+1}(T)^d}{u} \, \nu_T^{\frac12} \norm{\darcy,T}{\underline{v}_T}
    \quad \text{if $T \in \Thd(u)$}.
  \end{equation}

  Let us now consider the case $\CfT < 1$ (Stokes-dominated regime).
  We start with the second term as follows:
  \begin{align}
    \term_2(T)
    \overset{\eqref{eq:adT},\eqref{eq:def.Ih}}&=
    \nu_T \lambda_T \int_T \nproj{k}{T} (\PDT \huline{u}_T - u) \cdot \nproj{k}{T} (\PDT \underline{v}_T - v_T)
    \nonumber\\
    &\quad
    + \nu_T h_T \int_{\partial T} (\PDT \huline{u}_T - \cancel{\lproj{k}{\partial T}} u) \cdot (\PDT \underline{v}_T - v_{\partial T})
    \nonumber\\
    &\lesssim
    \nu_T \norm{L^{r'}(T)^d}{\PDT \huline{u}_T -u} ~ \norm{L^r(T)^d}{\PDT \underline{v}_T - v_T}
    \nonumber\\
    &\quad
    + \nu_T h_T \, h_T^{\frac{1}{r'}} \norm{L^{r'}(\partial T)}{ \PDT \huline{u}_T - u} ~ h_T^{\frac{1-r}{r}}\norm{L^r(\partial T)}{\PDT \underline{v}_T - v_{\partial T}}.
    \label{eq:consistency.darcy.Cf<1.T2-1}
  \end{align}
  where the cancellation of the projector is justified noticing that $(\PDT \underline{v}_T - v_{\partial T})_{|F}\in\Poly{k}(F)^d$ for all $F\in\FT$.
  In the second step, we have used $(r',r)$-H\"{o}lder inequalities together with the $L^q$-boundedness of $\nproj{k}{T}$ with $q \in \{ r', r\}$ and the fact that $1 \overset{\eqref{eq:r'}}= h_T^{\frac{1}{r'}} \, h_T^{\frac{1-r}{r}}$.
  Let us consider the last term. Inserting $\pm v_T$ and using a triangle and a discrete trace inequality, we have
  \[
  \begin{aligned}
    &h_T^{\frac{1-r}{r}}\norm{L^r(\partial T)}{\PDT \underline{v}_T - v_{\partial T}}
    \\
    &\quad
    \begin{aligned}[t]
       &\le h_T^{-1}\norm{L^r(T)}{\PDT \underline{v}_T - v_T}
       + h_T^{\frac{1-r}{r}}\norm{L^r(\partial T)}{v_T-v_{\partial T}}
       \\
       \overset{\eqref{eq:poincare.PDT},\eqref{eq:norm.1.q.h}}&\lesssim
       \norm{1,r,T}{\underline{v}_T}.
     \end{aligned}
  \end{aligned}
  \]
  Plugging this into \eqref{eq:consistency.darcy.Cf<1.T2-1} yields
  \begin{align}
    \term_2(T) &\lesssim
    \nu_T h_T \left(\norm{L^{r'}(T)^d}{\PDT \huline{u}_T -u} +  h_T^{\frac{1}{r'}} \norm{L^{r'}(\partial T)}{ \PDT \huline{u}_T - u}\right)\norm{1,r,T}{\underline{v}_T}
    \nonumber\\
    &\lesssim \nu_T h_T^{m+2} \seminorm{W^{m+1,r'}(T)^d}{u} ~ \norm{1,r,T}{\underline{v}_T}
    \nonumber\\
    \overset{\eqref{eq:CfT}}&\le
    \nu_T^{\frac1{r'}} h_T^{m+\frac2{r'}} \norm{L^\infty(T)}{|\nabla u|^{r-2}}^{\frac1r} \CfT^{\frac1r} \seminorm{W^{m+1,r'}(T)^d}{u}
    ~ \mu_T^{\frac1r} \norm{1,r,T}{\underline{v}_T},
    \label{eq:consistency.darcy.Cf<1.T2}
  \end{align}
  where, in the second step, we have invoked the approximation properties of $\PDT \circ \IT$ detailed in~\cite[Proposition~6]{Di-Pietro.Droniou:23} for the Hilbertian case (in the non-Hilbertian case, the proof relies on \cite[Lemma 1.43]{Di-Pietro.Droniou:20} using the bound \cite[(A.21)]{Di-Pietro.Droniou.ea:24}).
  The conclusion is obtained writing $\nu_T=\nu_T^{\frac{1}{r'}}\nu_T^{\frac1r}$ and using \eqref{eq:CfT}.

  For $\term_1(T)$ we treat the cases $k\ge 1$ and $k=0$ separately.
  We first use~\eqref{eq:pi.k-1.PDT=vT}, the definition \eqref{eq:def.Ih} of the interpolator, and the relation \eqref{eq:nproj.circ.lproj} (with $m=k$) to see that $\int_T(\PDT\huline{u}_T-u)\cdot\Psi=0$ for all $\Psi\in\Poly{k-1}(T)^d$.
  For $k\ge 1$, we can apply this relation to $\Psi= \lproj{0}{T} v_T\in \Poly{0}(T)^d \subset \Poly{k-1}(T)^d$ to get
  \begin{equation}\label{eq:consistency.darcy.Cf<1.k>0.T1}
    \begin{aligned}
      \term_1(T)
      &= \nu_T \int_T  (u - \PDT \huline{u}_T) \cdot (\PDT \underline{v}_T - \lproj{0}{T} v_T)
      \\
      &\le \nu_T \norm{L^{r'}(T)^d}{u - \PDT \huline{u}_T}
      \\
      &\quad \times \left(
      \norm{L^r(T)^d}{\PDT \underline{v}_T - v_T}
      + \norm{L^r(T)^d}{v_T - \lproj{0}{T} v_T}
      \right)
      \\
      &\lesssim
      \nu_T h_T^{m+1} \seminorm{W^{m+1,r'}(T)^d}{u}
      ~ h_T \norm{1,r,T}{\underline{v}_T}
      \\
      \overset{\eqref{eq:CfT}}&\le
      \nu_T^{\frac{1}{r'}} h_T^{m+\frac{2}{r'}} \norm{L^\infty(T)}{|\nabla u|^{r-2}}^{\frac1r}\CfT^{\frac1r} \seminorm{W^{m+1,r'}(T)^d}{u}
      ~ \mu_T^{\frac1r}\norm{1,r,T}{\underline{v}_T},
    \end{aligned}
  \end{equation}
  where, in the third step, we have used~\eqref{eq:poincare.PDT} for the first term in parenthesis together with a local Poincaré--Wirtinger inequality $\norm{L^r(T)^d}{v_T - \lproj{0}{T} v_T} \lesssim h_T \norm{L^r(T)^{d\times d}}{\nabla v_T}$ for the second.

  Consider now the case $k=0$. Since $\PDTz \underline{v}_T\in\Poly{0}(T)^d$, we can write $\int_T (u - \PDTz \huline{u}_T) \cdot \PDTz \underline{v}_T= \int_T (\lproj{0}{T}u - \PDTz \huline{u}_T) \cdot \PDTz \underline{v}_T$. Setting $\varphi_u\coloneq (\lproj{0}{T}u - \PDTz \huline{u}_T)\cdot (x - x_T)\in\Poly{1}(T)$ with $x_T\in T$,
  we have $\nabla \varphi_u=\lproj{0}{T}u - \PDTz \huline{u}_T$, and the definition \eqref{eq:PDT} of $\PDTz \underline{v}_T$ thus yields
  \begin{equation}\label{eq:darcy.est.T1.1}
    \begin{aligned}
      \term_1(T)
      &=-\nu_T\int_T\DTz\underline{v}_T \, \varphi_u
      +\nu_T\int_{\partial T} (v_{\partial T}\cdot n_{\partial T}) \, \varphi_u
      \\
      &\le \nu_T\norm{L^2(T)}{\DTz\underline{v}_T} ~ \norm{L^2(T)}{\varphi_u}
      + \nu_Th_T^{\frac12}\norm{L^2(\partial T)^d}{v_{\partial T}} ~ h_T^{-1}\norm{L^2(T)^d}{\varphi_u},
    \end{aligned}
  \end{equation}
  where we have used Cauchy--Schwarz and discrete trace inequalities in the conclusion.
  Using an $(\infty,2)$-H\"{o}lder inequality together with $| x - x_T | \le h_T$ followed by the approximation properties of $\PDTz\circ\underline{I}_h^0$ and $\lproj{0}{T}$, we have
  \begin{equation}\label{eq:est.phiw}
    \norm{L^2(T)}{\varphi_u}
    \le h_T {\norm{L^2(T)^d}{\lproj{0}{T}u - \PDTz \huline{u}_T}}
    \lesssim h_T^2\seminorm{H^1(T)^d}{u}.
  \end{equation}
  By \eqref{eq:equiv.norms.1}, $h_T^{\frac12}\norm{L^2(\partial T)^d}{v_{\partial T}}\lesssim  \norm{\darcy,T}{\underline{v}_T}$.
  Plugging this estimate into \eqref{eq:darcy.est.T1.1}, invoking \eqref{eq:est.phiw}, and recalling that $\DTz \underline{v}_T = 0$ by assumption, we get
  \begin{equation}\label{eq:est.T1.k=0}
    \term_1(T) \lesssim
    \nu_T^{\frac12}h_T\seminorm{H^1(T)^d}{u} ~ \nu_T^{\frac12}\norm{\darcy,T}{\underline{v}_T}.
  \end{equation}
  The conclusion follows from \eqref{eq:Edh:est.darcy}, \eqref{eq:consistency.darcy.Cf<1.T2}, \eqref{eq:consistency.darcy.Cf<1.k>0.T1}, \eqref{eq:est.T1.k=0}, and appropriate Cauchy--Schwarz and H\"older inequalities on the sums.
\end{proof}

\subsection{Consistency of the velocity-pressure coupling term}

\begin{lemma}[Regime-dependent consistency estimate for the velocity-pressure coupling term]\label{lem:Err.b.h}
  Let Assumption~\ref{ass:regularity} hold and define the consistency error linear form $\Err_{b,h} : \VhZ \to \Real$ for $b_h$ such that, for all $\underline{v}_h \in \VhZ,$
  \begin{equation}\label{eq:Err.b.h}
    \Err_{b,h}(\underline{v}_h)
    \coloneq
    \int_\Omega \nabla p\cdot\PDh\underline{v}_h
    - b_h (\underline{v}_h, \lproj{k}{h} p).
  \end{equation}
  Then, recalling~\eqref{eq:mathfrak.C}, it holds, for all $\underline{v}_h \in \Vh$,
  \begin{equation}\label{eq:Err.b.h:estimate}
    \Err_{b,h}(\underline{v}_h)
    \lesssim \mathfrak{C}_{\stokes}^{\frac1{r'}} \norm{\mu,r,h}{\underline{v}_h}
    + \mathfrak{C}_{\darcy}^{\frac12} \norm{\nu,h}{\underline{v}_h}.
  \end{equation}
\end{lemma}

\begin{proof}
  Expanding $b_h$ according to its definition~\eqref{eq:bh}, we have
  \[
  \Err_{b,h}(\underline{v}_h)
  = \sum_{T \in \Th} \left(
  \int_T \nabla p \cdot \PDT \underline{v}_T
  + \int_T \cancel{\lproj{k}{T}} p \, \DT \underline{v}_T
  - \int_{\partial T} p \, (v_{\partial T} \cdot n_{\partial T}),
  \right)
  \]
  where the cancellation of the projector is justified by its definition observing that $\DT \underline{v}_T \in \Poly{k}(T)$, and
  the last term adds to zero since both $p$ and the face components of $\underline{v}_h$ are single-valued at interfaces and $v_F=0$ whenever $F \in \Fhb$.
  Set, for the sake of brevity, $\check{p}_T \coloneq \lproj{k+1}{T} p$ for all $T \in \Th$.
  Applying the definition~\eqref{eq:PDT} of the Darcy velocity reconstruction to $(\varphi,w) = (\check{p}_T, 0)$, we get
  \[
  \int_T \nabla \check{p}_T \cdot \PDT \underline{v}_T
  + \int_T  \check{p}_T \, \DT \underline{v}_T
  - \int_{\partial T} \check{p}_{T} \, (v_{\partial T} \cdot n_{\partial T})
  = 0.
  \]
  Subtracting the above expression from $\Err_{b,h}(\underline{v}_h)$ and integrating by parts element by element, we arrive at
  \begin{multline}\label{eq:E.b.h:basic}
    \Err_{b,h}(\underline{v}_h)
    = \sum_{T \in \Th}
    \cancel{
      \int_T (p - \check{p}_T) \, ( \DT \underline{v}_T - \nabla \cdot \PDT \underline{v}_T )
    }
    \\
    + \sum_{T \in \Th} \underbrace{%
      \int_{\partial T} (\check{p}_T - p) \, ( v_{\partial T} - \PDT \underline{v}_T ) \cdot n_{\partial T},
    }_{\eqcolon \term(T)}
  \end{multline}
  where the cancellation comes from the fact that the second factor is in $\Poly{k}(T)$ and the first one is $L^2$-orthogonal to this space.

  For $T \in \Ths(w)$, we use an $(r,r',\infty)$-H\"{o}lder inequality together with $\norm{L^\infty(\partial T)}{n_{\partial T}} \le 1$ and $1 \overset{\eqref{eq:r'}}= h_T^{\frac1{r'}} \, h_T^{\frac{1-r}{r}}$ to write
  \begin{equation}\label{eq:E.b.h:T(T):stokes}
    \begin{aligned}
      \term(T)
      &\lesssim h_T^{\frac{1}{r'}}\norm{L^{r'}(\partial T)}{\check{p}_T - p} \, h_T^{\frac{1-r}{r}} \norm{L^r(\partial T)^d}{v_{\partial T} - \PDT \underline{v}_T}
      \\
      &\lesssim h_T^{m+1} \seminorm{W^{m+1,r'}(T)}{p} \\
      &\quad \times \left(
      h_T^{\frac{1-r}{r}} \norm{L^r(\partial T)^d}{v_{\partial T} - v_T}
      + h_T^{-1}\norm{L^r(T)^d}{\PDT \underline{v}_T - v_T}
      \right)
      \\
      \overset{\eqref{eq:norm.1.q.h},\,\eqref{eq:poincare.PDT}}&\lesssim
      h_T^{m+1} \mu_T^{-\frac1r} \seminorm{W^{m+1,r'}(T)}{p}
      ~ \mu_T^{\frac1r}\norm{1,r,T}{\underline{v}_T},
    \end{aligned}
  \end{equation}
  where we have used the approximation properties of $\lproj{k+1}{T}$ together with triangle and discrete trace inequalities in the second step.

  If $T \in \Thd(w)$, on the other hand, we have
  \begin{equation}\label{eq:E.b.h:T(T):darcy.prelim}
    \begin{aligned}
      \term(T)
      &\le h_T^{-\frac12} \norm{L^2(\partial T)}{\check{p}_T - p}\,
      h_T^{\frac12} \norm{L^2(\partial T)^d}{v_{\partial T} - \PDT \underline{v}_T}
      \\
      &\lesssim h_T^{m+1}  \seminorm{H^{m+2}(T)}{p} \,
    h_T^{\frac12}\norm{L^2(\partial T)^d}{v_{\partial T} - \PDT \underline{v}_T},
    \end{aligned}
  \end{equation}
  where we have used the approximation properties of $\lproj{k+1}{T}$ to conclude.
  Let $F\in\FT$. If $F\not\in\Fhb$, then \eqref{eq:adT} and \eqref{eq:norm.darcy} give
  \[
  h_T^{\frac12}\norm{L^2(F)^d}{v_F - \PDT \underline{v}_T}\le \norm{\darcy,T}{\underline{v}_T}.
  \]
  If $F\in\Fhb$, then $v_F=0$ so
  \[
  \begin{aligned}
    h_T^{\frac12}\norm{L^2(F)^d}{v_F- \PDT \underline{v}_T}
    &= h_T^{\frac12} \norm{L^2(F)^d}{\PDT \underline{v}_T}
    \\
    &\lesssim \norm{L^2(T)^d}{\PDT \underline{v}_T}
    \overset{\eqref{eq:adT},\,\eqref{eq:norm.darcy}}\le \norm{\darcy,T}{\underline{v}_T},
  \end{aligned}
  \]
  where the first inequality follows from a discrete trace inequality. All in all, coming back to \eqref{eq:E.b.h:T(T):darcy.prelim} we obtain
  \begin{equation}\label{eq:E.b.h:T(T):darcy}
    \term(T)
    \lesssim h_T^{m+1} \nu_T^{-\frac12} \seminorm{H^{m+2}(T)}{p} \, \nu_T^{\frac12}\norm{\darcy,T}{\underline{v}_T}.
  \end{equation}

  Plugging the estimates \eqref{eq:E.b.h:T(T):stokes} and \eqref{eq:E.b.h:T(T):darcy} into \eqref{eq:E.b.h:basic} and using, respectively, H\"{o}lder and Cauchy--Schwarz inequalities on the sums over $T \in \Ths(w)$ and $T \in \Thd(w)$, the conclusion follows.
\end{proof}

\subsection{Proof of Theorem~\ref{thm:err.estimate}}\label{sec:proof.theorem}

\begin{proof}[Proof of Theorem~\ref{thm:err.estimate}]
  Taking $\underline{v}_h = \underline{e}_h$ in \eqref{eq:discrete:momentum},
  $q_h = \epsilon_h \coloneq p_h - \widehat{p}_h$ with $\widehat{p}_h \coloneq \lproj{k}{h} p$ in \eqref{eq:discrete:mass},
  summing the two equations, and subtracting the quantity
  $a_{\mu,h}(\huline{u}_h, \underline{e}_h) + a_{\nu,h}(\huline{u}_h, \underline{e}_h) + b_h(\underline{e}_h, \widehat{p}_h) - b_h(\huline{u}_h, \epsilon_h)$ from both sides of the resulting expression, we end up with
  \begin{equation}\label{eq:error:basic}
    a_{\mu,h}(\underline{u}_h, \underline{e}_h)
    - a_{\mu,h}(\huline{u}_h, \underline{e}_h)
    + a_{\nu,h}(\underline{e}_h, \underline{e}_h)
    = \Err_h(u,p; \underline{e}_h)
  \end{equation}
  with
  \begin{multline*}
  \Err_h(u,p; \underline{e}_h)
  \coloneq
  \int_\Omega f \cdot \PDh \underline{e}_h
  - a_{\mu,h}(\huline{u}_h, \underline{e}_h)
  - a_{\nu,h}(\huline{u}_h, \underline{e}_h)
  - b_h(\underline{e}_h, \widehat{p}_h)
  \\
  +
  \cancel{%
    \int_\Omega g \, \epsilon _h
    + b_h(\huline{u}_h, \epsilon_h)
  },
  \end{multline*}
  where the cancellation follows from \eqref{eq:fortin} after observing that $g = \nabla \cdot u$ almost everywhere in $\Omega$.

  A lower bound for the left-hand side of~\eqref{eq:error:basic} is provided by~\eqref{eq:stability}.
  Let us now find an upper bound for $\Err_h(u,p; \underline{e}_h)$ in terms of the error components defined by~\eqref{eq:mathfrak.S},~\eqref{eq:mathfrak.D}, and~\eqref{eq:mathfrak.C}.
  Noticing that $f = - \nabla \cdot \sigma(\nabla u) + \nu u + \nabla p$ almost everywhere in $\Omega$ and recalling the definitions \eqref{eq:Err.stokes.h}, \eqref{eq:Err.darcy.h}, and \eqref{eq:Err.b.h} of the Stokes, Darcy, and velocity-pressure coupling consistency errors, we have
  \[
  \Err_h(u,p; \underline{e}_h)
  = \Err_{\stokes,h}(\underline{e}_h)
  + \Err_{\darcy,h}(\underline{e}_h)
  + \Err_{b,h}(p;\underline{e}_h).
  \]
  Collecting the bounds~\eqref{eq:Err.stokes.h:estimate},~\eqref{eq:Err.darcy.h:estimate} (further noticing that $\DT \underline{e}_T = \DT\underline{u}_T-\DT\IT u\overset{\eqref{eq:discrete:mass},\eqref{eq:DT.proj.div}}=\lproj{k}{T}g - {\lproj{k}{T}(\nabla\cdot u)} = 0$ for all $T \in \Th$), and~\eqref{eq:Err.b.h:estimate}, we infer the existence of $C_{\rm cons} > 0$ independent of $h$, $\mu$, and $\nu$ such that
  \begin{equation}\label{eq:consistency}
    \begin{aligned}
      &C_{\rm cons}^{-1} \Err_h(u,p; \underline{e}_h)
      \\
      &\quad\le
      \left(
      \mathfrak{S}_{\stokes,1}^{\frac1{r'}}
      + \mathfrak{S}_{\stokes,2}^{\frac1r}
      + \mathfrak{D}_{\stokes,1}^{\frac1{r'}}
      + \mathfrak{C}_{\stokes}^{\frac1{r'}}
      \right) \norm{\mu,r,h}{\underline{e}_h}
      \\
      &\qquad
      + \left(
      \mathfrak{S}_{\darcy,1}^{\frac1{r'}}
      + \mathfrak{S}_{\darcy,2}^{\frac12}
      +\mathfrak{D}_{\stokes,2}^{\frac12}
      + \mathfrak{D}_{\darcy}^{\frac12}
      + \mathfrak{C}_{\darcy}^{\frac12}
      \right) \norm{\nu,h}{\underline{e}_h}
      \\
      &\quad\le
      \frac{C_{\rm stab}}{2 C_{\rm cons}} \left(
      \alpha_{\mu} \norm{\mu,r,h}{\underline{e}_h}^{q_r}
      + \norm{\nu,h}{\underline{e}_h}^2
      \right)
      \\
      &\qquad
      + \left(
      \frac{2 C_{\rm cons}}{\alpha_\mu C_{\rm stab} q_r}
      \right)^{\frac{q_r'}{q_r}} \frac{1}{q_r'} \left(
      \mathfrak{S}_{\stokes,1}^{\frac1{r'}}
      + \mathfrak{S}_{\stokes,2}^{\frac1r}
      + \mathfrak{D}_{\stokes,1}^{\frac1{r'}}
      + \mathfrak{C}_{\stokes}^{\frac1{r'}}
      \right)^{q_r'}
      \\
      &\qquad
      + \frac{C_{\rm cons}}{2 C_{\rm stab}} \left(
      \mathfrak{S}_{\darcy,1}^{\frac1{r'}}
      + \mathfrak{S}_{\darcy,2}^{\frac12}
      + \mathfrak{D}_{\stokes,2}^{\frac12}
      + \mathfrak{D}_{\darcy}^{\frac12}
      + \mathfrak{C}_{\darcy}^{\frac12}
      \right)^2,
    \end{aligned}
  \end{equation}
  where the second step follows from the generalized Young's inequality $ab \le \frac{\epsilon^q}{q} a^q + \frac{1}{\epsilon^{q'}q'} b^{q'}$ with, respectively,
  $q = q_r$ and $\epsilon = \left(\frac{\alpha_\mu C_{\rm stab} q_r}{2 C_{\rm cons}}\right)^{\frac1{q_r}}$ for the first term
  and $q = 2$ and $\epsilon = \left( \frac{C_{\rm stab}}{C_{\rm cons}} \right)^{\frac12}$ for the second term.
  Plugging~\eqref{eq:stability} and~\eqref{eq:consistency} into~\eqref{eq:error:basic} and using standard inequalities for sums of powers yields~\eqref{eq:err.estimate}.
\end{proof}

%------------------------------------------------------------------------------%

\section{Numerical results}\label{sec:numerical.results}

In this section we numerically validate the proposed method \eqref{eq:discrete} for different values of $(r,\mu,\nu)$
on both standard and polygonal meshes.
Our implementation is based on the \texttt{HArDCore} library\footnote{\url{https://github.com/jdroniou/HArDCore}} and makes extensive use of the linear algebra \texttt{Eigen} open-source library \cite{Guennebaud.Jacob.ea:10}. Since we are dealing with a nonlinear scheme,
we use the  Newton's algorithm with a stopping tolerance of $10^{-10}$; the linearised equations are solved by means of the direct solver Pardiso \cite{Schenk.Gartner.ea:01}.

Following \cite{Botti.Di-Pietro.ea:18,Di-Pietro.Droniou:23}, we solve problem \eqref{eq:strong} on the unit square domain $\Omega = (0,1)^2$ for constant parameters $\mu$ and $\nu$, and define the global friction coefficient $\CfOmg\coloneq \frac{\nu}{\mu}$ with the convention that $\CfOmg=+ \infty$ when $\mu=0$. For the sake of simplicity, in defining $\CfOmg$ we have taken characteristic length equal to $1$.
Letting $\chi_\sstokes(\CfOmg)\coloneq\exp(-\CfOmg)$, the pressure and velocity are set as follows:
$p(x,y) = \sin(\pi x) \sin(\pi y)$ and
$u = \chi_\sstokes(\CfOmg)u_\sstokes + (1-\chi_\sstokes(\CfOmg)) u_\ddarcy$,
where $u_\sstokes$ and $u_\ddarcy$ are the velocities obtained in the $r$-Stokes ($\CfOmg=0$) and Darcy ($\CfOmg=+\infty$) limits, and are given by
\[
\begin{aligned}
  u_\sstokes(x,y)
  &=\begin{bmatrix}
  \phantom{-}\sin(\pi x)\cos(\pi y)\\
  -\cos(\pi x)\sin(\pi y)
  \end{bmatrix}\quad\forall (x,y)\in\Omega,
  \\
  u_\ddarcy
  &=\begin{cases}
  -\nu^{-1}\nabla p & \text{if $\nu>0$} ,\\
  {0} & \text{otherwise}.
  \end{cases}
\end{aligned}
\]
Observe that $\nabla \cdot u_\sstokes=0$ and that $\nu u_\ddarcy + \nabla p =0$,
which are the expected relations, respectively,
for a solution of the incompressible $r$-Stokes equation
and for a solution of the Darcy equation in mixed form.
We consider computations over two $h$-refined mesh families
(simplicial and hexagonal).
\begin{figure}
  \centering
  \includegraphics[height=4cm]{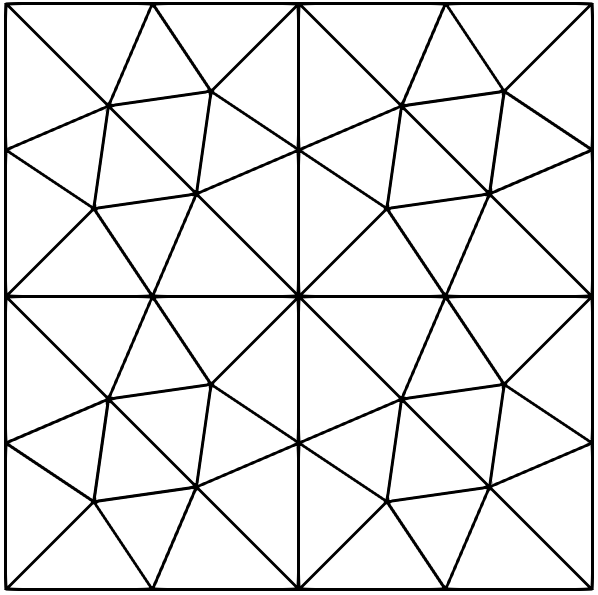}
  \hspace{1cm}
  \includegraphics[height=4cm]{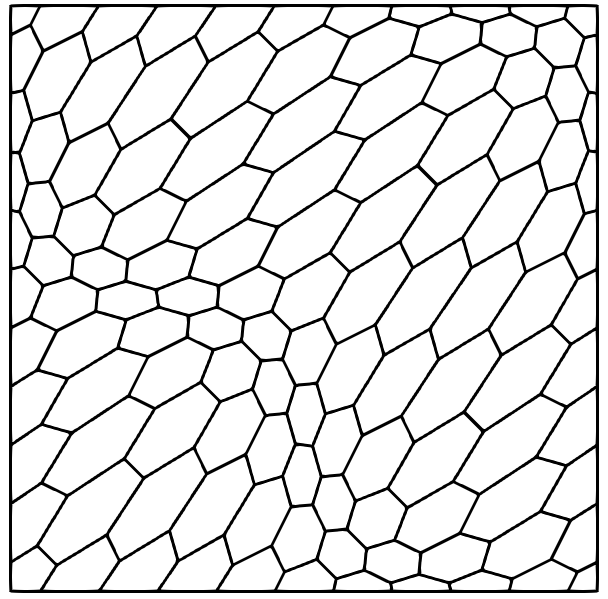}
  \caption{Coarsest triangular and hexagonal meshes used in the numerical tests.\label{fig:meshes:coarsest}}
\end{figure}
Figure~\ref{fig:meshes:coarsest} shows the coarsest mesh for each family. Recalling the definition of the velocity error $\underline{e}_h \coloneq \underline{u}_h - \huline{u}_h$,
we monitor the quantity $\EnergyError$ estimated in Theorem~\ref{thm:err.estimate}.

\begin{table}\centering
  \subcaptionbox{Stokes-dominated}[0.4\textwidth]{
    \begin{tabular}{c|ccc}
      \toprule
      \diagbox{$r$}{$k$} & 0 & 1 & 2 \\
      \midrule
      1.75 & 1.50 & 3.00 & 4.50 \\
      2 & 2.00 & 4.00 & 6.00 \\
      3 & 1.50 & 3.00 & 4.50 \\
      4 & 1.33 & 2.67 & 4.00 \\
      \bottomrule
    \end{tabular}
  }
  \subcaptionbox{Darcy-dominated}[0.4\textwidth]{
    \begin{tabular}{c|ccc}
      \toprule
      \diagbox{$r$}{$k$} & 0 & 1 & 2 \\
      \midrule
      1.75 & 1.50 & 3.00 & 4.50 \\
      2 & 2.00 & 4.00 & 6.00 \\
      3 & 2.00 & 4.00 & 6.00 \\
      4 & 2.00 & 4.00 & 6.00 \\
      \bottomrule
    \end{tabular}
  }
  \caption{Convergence rates for $\alpha_\mu \norm{\mu,r,h}{\underline{e}_h}^{q_r} + \norm{\nu,h}{\underline{e}_h}^2$ predicted by Theorem \ref{thm:err.estimate}, for the values of $r$ and $k$ used in the numerical tests.\label{tab:convergence.rates:numerical.tests}}
\end{table}

\begin{figure}\centering
  \subcaptionbox[B]{$(r,\mu,\nu) = (3,1,1)$.\label{fig:tri:(3,1,1)}}[0.75\textwidth]{\centering
    \includegraphics[height=3.75cm]{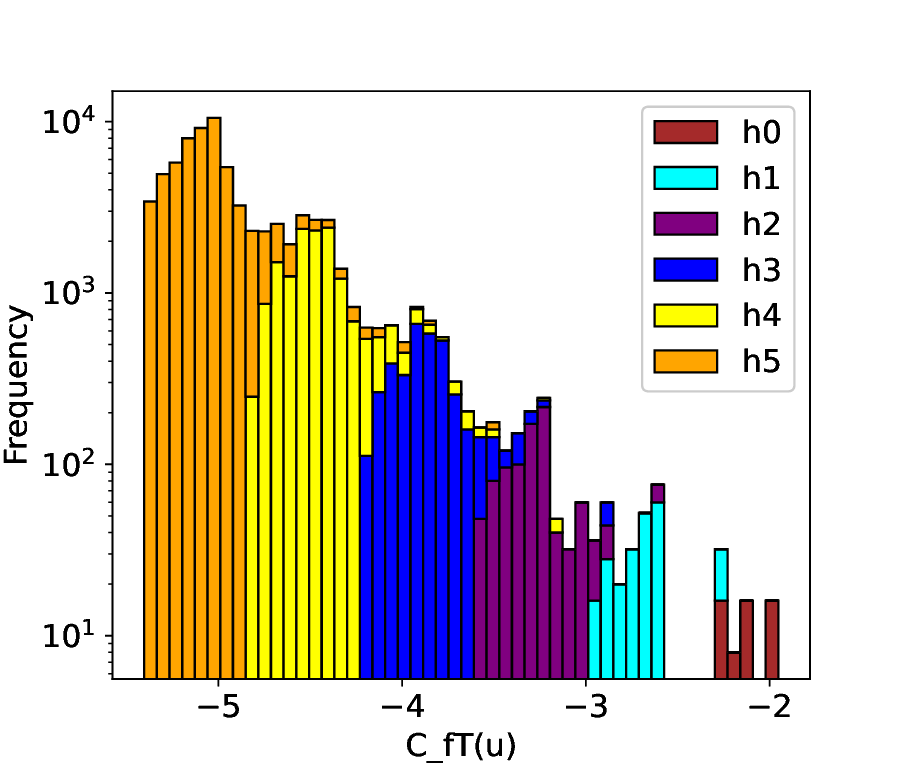}
    \hfill
    \begin{tikzpicture}[scale=0.525]
      \begin{loglogaxis}[legend pos=south east]
        \addplot table[col sep=comma,x=MeshSize,y=EnergyError] {results/test1/r3/tri/tri_k0_mu1_nu2/data_rates.dat}
        node[pos = 0.1, above=2pt]{2.78} % EnergyErrorRate
        node[pos = 0.3, above=2pt]{2.92} % EnergyErrorRate
        node[pos = 0.5, above=2pt]{2.97} % EnergyErrorRate
        node[pos = 0.7, above=2pt]{2.99} % EnergyErrorRate
        node[pos = 0.9, above=2pt]{3.00}; % EnergyErrorRate
        \addlegendentry{$k=0$}

        \addplot table[col sep=comma,x=MeshSize,y=EnergyError] {results/test1/r3/tri/tri_k1_mu1_nu2/data_rates.dat}
        node[pos = 0.1, above=2pt]{3.28} % EnergyErrorRate
        node[pos = 0.3, above=2pt]{3.01} % EnergyErrorRate
        node[pos = 0.5, above=2pt]{3.00} % EnergyErrorRate
        node[pos = 0.7, above=2pt]{3.00} % EnergyErrorRate
        node[pos = 0.9, above=2pt]{3.00}; % EnergyErrorRate
        \addlegendentry{$k=1$}

        \addplot table[col sep=comma,x=MeshSize,y=EnergyError] {results/test1/r3/tri/tri_k2_mu1_nu2/data_rates.dat}
        node[pos = 0.1, above=2pt]{4.09} % EnergyErrorRate
        node[pos = 0.3, above=2pt]{4.26} % EnergyErrorRate
        node[pos = 0.5, above=2pt]{4.35} % EnergyErrorRate
        node[pos = 0.7, above=2pt]{4.40} % EnergyErrorRate
        node[pos = 0.9, above=2pt]{4.44}; % EnergyErrorRate
        \addlegendentry{$k=2$}

      \end{loglogaxis}
    \end{tikzpicture}
  }\vspace{0.1cm}\\
  %%%tri r=3, mu=10e-4,nu=1
  \subcaptionbox[B]{$(r,\mu,\nu) = (3,10^{-3},1)$.\label{fig:tri:(3,1e-3,1)}}[0.75\textwidth]{\centering
    \includegraphics[height=3.75cm]{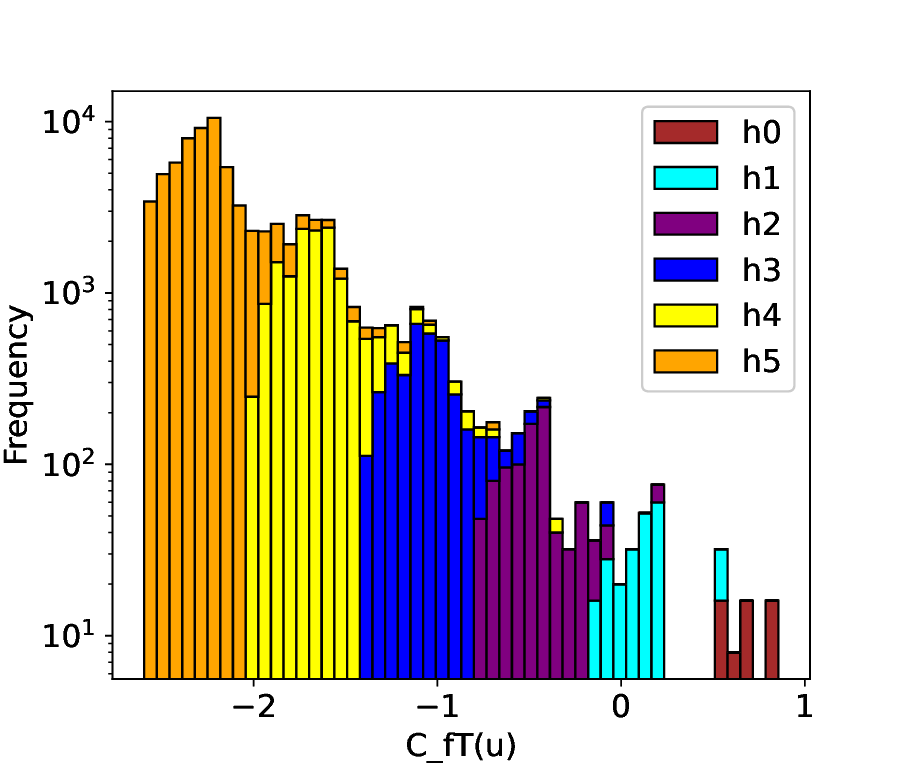}
    \hfill
    \begin{tikzpicture}[scale=0.525]
      \begin{loglogaxis}[legend pos=south east]
        \addplot table[col sep=comma,x=MeshSize,y=EnergyError] {results/test1/r3/tri/tri_k0_mu4_nu2/data_rates.dat}
        node[pos = 0.1, above=2pt]{2.41} % EnergyErrorRate
        node[pos = 0.3, above=2pt]{2.79} % EnergyErrorRate
        node[pos = 0.5, above=2pt]{2.88} % EnergyErrorRate
        node[pos = 0.7, above=2pt]{2.92} % EnergyErrorRate
        node[pos = 0.9, above=2pt]{2.95}; % EnergyErrorRate
        \addlegendentry{$k=0$}

        \addplot table[col sep=comma,x=MeshSize,y=EnergyError] {results/test1/r3/tri/tri_k1_mu4_nu2/data_rates.dat}
        node[pos = 0.15, above=2pt]{4.32} % EnergyErrorRate
        node[pos = 0.40, above=2pt]{4.02} % EnergyErrorRate
        node[pos = 0.60, above=2pt]{3.81} % EnergyErrorRate
        node[pos = 0.80, above=2pt]{3.32} % EnergyErrorRate
        node[pos = 0.95, above=2pt]{2.80}; % EnergyErrorRate
        \addlegendentry{$k=1$}

        \addplot table[col sep=comma,x=MeshSize,y=EnergyError] {results/test1/r3/tri/tri_k2_mu4_nu2/data_rates.dat}
        node[pos = 0.15, above=2pt]{5.93} % EnergyErrorRate
        node[pos = 0.40, above=2pt]{4.83} % EnergyErrorRate
        node[pos = 0.60, above=2pt]{3.99} % EnergyErrorRate
        node[pos = 0.80, above=2pt]{3.62} % EnergyErrorRate
        node[pos = 0.95, above=2pt]{3.54}; % EnergyErrorRate
        \addlegendentry{$k=2$}
      \end{loglogaxis}
    \end{tikzpicture}
  }\vspace{0.1cm}\\
  %%%%tri r=3, mu=10e-6,nu=1
  \subcaptionbox[B]{$(r,\mu,\nu) = (3,10^{-5},1)$.\label{fig:tri:(3,1e-5,1)}}[0.75\textwidth]{\centering
    \includegraphics[height=3.75cm]{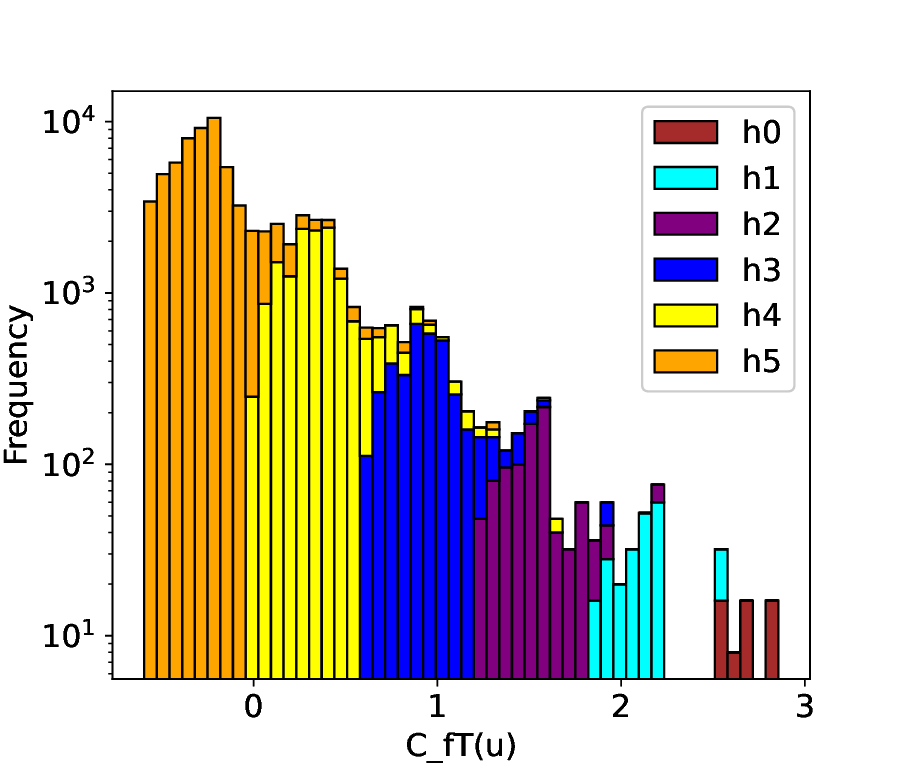}
    \hfill
    \begin{tikzpicture}[scale=0.525]
      \begin{loglogaxis}[legend pos=south east]
        \addplot table[col sep=comma,x=MeshSize,y=EnergyError] {results/test1/r3/tri/tri_k0_mu6_nu2/data_rates.dat}
        node[pos = 0.1, above=2pt]{1.84} % EnergyErrorRate
        node[pos = 0.3, above=2pt]{1.99} % EnergyErrorRate
        node[pos = 0.5, above=2pt]{2.22} % EnergyErrorRate
        node[pos = 0.7, above=2pt]{2.46} % EnergyErrorRate
        node[pos = 0.9, above=2pt]{2.71}; % EnergyErrorRate
        \addlegendentry{$k=0$}

        \addplot table[col sep=comma,x=MeshSize,y=EnergyError] {results/test1/r3/tri/tri_k1_mu6_nu2/data_rates.dat}
        node[pos = 0.1, above=2pt]{3.99} % EnergyErrorRate
        node[pos = 0.3, above=2pt]{4.25} % EnergyErrorRate
        node[pos = 0.5, above=2pt]{4.55} % EnergyErrorRate
        node[pos = 0.7, above=2pt]{4.40} % EnergyErrorRate
        node[pos = 0.9, above=2pt]{3.94}; % EnergyErrorRate
        \addlegendentry{$k=1$}

        \addplot table[col sep=comma,x=MeshSize,y=EnergyError] {results/test1/r3/tri/tri_k2_mu6_nu2/data_rates.dat}
        node[pos = 0.1, above=2pt]{5.98} % EnergyErrorRate
        node[pos = 0.3, above=2pt]{6.52} % EnergyErrorRate
        node[pos = 0.5, above=2pt]{6.93} % EnergyErrorRate
        node[pos = 0.7, above=2pt]{6.51} % EnergyErrorRate
        node[pos = 0.9, above=2pt]{5.91}; % EnergyErrorRate
        \addlegendentry{$k=2$}

      \end{loglogaxis}
    \end{tikzpicture}
  }\vspace{0.1cm}\\
  %tri r=4, mu=10e-6,nu=1
  \subcaptionbox[B]{$(r,\mu,\nu) = (4,10^{-5},1)$.\label{fig:tri:(4,1e-5,1)}}[0.75\textwidth]{\centering
    \includegraphics[height=3.75cm]{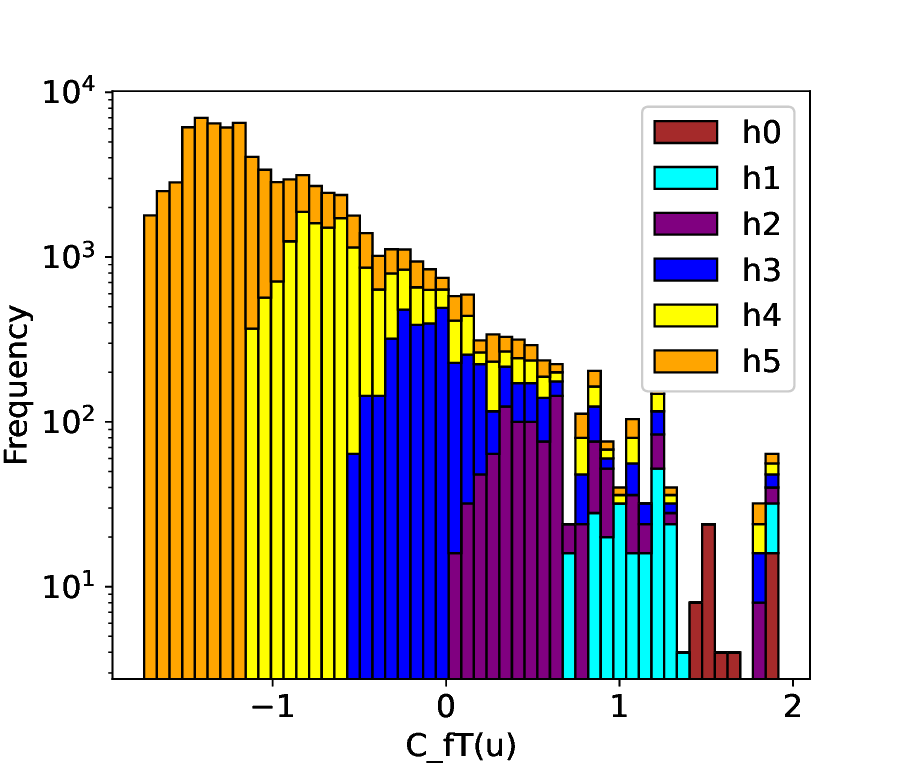}
    \hfill
    \begin{tikzpicture}[scale=0.525]
      \begin{loglogaxis}[legend pos=south east]
        \addplot table[col sep=comma,x=MeshSize,y=EnergyError] {results/test1/r4/tri/tri_k0_mu6_nu2/data_rates.dat}
        node[pos = 0.1, above=2pt]{1.93} % EnergyErrorRate
        node[pos = 0.3, above=2pt]{2.20} % EnergyErrorRate
        node[pos = 0.5, above=2pt]{2.40} % EnergyErrorRate
        node[pos = 0.7, above=2pt]{2.47} % EnergyErrorRate
        node[pos = 0.9, above=2pt]{2.52}; % EnergyErrorRate
        \addlegendentry{$k=0$}

        \addplot table[col sep=comma,x=MeshSize,y=EnergyError] {results/test1/r4/tri/tri_k1_mu6_nu2/data_rates.dat}
        node[pos = 0.1, above=2pt]{4.32} % EnergyErrorRate
        node[pos = 0.3, above=2pt]{4.42} % EnergyErrorRate
        node[pos = 0.5, above=2pt]{4.09} % EnergyErrorRate
        node[pos = 0.7, above=2pt]{3.87} % EnergyErrorRate
        node[pos = 0.9, above=2pt]{3.77}; % EnergyErrorRate
        \addlegendentry{$k=1$}

        \addplot table[col sep=comma,x=MeshSize,y=EnergyError] {results/test1/r4/tri/tri_k2_mu6_nu2/data_rates.dat}
        node[pos = 0.150, above=2pt]{6.59} % EnergyErrorRate
        node[pos = 0.350, above=2pt]{6.57} % EnergyErrorRate
        node[pos = 0.575, above=2pt]{5.92} % EnergyErrorRate
        node[pos = 0.775, above=2pt]{5.53} % EnergyErrorRate
        node[pos = 0.925, above=2pt]{5.00}; % EnergyErrorRate
        \addlegendentry{$k=2$}

      \end{loglogaxis}
    \end{tikzpicture}
  }\vspace{0.1cm}\\
  \caption{Convergence results using the triangular mesh family for different values of $(r,\mu,\nu)$ and for $k\in\{0,1,2\}$. 
    Left column: Distribution of $\CfT$ for all the meshes.
    Right column:  The quantity $\EnergyError[r]$ as a function of $h$ (the numbers on the graphs indicate the slopes between two data points).}
  \label{fig:nt1.tri}
\end{figure}
%
 %Regime test

\begin{figure}\centering
  %%%hexa1 r=3, mu=10e-2,nu=1
  \subcaptionbox[B]{$(r,\mu,\nu) = (3,10^{-1},1)$.\label{fig:hexa:(3,1e-1,1)}}[0.75\textwidth]{\centering
    \includegraphics[height=3.75cm]{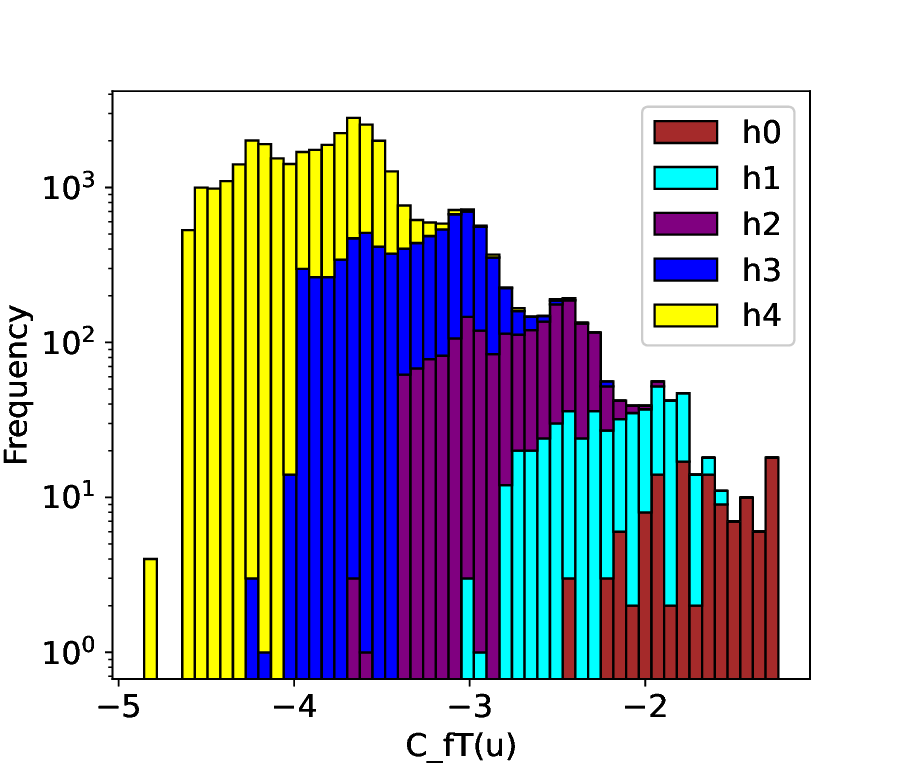}
    \hfill
    \begin{tikzpicture}[scale=0.525]
      \begin{loglogaxis}[legend pos=south east]
        \addplot table[col sep=comma,x=MeshSize,y=EnergyError] {results/test1/r3/hexa1/hexa1_k0_mu2_nu2/data_rates.dat}
        node[pos = 0.125, above=2pt]{0.90} % EnergyErrorRate
        node[pos = 0.375, above=2pt]{0.95} % EnergyErrorRate
        node[pos = 0.625, above=2pt]{1.11} % EnergyErrorRate
        node[pos = 0.875, above=2pt]{1.26}; % EnergyErrorRate
        \addlegendentry{$k=0$}

        \addplot table[col sep=comma,x=MeshSize,y=EnergyError] {results/test1/r3/hexa1/hexa1_k1_mu2_nu2/data_rates.dat}
        node[pos = 0.125, above=2pt]{2.47} % EnergyErrorRate
        node[pos = 0.375, above=2pt]{2.55} % EnergyErrorRate
        node[pos = 0.625, above=2pt]{2.75} % EnergyErrorRate
        node[pos = 0.875, above=2pt]{2.87}; % EnergyErrorRate
        \addlegendentry{$k=1$}

        \addplot table[col sep=comma,x=MeshSize,y=EnergyError] {results/test1/r3/hexa1/hexa1_k2_mu2_nu2/data_rates.dat}
        node[pos = 0.125, above=2pt]{4.08} % EnergyErrorRate
        node[pos = 0.375, above=2pt]{4.14} % EnergyErrorRate
        node[pos = 0.625, above=2pt]{4.27} % EnergyErrorRate
        node[pos = 0.875, above=2pt]{4.34}; % EnergyErrorRate
        \addlegendentry{$k=2$}

      \end{loglogaxis}
    \end{tikzpicture}
  }\vspace{0.1cm}\\
  %%%hexa1 r=3, mu=10e-5,nu=1
  \subcaptionbox[B]{$(r,\mu,\nu) = (3,10^{-4},1)$.\label{fig:hexa:(3,1e-4,1)}}[0.75\textwidth]{\centering
    \includegraphics[height=3.75cm]{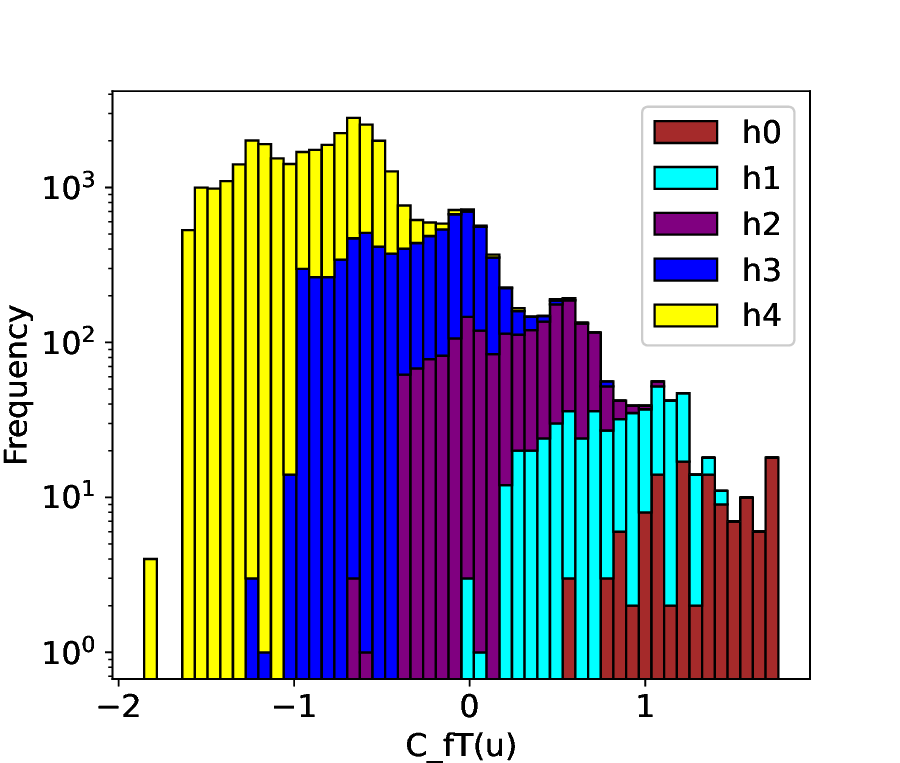}
    \hfill
    \begin{tikzpicture}[scale=0.525]
      \begin{loglogaxis}[legend pos=south east]
        \addplot table[col sep=comma,x=MeshSize,y=EnergyError] {results/test1/r3/hexa1/hexa1_k0_mu5_nu2/data_rates.dat}
        node[pos = 0.150, above=2pt]{2.79} % EnergyErrorRate
        node[pos = 0.425, above=2pt]{3.34} % EnergyErrorRate
        node[pos = 0.675, above=2pt]{2.43} % EnergyErrorRate
        node[pos = 0.900, above=2pt]{1.25}; % EnergyErrorRate
        \addlegendentry{$k=0$}

        \addplot table[col sep=comma,x=MeshSize,y=EnergyError] {results/test1/r3/hexa1/hexa1_k1_mu5_nu2/data_rates.dat}
        node[pos = 0.150, above=2pt]{4.70} % EnergyErrorRate
        node[pos = 0.425, above=2pt]{5.60} % EnergyErrorRate
        node[pos = 0.675, above=2pt]{5.59} % EnergyErrorRate
        node[pos = 0.900, above=2pt]{2.72}; % EnergyErrorRate
        \addlegendentry{$k=1$}

        \addplot table[col sep=comma,x=MeshSize,y=EnergyError] {results/test1/r3/hexa1/hexa1_k2_mu5_nu2/data_rates.dat}
        node[pos = 0.150, above=2pt]{6.97} % EnergyErrorRate
        node[pos = 0.425, above=2pt]{7.39} % EnergyErrorRate
        node[pos = 0.675, above=2pt]{5.94} % EnergyErrorRate
        node[pos = 0.900, above=2pt]{4.13}; % EnergyErrorRate
        \addlegendentry{$k=2$}

      \end{loglogaxis}
    \end{tikzpicture}
  }\vspace{0.1cm}\\
  %%hexa1 r=4, mu=10e-3,nu=1
  \subcaptionbox[B]{$(r,\mu,\nu) = (4,10^{-2},1)$.\label{fig:hexa:(4,1e-2,1)}}[0.75\textwidth]{\centering
    \includegraphics[height=3.75cm]{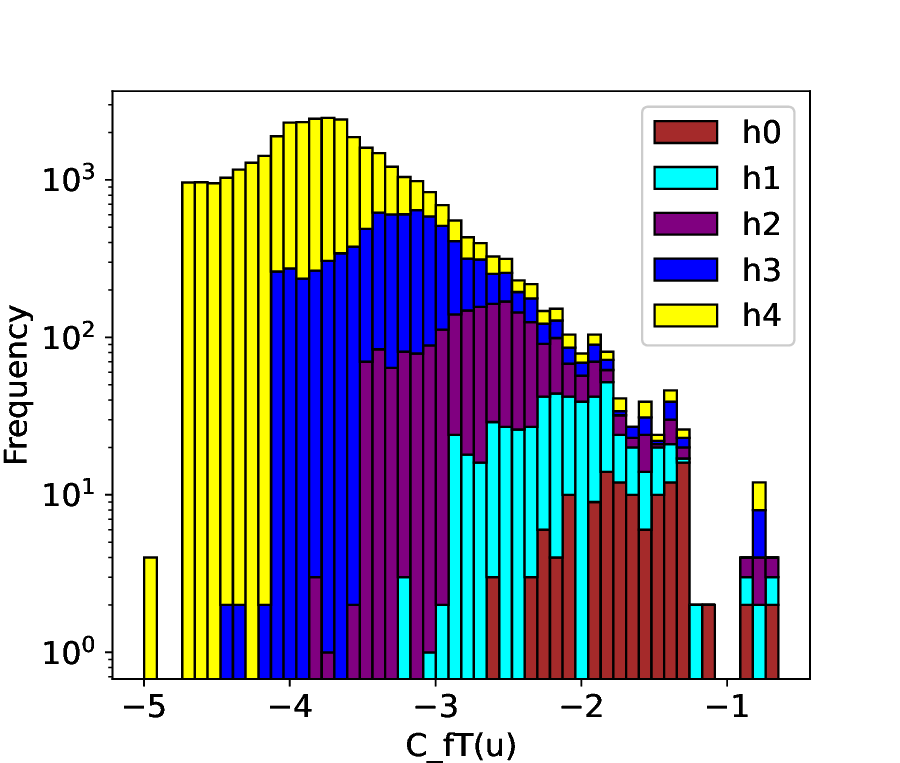}
    \hfill
    \begin{tikzpicture}[scale=0.525]
      \begin{loglogaxis}[legend pos=south east]
        \addplot table[col sep=comma,x=MeshSize,y=EnergyError] {results/test1/r4/hexa1/hexa1_k0_mu3_nu2/data_rates.dat}
        node[pos = 0.125, above=2pt]{0.71} % EnergyErrorRate
        node[pos = 0.375, above=2pt]{0.80} % EnergyErrorRate
        node[pos = 0.625, above=2pt]{0.90} % EnergyErrorRate
        node[pos = 0.875, above=2pt]{1.02}; % EnergyErrorRate
        \addlegendentry{$k=0$}

        \addplot table[col sep=comma,x=MeshSize,y=EnergyError] {results/test1/r4/hexa1/hexa1_k1_mu3_nu2/data_rates.dat}
        node[pos = 0.125, above=2pt]{2.38} % EnergyErrorRate
        node[pos = 0.375, above=2pt]{2.36} % EnergyErrorRate
        node[pos = 0.625, above=2pt]{2.48} % EnergyErrorRate
        node[pos = 0.875, above=2pt]{2.57}; % EnergyErrorRate
        \addlegendentry{$k=1$}

        \addplot table[col sep=comma,x=MeshSize,y=EnergyError] {results/test1/r4/hexa1/hexa1_k2_mu3_nu2/data_rates.dat}
        node[pos = 0.125, above=2pt]{3.79} % EnergyErrorRate
        node[pos = 0.375, above=2pt]{3.73} % EnergyErrorRate
        node[pos = 0.625, above=2pt]{3.89} % EnergyErrorRate
        node[pos = 0.875, above=2pt]{3.95}; % EnergyErrorRate
        \addlegendentry{$k=2$}

      \end{loglogaxis}
    \end{tikzpicture}
  }\vspace{0.1cm}\\
  %%%hexa1 r=4, mu=10e-6,nu=1
  \subcaptionbox[B]{$(r,\mu,\nu) = (4,10^{-5},1)$.\label{fig:hexa:(4,1e-5,1)}}[0.75\textwidth]{\centering
    \includegraphics[height=3.75cm]{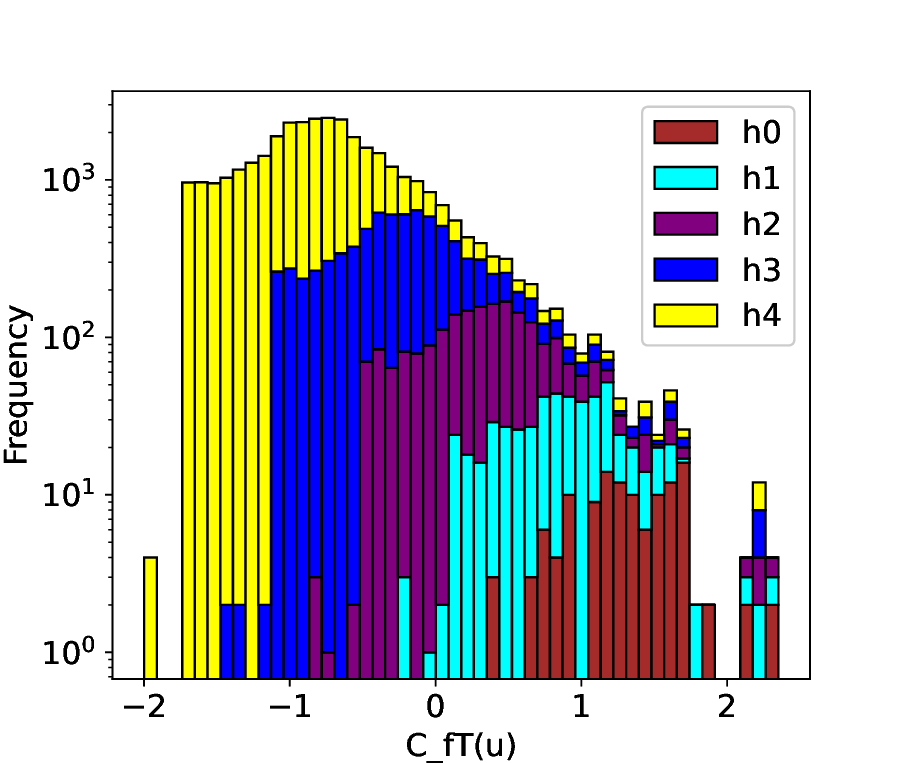}
    \hfill
    \begin{tikzpicture}[scale=0.525]
      \begin{loglogaxis}[legend pos=south east]
        \addplot table[col sep=comma,x=MeshSize,y=EnergyError] {results/test1/r4/hexa1/hexa1_k0_mu6_nu2/data_rates.dat}
        node[pos = 0.150, above=2pt]{2.84} % EnergyErrorRate
        node[pos = 0.475, above=2pt]{3.07} % EnergyErrorRate
        node[pos = 0.750, above=2pt]{1.75} % EnergyErrorRate
        node[pos = 0.925, above=2pt]{0.67}; % EnergyErrorRate
        \addlegendentry{$k=0$}

        \addplot table[col sep=comma,x=MeshSize,y=EnergyError] {results/test1/r4/hexa1/hexa1_k1_mu6_nu2/data_rates.dat}
        node[pos = 0.150, above=2pt]{4.69} % EnergyErrorRate
        node[pos = 0.475, above=2pt]{5.35} % EnergyErrorRate
        node[pos = 0.750, above=2pt]{4.69} % EnergyErrorRate
        node[pos = 0.925, above=2pt]{2.41}; % EnergyErrorRate
        \addlegendentry{$k=1$}

        \addplot table[col sep=comma,x=MeshSize,y=EnergyError] {results/test1/r4/hexa1/hexa1_k2_mu6_nu2/data_rates.dat}
        node[pos = 0.150, above=2pt]{6.97} % EnergyErrorRate
        node[pos = 0.475, above=2pt]{6.88} % EnergyErrorRate
        node[pos = 0.750, above=2pt]{5.10} % EnergyErrorRate
        node[pos = 0.925, above=2pt]{4.10}; % EnergyErrorRate
        \addlegendentry{$k=2$}

      \end{loglogaxis}
    \end{tikzpicture}
  }\vspace{0.1cm}\\
  \caption{Convergence results using the hexagonal mesh family for different values of $(r,\mu,\nu)$ and for $k\in\{0,1,2\}$.
    Left column: Distribution of $\CfT$ for all the meshes.
    Right column:  The quantity $\EnergyError[r]$ as a function of $h$ (the numbers on the graphs indicate the slopes between two data points).}
  \label{fig:nt1.hexa1}
\end{figure}
 %Regime test
\begin{figure}\centering
  \subcaptionbox[B]{$(\mu,\nu) = (1,10^{-5})$.\label{fig:(r1_75,1,10pm5)}}[0.95\textwidth]{\centering
\hspace*{\fill}
      %Tri
        \begin{tikzpicture}[scale=0.55]
          \begin{loglogaxis}[legend pos=south east]
              \addplot table[col sep=comma,x=MeshSize,y=EnergyError] {results/test1/r1_75/tri/tri_k0_mu1_nu7/data_rates.dat}
                  node[pos = 0.1, above=2pt]{2.04} % EnergyErrorRate
                  node[pos = 0.3, above=2pt]{2.01} % EnergyErrorRate
                  node[pos = 0.5, above=2pt]{2.00} % EnergyErrorRate
                  node[pos = 0.7, above=2pt]{2.00} % EnergyErrorRate
                  node[pos = 0.9, above=2pt]{2.00}; % EnergyErrorRate
               \addlegendentry{$k=0$}
               
              \addplot table[col sep=comma,x=MeshSize,y=EnergyError] {results/test1/r1_75/tri/tri_k1_mu1_nu7/data_rates.dat}
                  node[pos = 0.15, above=2pt]{3.90} % EnergyErrorRate
                  node[pos = 0.35, above=2pt]{3.94} % EnergyErrorRate
                  node[pos = 0.55, above=2pt]{3.95} % EnergyErrorRate
                  node[pos = 0.75, above=2pt]{3.96} % EnergyErrorRate
                  node[pos = 0.95, above=2pt]{3.97}; % EnergyErrorRate
               \addlegendentry{$k=1$}
               
              \addplot table[col sep=comma,x=MeshSize,y=EnergyError] {results/test1/r1_75/tri/tri_k2_mu1_nu7/data_rates.dat}
                  node[pos = 0.15, above=2pt]{4.89} % EnergyErrorRate
                  node[pos = 0.35, above=2pt]{4.60} % EnergyErrorRate
                  node[pos = 0.55, above=2pt]{4.42} % EnergyErrorRate
                  node[pos = 0.75, above=2pt]{4.34} % EnergyErrorRate
                  node[pos = 0.95, above=2pt]{4.31}; % EnergyErrorRate
               \addlegendentry{$k=2$}
               
          \end{loglogaxis}
        \end{tikzpicture}
        \hfill
        %Hexa
        \begin{tikzpicture}[scale=0.55]
          \begin{loglogaxis}[legend pos=south east]
              \addplot table[col sep=comma,x=MeshSize,y=EnergyError] {results/test1/r1_75/hexa/hexa1_k0_mu1_nu7/data_rates.dat}
                  node[pos = 0.125, above=2pt]{2.14} % EnergyErrorRate
                  node[pos = 0.375, above=2pt]{2.01} % EnergyErrorRate
                  node[pos = 0.625, above=2pt]{1.99} % EnergyErrorRate
                  node[pos = 0.875, above=2pt]{1.99}; % EnergyErrorRate
               \addlegendentry{$k=0$}
               
              \addplot table[col sep=comma,x=MeshSize,y=EnergyError] {results/test1/r1_75/hexa/hexa1_k1_mu1_nu7/data_rates.dat}
                  node[pos = 0.125, above=2pt]{4.22} % EnergyErrorRate
                  node[pos = 0.375, above=2pt]{4.09} % EnergyErrorRate
                  node[pos = 0.625, above=2pt]{4.07} % EnergyErrorRate
                  node[pos = 0.875, above=2pt]{4.04}; % EnergyErrorRate
               \addlegendentry{$k=1$}
               
              \addplot table[col sep=comma,x=MeshSize,y=EnergyError] {results/test1/r1_75/hexa/hexa1_k2_mu1_nu7/data_rates.dat}
                  node[pos = 0.125, above=2pt]{5.28} % EnergyErrorRate
                  node[pos = 0.375, above=2pt]{5.60} % EnergyErrorRate
                  node[pos = 0.625, above=2pt]{5.74} % EnergyErrorRate
                  node[pos = 0.875, above=2pt]{5.62}; % EnergyErrorRate
               \addlegendentry{$k=2$}
               
          \end{loglogaxis}
        \end{tikzpicture}
     \hspace*{\fill}
   }\vspace{0.1cm}\\
  \subcaptionbox[B]{$(\mu,\nu) = (1,1)$.\label{fig:(r1_75,1,1)}}[0.95\textwidth]{\centering
    \hspace*{\fill}
        %Tri
        \begin{tikzpicture}[scale=0.55]
          \begin{loglogaxis}[legend pos=south east]
              \addplot table[col sep=comma,x=MeshSize,y=EnergyError] {results/test1/r1_75/tri/tri_k0_mu1_nu2/data_rates.dat}
                  node[pos = 0.1, above=2pt]{2.04} % EnergyErrorRate
                  node[pos = 0.3, above=2pt]{2.01} % EnergyErrorRate
                  node[pos = 0.5, above=2pt]{2.00} % EnergyErrorRate
                  node[pos = 0.7, above=2pt]{2.00} % EnergyErrorRate
                  node[pos = 0.9, above=2pt]{2.00}; % EnergyErrorRate
               \addlegendentry{$k=0$}
               
              \addplot table[col sep=comma,x=MeshSize,y=EnergyError] {results/test1/r1_75/tri/tri_k1_mu1_nu2/data_rates.dat}
                  node[pos = 0.15, above=2pt]{3.90} % EnergyErrorRate
                  node[pos = 0.35, above=2pt]{3.94} % EnergyErrorRate
                  node[pos = 0.55, above=2pt]{3.96} % EnergyErrorRate
                  node[pos = 0.75, above=2pt]{3.97} % EnergyErrorRate
                  node[pos = 0.95, above=2pt]{3.97}; % EnergyErrorRate
               \addlegendentry{$k=1$}
               
              \addplot table[col sep=comma,x=MeshSize,y=EnergyError] {results/test1/r1_75/tri/tri_k2_mu1_nu2/data_rates.dat}
                  node[pos = 0.15, above=2pt]{4.95} % EnergyErrorRate
                  node[pos = 0.35, above=2pt]{4.64} % EnergyErrorRate
                  node[pos = 0.55, above=2pt]{4.45} % EnergyErrorRate
                  node[pos = 0.75, above=2pt]{4.35} % EnergyErrorRate
                  node[pos = 0.95, above=2pt]{4.31}; % EnergyErrorRate
               \addlegendentry{$k=2$}
               
          \end{loglogaxis}
        \end{tikzpicture}
       \hfill
        %Hexa
        \begin{tikzpicture}[scale=0.55]
          \begin{loglogaxis}[legend pos=south east]
              \addplot table[col sep=comma,x=MeshSize,y=EnergyError] {results/test1/r1_75/hexa/hexa1_k0_mu1_nu2/data_rates.dat}
                  node[pos = 0.125, above=2pt]{2.17} % EnergyErrorRate
                  node[pos = 0.375, above=2pt]{2.02} % EnergyErrorRate
                  node[pos = 0.625, above=2pt]{2.00} % EnergyErrorRate
                  node[pos = 0.875, above=2pt]{1.99}; % EnergyErrorRate
               \addlegendentry{$k=0$}
               
              \addplot table[col sep=comma,x=MeshSize,y=EnergyError] {results/test1/r1_75/hexa/hexa1_k1_mu1_nu2/data_rates.dat}
                  node[pos = 0.125, above=2pt]{4.24} % EnergyErrorRate
                  node[pos = 0.375, above=2pt]{4.09} % EnergyErrorRate
                  node[pos = 0.625, above=2pt]{4.07} % EnergyErrorRate
                  node[pos = 0.875, above=2pt]{4.04}; % EnergyErrorRate
               \addlegendentry{$k=1$}
               
              \addplot table[col sep=comma,x=MeshSize,y=EnergyError] {results/test1/r1_75/hexa/hexa1_k2_mu1_nu2/data_rates.dat}
                  node[pos = 0.125, above=2pt]{5.34} % EnergyErrorRate
                  node[pos = 0.375, above=2pt]{5.66} % EnergyErrorRate
                  node[pos = 0.625, above=2pt]{5.76} % EnergyErrorRate
                  node[pos = 0.875, above=2pt]{5.60}; % EnergyErrorRate
               \addlegendentry{$k=2$}
               
          \end{loglogaxis}
        \end{tikzpicture}
    \hspace*{\fill}
   }\vspace{0.1cm}\\
  \subcaptionbox[B]{$(\mu,\nu) = (10^{-3},1)$.\label{fig:(r1_75,10pm3,1)}}[0.95\textwidth]{\centering
    \hspace*{\fill}
        %Tri
        \begin{tikzpicture}[scale=0.55]
          \begin{loglogaxis}[legend pos=south east]
              \addplot table[col sep=comma,x=MeshSize,y=EnergyError] {results/test1/r1_75/tri/tri_k0_mu4_nu2/data_rates.dat}
                  node[pos = 0.1, above=2pt]{2.18} % EnergyErrorRate
                  node[pos = 0.3, above=2pt]{2.87} % EnergyErrorRate
                  node[pos = 0.5, above=2pt]{3.53} % EnergyErrorRate
                  node[pos = 0.7, above=2pt]{3.52} % EnergyErrorRate
                  node[pos = 0.9, above=2pt]{3.64}; % EnergyErrorRate
               \addlegendentry{$k=0$}
               
              \addplot table[col sep=comma,x=MeshSize,y=EnergyError] {results/test1/r1_75/tri/tri_k1_mu4_nu2/data_rates.dat}
                  node[pos = 0.1, above=2pt]{4.46} % EnergyErrorRate
                  node[pos = 0.3, above=2pt]{4.92} % EnergyErrorRate
                  node[pos = 0.5, above=2pt]{5.46} % EnergyErrorRate
                  node[pos = 0.7, above=2pt]{5.84} % EnergyErrorRate
                  node[pos = 0.9, above=2pt]{5.89}; % EnergyErrorRate
               \addlegendentry{$k=1$}
               
              \addplot table[col sep=comma,x=MeshSize,y=EnergyError] {results/test1/r1_75/tri/tri_k2_mu4_nu2/data_rates.dat}
                  node[pos = 0.15, above=2pt]{6.57} % EnergyErrorRate
                  node[pos = 0.35, above=2pt]{7.08} % EnergyErrorRate
                  node[pos = 0.55, above=2pt]{7.34} % EnergyErrorRate
                  node[pos = 0.75, above=2pt]{6.55} % EnergyErrorRate
                  node[pos = 0.95, above=2pt]{5.12}; % EnergyErrorRate
               \addlegendentry{$k=2$}
               
          \end{loglogaxis}
        \end{tikzpicture}
  \hfill
     %Hexa
      \begin{tikzpicture}[scale=0.55]
          \begin{loglogaxis}[legend pos=south east]
              \addplot table[col sep=comma,x=MeshSize,y=EnergyError] {results/test1/r1_75/hexa/hexa1_k0_mu4_nu2/data_rates.dat}
                  node[pos = 0.125, above=2pt]{2.80} % EnergyErrorRate
                  node[pos = 0.375, above=2pt]{3.31} % EnergyErrorRate
                  node[pos = 0.625, above=2pt]{3.90} % EnergyErrorRate
                  node[pos = 0.875, above=2pt]{4.89}; % EnergyErrorRate
               \addlegendentry{$k=0$}
               
              \addplot table[col sep=comma,x=MeshSize,y=EnergyError] {results/test1/r1_75/hexa/hexa1_k1_mu4_nu2/data_rates.dat}
                  node[pos = 0.125, above=2pt]{4.76} % EnergyErrorRate
                  node[pos = 0.375, above=2pt]{6.54} % EnergyErrorRate
                  node[pos = 0.625, above=2pt]{8.82} % EnergyErrorRate
                  node[pos = 0.875, above=2pt]{8.14}; % EnergyErrorRate
               \addlegendentry{$k=1$}
               
              \addplot table[col sep=comma,x=MeshSize,y=EnergyError] {results/test1/r1_75/hexa/hexa1_k2_mu4_nu2/data_rates.dat}
                  node[pos = 0.125, above=2pt]{6.93} % EnergyErrorRate
                  node[pos = 0.375, above=2pt]{8.18} % EnergyErrorRate
                  node[pos = 0.625, above=2pt]{10.93} % EnergyErrorRate
                  node[pos = 0.875, above=2pt]{8.74}; % EnergyErrorRate
               \addlegendentry{$k=2$}
               
          \end{loglogaxis}
        \end{tikzpicture}   \hspace*{\fill}
   }\vspace{0.1cm}\\
  \subcaptionbox[B]{$(\mu,\nu) = (10^{-6},1)$.\label{fig:(r1_75,10pm6,1)}}[0.95\textwidth]{\centering
    \hspace*{\fill}
        %Tri
         \begin{tikzpicture}[scale=0.55]
          \begin{loglogaxis}[legend pos=south east]
              \addplot table[col sep=comma,x=MeshSize,y=EnergyError] {results/test1/r1_75/tri/tri_k0_mu7_nu2/data_rates.dat}
                  node[pos = 0.1, above=2pt]{1.83} % EnergyErrorRate
                  node[pos = 0.3, above=2pt]{1.92} % EnergyErrorRate
                  node[pos = 0.5, above=2pt]{2.01} % EnergyErrorRate
                  node[pos = 0.7, above=2pt]{2.04} % EnergyErrorRate
                  node[pos = 0.9, above=2pt]{2.04}; % EnergyErrorRate
               \addlegendentry{$k=0$}
               
              \addplot table[col sep=comma,x=MeshSize,y=EnergyError] {results/test1/r1_75/tri/tri_k1_mu7_nu2/data_rates.dat}
                  node[pos = 0.1, above=2pt]{3.88} % EnergyErrorRate
                  node[pos = 0.3, above=2pt]{3.95} % EnergyErrorRate
                  node[pos = 0.5, above=2pt]{3.97} % EnergyErrorRate
                  node[pos = 0.7, above=2pt]{3.99} % EnergyErrorRate
                  node[pos = 0.9, above=2pt]{3.99}; % EnergyErrorRate
               \addlegendentry{$k=1$}
               
              \addplot table[col sep=comma,x=MeshSize,y=EnergyError] {results/test1/r1_75/tri/tri_k2_mu7_nu2/data_rates.dat}
                  node[pos = 0.15, above=2pt]{5.66} % EnergyErrorRate
                  node[pos = 0.35, above=2pt]{5.91} % EnergyErrorRate
                  node[pos = 0.55, above=2pt]{5.98} % EnergyErrorRate
                  node[pos = 0.75, above=2pt]{6.00} % EnergyErrorRate
                  node[pos = 0.95, above=2pt]{5.83}; % EnergyErrorRate
               \addlegendentry{$k=2$}
               
          \end{loglogaxis}
        \end{tikzpicture}
  \hfill
        %Hex5
        \begin{tikzpicture}[scale=0.55]
          \begin{loglogaxis}[legend pos=south east]
              \addplot table[col sep=comma,x=MeshSize,y=EnergyError] {results/test1/r1_75/hexa/hexa1_k0_mu7_nu2/data_rates.dat}
                  node[pos = 0.125, above=2pt]{2.32} % EnergyErrorRate
                  node[pos = 0.375, above=2pt]{2.49} % EnergyErrorRate
                  node[pos = 0.625, above=2pt]{2.64} % EnergyErrorRate
                  node[pos = 0.875, above=2pt]{2.78}; % EnergyErrorRate
               \addlegendentry{$k=0$}
               
              \addplot table[col sep=comma,x=MeshSize,y=EnergyError] {results/test1/r1_75/hexa/hexa1_k1_mu7_nu2/data_rates.dat}
                  node[pos = 0.125, above=2pt]{3.90} % EnergyErrorRate
                  node[pos = 0.375, above=2pt]{3.90} % EnergyErrorRate
                  node[pos = 0.625, above=2pt]{3.96} % EnergyErrorRate
                  node[pos = 0.875, above=2pt]{3.99}; % EnergyErrorRate
               \addlegendentry{$k=1$}
               
              \addplot table[col sep=comma,x=MeshSize,y=EnergyError] {results/test1/r1_75/hexa/hexa1_k2_mu7_nu2/data_rates.dat}
                  node[pos = 0.125, above=2pt]{6.09} % EnergyErrorRate
                  node[pos = 0.375, above=2pt]{6.00} % EnergyErrorRate
                  node[pos = 0.625, above=2pt]{6.00} % EnergyErrorRate
                  node[pos = 0.875, above=2pt]{6.00}; % EnergyErrorRate
               \addlegendentry{$k=2$}
               
          \end{loglogaxis}
        \end{tikzpicture}
     \hspace*{\fill}
   }\vspace{0.1cm}\\
    \caption{Convergence results for $r=1.75$ and different values of $(\mu,\nu)$ and for $k\in\{0,1,2\}$. 
    Left column: Triangular mesh family.
    Right column: Hexagonal mesh family.  In both columns we show the quantity $\EnergyError[2]$ as a function of $h$ (the numbers on the graphs indicate the slopes between two data points).}
\label{fig:ntest1:r1.75}
\end{figure}

%------------------------------------------------------------------------------%

\begin{table}
  \centering
  \resizebox{\textwidth}{!}{%
    \begin{tabular}{cccc|ccc|ccc|ccc}
      \toprule
      $h$      &
      Max & Min &  Avg &
      Max & Min &  Avg &
      Max & Min &  Avg &
      Max & Min &  Avg
      \\
      \midrule
      \multicolumn{1}{c}{}&
      \multicolumn{3}{c|}{Triangular,  $(\mu,\nu)=(1,10^{-5})$}&
      \multicolumn{3}{c|}{Triangular,  $(\mu,\nu)=(1,1)$}&
      \multicolumn{3}{c}{Triangular,  $(\mu,\nu)=(10^{-3},1)$}&
      \multicolumn{3}{c}{Triangular,  $(\mu,\nu)=(10^{-6},1)$}    
      \\
      \midrule
      $h_0$    &  1.1E-06  & 6.3E-07 & 8.8E-07   &  1.0E-01 & 5.5E-02 & 7.8E-02       & 1.1E+02 & 6.2E+01 & 8.7E+01  &  1.1E+05 & 6.2E+04 & 8.7E+04 \\
      $h_1$    &  3.0E-07 & 1.3E-07 & 2.2E-07   &  2.6E-02 & 1.2E-02 & 1.9E-02        & 2.9E+01 & 1.3E+01 & 2.2E+01  &  2.9E+04 & 1.3E+04 & 2.2E+04 \\
      $h_2$    &  7.6E08  & 2.9E-08 & 5.7E-08   &  6.7E-03 & 2.5E-03 & 5.0E-03        & 7.5E+00 & 2.8E+00 & 5.6E+00  &  7.5E+03 & 2.8E+03 & 5.6E+03 \\
      $h_3$    &  1.9E08  & 6.1E-09 & 1.4E-08   &  1.6E-03 & 5.4E-04 & 1.2E-03        & 1.8E+00 & 0.6E-01 & 1.4E+00  &  1.8E+03 & 6.0E+02 & 1.4E+03 \\
      $h_4$    &  4.7E-09  & 1.2E-09 & 3.5E-09   &  4.2E-04 & 1.1E-04 & 3.1E-04       &  4.7E-01 & 1.2E-01 & 3.5E-01 &  4.7E+02 & 1.2E+02 & 3.5E+02 \\
      $h_5$    &  1.1E-09  & 2.7E-10 & 8.9E-10   &  1.0E-04 & 2.4E-05 & 7.9E-05       & 1.1E-01 &  2.6E-02 & 8.8E-02 &  1.1E+02 & 2.6E+01 & 8.8E+01 \\
      \midrule
      \multicolumn{1}{c}{}&
      \multicolumn{3}{c|}{Hexagonal,  $(\mu,\nu)=(1,10^{-5})$}&
      \multicolumn{3}{c|}{Hexagonal,  $(\mu,\nu)=(1,1)$}&
      \multicolumn{3}{c}{Hexagonal,  $(\mu,\nu)=(10^{-3},1)$}&
      \multicolumn{3}{c}{Hexagonal,  $(\mu,\nu)=(10^{-6},1)$}    
      \\
      \midrule
      $h_0$    &  1.0E-06 &  9.7E-08 & 4.2E-07    &  9.1E-02 &  8.6E-03 & 3.7E-02    &  1.0E+02 & 9.6E+00 & 4.1E+01 & 1.0E+05 & 9.6E+03 & 4.1E+04   \\
      $h_1$    &  2.9E-07 &  2.4E-08 & 1.2E-07    &  2.5E-02 &  2.1E-03 & 1.1E-02    &  2.9E+01 & 2.4E+00 & 1.2E+01 & 2.9E+04 & 2.4E+03 & 1.2E+04   \\
      $h_2$    &  7.7E-08 &  5.7E-09 & 3.3E-08    &  6.8E-03 &  5.1E-04 & 2.9E-03    &  7.6E+00 & 5.7E-01 & 3.3E+00 & 7.6E+03 & 5.7E+02 & 3.3E+03   \\
      $h_3$    &  1.9E-08 &  1.2E-09 & 8.6E-09    &  1.7E-03 &  1.1E-04 & 7.6E-04    &  1.9E+00 & 1.2E-01 & 8.5E-01 & 1.9E+03 & 1.2E+02 &  8.5E+02  \\
      $h_4$    &  4.9E-09 &  2.6E-10 & 2.1E-09    &  4.3E-04 &  2.3E-05 & 1.9E-04    &  4.8E-01 & 2.6E-02 & 2.1E-01 & 4.8E+02 & 2.6E+01 & 2.1E+02   \\
      \bottomrule
    \end{tabular}
  }
  \caption{Maximum, minimum and average $\CfT$ values for $r=1.75$  and $(\mu,\nu) \in \{ (1,10^{-5}), (1,1), (10^{-3},1), (10^{-6},1) \}$ on the triangular and hexagonal mesh families. 
    \label{tbl:ntest1:r1.75:1}
  }
\end{table}

 %Regime test

  Figures~\ref{fig:nt1.tri} and~\ref{fig:nt1.hexa1} display in the left column the frequency distribution of $\CfT$ for each mesh  (using logarithm scaling in both axes)
where each color corresponds to a different mesh with mesh size $h_i$.
Recalling the definition \eqref{eq:CfT} of $\CfT$, if the frequency distribution of $\CfT$  for a given mesh is concentrated to the left (resp., right) of $0$, then most elements are in the Stokes-dominated (resp., Darcy-dominated) regime.

  In the right column of Figures~\ref{fig:nt1.tri} and~\ref{fig:nt1.hexa1} we show the convergence of the error for polynomial degrees $k$ ranging from $0$ to $2$.
  A first general remark concerning the triangular mesh is that we observe superconvergence for $k = 0$.
  This phenomenon appears to be linked to the use of the Darcy potential reconstruction in the right-hand side of \eqref{eq:discrete:momentum}, and is not observed when the element component of the test function is used instead.
  We postpone the theoretical study of this behaviour to a future work.
  On hexagonal meshes, on the other hand, a decreased convergence rate is sometimes observed for $k = 0$ in Figure \ref{fig:nt1.hexa1}.
  A similar behaviour was observed in \cite{Anderson.Droniou:18} for a different non-linear model representing miscible flows in porous media, with a rapid improvement starting from $k = 1$.
  A possible explanation is that, while the HHO space does have the required consistency properties, it is still quite ``rigid'' for $k = 0$ when the incompressibility constraint is accounted for.
  Aside from the case $k = 0$ commented above, the numerical results in Figures~\ref{fig:nt1.tri} and~\ref{fig:nt1.hexa1} show good agreement with the theory.
  Specifically, in Figures~\ref{fig:tri:(3,1,1)},~\ref{fig:hexa:(3,1e-1,1)} and~\ref{fig:hexa:(4,1e-2,1)}, we observe orders of convergence close to the ones in the left panel of Table~\ref{tab:convergence.rates:numerical.tests}, consistently with the fact that $\CfT \ll 1$ for most $T \in \Th$.
  Slope changes are observed in Figures~\ref{fig:tri:(3,1e-3,1)},~\ref{fig:tri:(4,1e-5,1)},~\ref{fig:hexa:(3,1e-4,1)}, and~\ref{fig:hexa:(4,1e-5,1)}, where we start with convergence rates close to or slightly higher than the ones expected in the Darcy-dominated regime (see the right panel of Table~\ref{tab:convergence.rates:numerical.tests}) and finish with convergence rates close to the ones expected in the Stokes-dominated regime.
  This behaviour is consistent with the fact that the distribution of $\CfT$ crosses the threshold $1$ identifying the two regimes as the mesh refinement progresses, and shows that the Stokes regimes is eventually observed when the mesh size is small enough, as explained in Remark~\ref{rem:conv.rates}.
  Finally, orders of convergence corresponding to the Darcy-dominated regime are observed in Figure~\ref{fig:tri:(3,1e-5,1)}, showing that, in practical situations, one may not achieve the asymptotic regime even for very fine meshes, corroborating the need of error estimates accounting for pre-asymptotic regimes.
  In addition, we observe that the computed convergence rate for the case $k = 1$ in Figure~\ref{fig:hexa:(4,1e-5,1)} for the finest  hexagonal mesh drops to 2.41, slightly less than the theoretical value 2.65.
  This behaviour probably can be improved using specialized convex numerical optimization algorithms, but we leave this issue for future research.
  
  Figure~\ref{fig:ntest1:r1.75} shows the convergence rates for the case $r=1.75$ (using Picard's algorithm instead of Newton's), and different values of $(\mu,\nu)$.  
  We provide in Table \ref{tbl:ntest1:r1.75:1} the maximum, minimum, and average values $\CfT$ values for each mesh and each case shown in Figure~\ref{fig:ntest1:r1.75}.
  It is observed that in almost all cases, the convergence rates reported in Figure~\ref{fig:ntest1:r1.75} exceed those predicted by the theory (see the Table~\ref{tab:convergence.rates:numerical.tests}), with the exception of those in Figures~\ref{fig:(r1_75,1,10pm5)} and \ref{fig:(r1_75,1,1)} for triangular meshes with $k=2$ (these cases are Stokes dominated, see Table \ref{tbl:ntest1:r1.75:1}).
  However, the observed rates are still closed to the expected value $4.5$, the difference being most likely due to machine precision limitation.
  We have also run tests for $r \in \{1.25, 1.5\}$ (not reported here), values for which better convergence rates than predicted by the theory are observed in practice.
  This improved convergence is most likely linked to the fact that we are testing against a smooth solution.
  In the context of the scalar $p$-Laplace equation, improved estimates for HHO discretisations can indeed be derived when the solution does not lead to flux degeneration~\cite{Di-Pietro.Droniou.ea:21}.

%------------------------------------------------------------------------------%

\bibliographystyle{siamplain}
\bibliography{p-brinkman-poly}

@Article{         Allaire:91,
  author        = {Allaire, G.},
  title         = {Homogenization of the {Navier--Stokes} equations in open
                  set perforated with tiny holes. Part {I}: Abstract
                  framework, a volume distribution of holes},
  journal       = {Arch. Rat. Mech. Anal.},
  volume        = {113},
  pages         = {209--259},
  year          = {1991},
  doi           = {10.1007/BF00375065}
}

@Article{         Allaire:91*1,
  author        = {Allaire, G.},
  title         = {Homogenization of the {Navier--Stokes} equations in open
                  set perforated with tiny holes. Part {II}: Non-critical
                  size of the holes for a volume distribution and a surface
                  distribution of holes},
  journal       = {Arch. Rat. Mech. Anal.},
  volume        = {113},
  pages         = {261--298},
  year          = {1991},
  doi           = {10.1007/BF00375066}
}

@InBook{          Allaire:91*2,
  author        = {Allaire, G.},
  chapter       = {Homogenization of the {Navier--Stokes} equations and
                  derivation of {Brinkman}'s law},
  series        = {Applied Mathematics for Engineering Sciences},
  editor        = {C. Carasso et al.},
  pages         = {7--20},
  publisher     = {Cépaduès Editions},
  address       = {Toulouse},
  year          = {1991}
}

@Article{         Anaya.Gatica.ea:15,
  author        = {Anaya, V. and Gatica, G. N. and Mora, D. and Ruiz-Baier,
                  R.},
  title         = {An augmented velocity-vorticity-pressure formulation for
                  the {B}rinkman equations},
  journal       = {Internat. J. Numer. Methods Fluids},
  volume        = {79},
  year          = {2015},
  number        = {3},
  pages         = {109--137},
  doi           = {10.1002/fld.4041}
}

@Article{         Anderson.Droniou:18,
  title         = {An arbitrary order scheme on generic meshes for miscible
                  displacements in porous media},
  author        = {Anderson, D. and Droniou, J.},
  journal       = {SIAM J. Sci. Comput.},
  volume        = {40},
  number        = {4},
  pages         = {B1020--B1054},
  year          = {2018},
  doi           = {10.1137/17M1138807}
}

@Article{         Anguiano.Suarez-Grau:21,
  author        = {Anguiano, M. and Su\'{a}rez-Grau, F. J.},
  title         = {Lower-dimensional nonlinear {Brinkman}’s law for
                  {non-Newtonian} flows in a thin porous medium},
  journal       = {Mediterranean Journal of Mathematics},
  year          = {2021},
  volume        = {18},
  number        = {175},
  doi           = {10.1007/s00009-021-01814-5}
}

@Article{         Araya.Harder.ea:17,
  author        = {Araya, R. and Harder, C. and Poza, A. H. and Valentin,
                  F.},
  title         = {Multiscale hybrid-mixed method for the {S}tokes and
                  {B}rinkman equations---the method},
  journal       = {Comput. Methods Appl. Mech. Engrg.},
  volume        = {324},
  year          = {2017},
  pages         = {29--53},
  doi           = {10.1016/j.cma.2017.05.027}
}

@Article{         Barrett.Liu:94,
  title         = {Quasi-norm error bounds for the finite element
                  approximation of a non-Newtonian flow},
  author        = {Barrett, J.~W. and Liu, W.~B.},
  journal       = {Numer. Math.},
  volume        = {68},
  number        = {4},
  pages         = {437--456},
  year          = {1994}
}

@Article{         Beirao-da-Veiga.Di-Pietro.ea:24,
  author        = {Beir\~{a}o da Veiga, L. and Di Pietro, D. A. and Haile, K.
                  B.},
  title         = {A {P\'{e}clet}-robust discontinuous {Galerkin} method for
                  nonlinear diffusion with advection},
  year          = {2024},
  journal       = {Math. Models Methods Appl. Sci.},
  volume        = {34},
  number        = {9},
  pages         = {1781--1807},
  doi           = {10.1142/S0218202524500350}
}

@Article{         Belenki.Berselli.ea:12,
  author        = {Belenki, L. and Berselli, L. C. and Diening, L. and
                  R\r{u}\v{z}i\v{c}ka, M.},
  title         = {On the finite element approximation of {$p$}-{S}tokes
                  systems},
  journal       = {SIAM J. Numer. Anal.},
  volume        = {50},
  year          = {2012},
  number        = {2},
  pages         = {373--397},
  doi           = {10.1137/10080436X}
}

@Article{         Botti.Castanon-Quiroz.ea:21,
  author        = {Botti, M. and Casta\~{n}\'{o}n Quiroz, D. and Di Pietro,
                  D. A. and Harnist, A.},
  title         = {A {Hybrid High-Order} method for creeping flows of
                  non-{Newtonian} fluids},
  journal       = {ESAIM: Math. Model Numer. Anal.},
  year          = {2021},
  volume        = {55},
  number        = {5},
  pages         = {2045--2073},
  doi           = {10.1051/m2an/2021051}
}

@Article{         Botti.Di-Pietro.ea:17,
  author        = {Botti, M. and Di Pietro, D. A. and Sochala, P.},
  title         = {A {Hybrid High-Order} method for nonlinear elasticity},
  year          = {2017},
  journal       = {SIAM J. Numer. Anal.},
  volume        = {55},
  number        = {6},
  pages         = {2687--2717},
  doi           = {10.1137/16M1105943}
}

@Article{         Botti.Di-Pietro.ea:18,
  author        = {Botti, Lorenzo and Di Pietro, Daniele A. and Droniou,
                  J\'{e}r\^{o}me},
  title         = {A {H}ybrid {H}igh-{O}rder discretisation of the {B}rinkman
                  problem robust in the {D}arcy and {S}tokes limits},
  journal       = {Comput. Methods Appl. Mech. Engrg.},
  volume        = {341},
  year          = {2018},
  pages         = {278--310},
  doi           = {10.1016/j.cma.2018.07.004}
}

@Article{         Botti.Di-Pietro.ea:19*3,
  author        = {Botti, M. and Di Pietro, D. A. and Guglielmana, A.},
  title         = {A low-order nonconforming method for linear elasticity on
                  general meshes},
  journal       = {Comput. Meth. Appl. Mech. Engrg.},
  year          = {2019},
  volume        = {354},
  pages         = {96--118},
  doi           = {10.1016/j.cma.2019.05.031}
}

@Article{         Burman.Hansbo:05,
  author        = {Burman, Erik and Hansbo, Peter},
  title         = {Stabilized {C}rouzeix-{R}aviart element for the
                  {D}arcy-{S}tokes problem},
  journal       = {Numer. Methods Partial Differential Equations},
  volume        = {21},
  year          = {2005},
  number        = {5},
  pages         = {986--997},
  doi           = {10.1002/num.20076}
}

@Article{         Burman.Hansbo:07,
  author        = {Burman, E. and Hansbo, P.},
  title         = {A unified stabilized method for {S}tokes' and {D}arcy's
                  equations},
  journal       = {J. Comput. Appl. Math.},
  volume        = {198},
  year          = {2007},
  number        = {1},
  pages         = {35--51},
  doi           = {10.1016/j.cam.2005.11.022}
}

@Article{         Caceres.Gatica.ea:17,
  title         = {A mixed virtual element method for the {Brinkman}
                  problem},
  author        = {C\'{a}ceres, E. and Gatica, G. N. and Sequeira, F. A.},
  journal       = {Math. Models Methods Appl. Sci.},
  volume        = {27},
  number        = {4},
  pages         = {707--743},
  year          = {2017}
}

@Article{         Di-Pietro.Droniou.ea:21,
  author        = {Di Pietro, D. A. and Droniou, J. and Harnist, A.},
  title         = {Improved error estimates for {Hybrid High-Order}
                  discretizations of {Leray--Lions} problems},
  year          = {2021},
  journal       = {Calcolo},
  volume        = {58},
  number        = {19},
  doi           = {10.1007/s10092-021-00410-z}
}

@Article{         Di-Pietro.Droniou.ea:24,
  title         = {A pressure-robust {Discrete de Rham} scheme for the
                  {Navier--Stokes} equations},
  author        = {Di Pietro, D. A. and Droniou, J. and Qian, J. J.},
  journal       = {Comput. Meth. Appl. Mech. Engrg.},
  year          = {2024},
  volume        = {421},
  number        = {116765},
  doi           = {10.1016/j.cma.2024.116765}
}

@Article{         Di-Pietro.Droniou:17,
  author        = {Di Pietro, D. A. and Droniou, J.},
  title         = {A {Hybrid High-Order} method for {Leray--Lions} elliptic
                  equations on general meshes},
  journal       = {Math. Comp.},
  year          = {2017},
  volume        = {86},
  number        = {307},
  pages         = {2159--2191},
  doi           = {10.1090/mcom/3180}
}

@Article{         Di-Pietro.Droniou:17*1,
  author        = {Di Pietro, D. A. and Droniou, J.},
  title         = {{$W^{s,p}$}-approximation properties of elliptic
                  projectors on polynomial spaces, with application to the
                  error analysis of a {Hybrid High-Order} discretisation of
                  {Leray--Lions} problems},
  year          = {2017},
  volume        = {27},
  number        = {5},
  pages         = {879--908},
  journal       = {Math. Models Methods Appl. Sci.},
  doi           = {10.1142/S0218202517500191}
}

@Article{         Di-Pietro.Droniou:18,
  author        = {Di Pietro, D. A. and Droniou, J.},
  title         = {A third {Strang} lemma for schemes in fully discrete
                  formulation},
  year          = {2018},
  journal       = {Calcolo},
  volume        = {55},
  number        = {40},
  doi           = {10.1007/s10092-018-0282-3}
}

@Book{            Di-Pietro.Droniou:20,
  author        = {Di Pietro, D. A. and Droniou, J.},
  title         = {The {Hybrid High-Order} method for polytopal meshes},
  subtitle      = {Design, analysis, and applications},
  publisher     = {Springer International Publishing},
  year          = {2020},
  series        = {Modeling, Simulation and Application},
  number        = {19},
  doi           = {10.1007/978-3-030-37203-3}
}

@Article{         Di-Pietro.Droniou:23,
  title         = {A polytopal method for the {Brinkman} problem robust in
                  all regimes},
  author        = {Di Pietro, D. A. and Droniou, J.},
  journal       = {Comput. Meth. Appl. Mech. Engrg.},
  year          = {2023},
  volume        = {409},
  number        = {115981},
  doi           = {10.1016/j.cma.2023.115981}
}

@Article{         Di-Pietro.Droniou:23*1,
  author        = {Di Pietro, Daniele A. and Droniou, J\'er\^ome},
  title         = {An arbitrary-order discrete {de Rham} complex on
                  polyhedral meshes: {Exactness}, {Poincar\'e} inequalities,
                  and consistency},
  journal       = {Found. Comput. Math.},
  volume        = {23},
  pages         = {85--164},
  year          = {2023},
  doi           = {10.1007/s10208-021-09542-8}
}

@Article{         Di-Pietro.Ern.ea:16,
  author        = {Di Pietro, D. A. and Ern, A. and Linke, A. and Schieweck,
                  F.},
  title         = {A discontinuous skeletal method for the
                  viscosity-dependent {Stokes} problem},
  journal       = {Comput. Meth. Appl. Mech. Engrg.},
  year          = {2016},
  volume        = {306},
  pages         = {175--195},
  doi           = {10.1016/j.cma.2016.03.033}
}

@Article{         Di-Pietro.Ern:15,
  author        = {Di Pietro, D. A. and Ern, A.},
  title         = {A hybrid high-order locking-free method for linear
                  elasticity on general meshes},
  journal       = {Comput. Meth. Appl. Mech. Engrg.},
  year          = {2015},
  volume        = {283},
  pages         = {1--21},
  doi           = {10.1016/j.cma.2014.09.009}
}

@Article{         Di-Pietro.Ern:17,
  author        = {Di Pietro, D. A. and Ern, A.},
  title         = {Arbitrary-order mixed methods for heterogeneous
                  anisotropic diffusion on general meshes},
  journal       = {IMA J. Numer. Anal.},
  year          = {2017},
  volume        = {37},
  number        = {1},
  pages         = {40--63},
  doi           = {10.1093/imanum/drw003}
}

@Article{         Diening.Hirn.ea:25,
  author        = {Diening, Lars and Hirn, Adrian and Kreuzer, Christian and
                  Zanotti, Pietro},
  title         = {Pressure-robust finite element discretizations of the
                  nonlinear {S}tokes equations},
  journal       = {Math. Models Methods Appl. Sci.},
  volume        = {35},
  year          = {2025},
  number        = {12},
  pages         = {2661--2693},
  doi           = {10.1142/S0218202525500484}
}

@Article{         Diening.Koner.ea:14,
  author        = {Diening, Lars and K\"{o}ner, Dietmar and Ru\v{z}i\v{c}ka,
                  Michael and Toulopoulos, Ioannis},
  title         = {A local discontinuous {G}alerkin approximation for systems
                  with {$p$}-structure},
  journal       = {IMA J. Numer. Anal.},
  volume        = {34},
  year          = {2014},
  number        = {4},
  pages         = {1447--1488},
  doi           = {10.1093/imanum/drt040}
}

@Article{         Evans.Hughes:13,
  author        = {Evans, J. A. and Hughes, T. J. R.},
  title         = {Isogeometric divergence-conforming {B}-splines for the
                  {D}arcy-{S}tokes-{B}rinkman equations},
  journal       = {Math. Models Methods Appl. Sci.},
  volume        = {23},
  year          = {2013},
  number        = {4},
  pages         = {671--741},
  doi           = {10.1142/S0218202512500583}
}

@Article{         Gatica.Munar.ea:18,
  author        = {Gatica, G. N. and Munar, M. and Sequeira, F. A.},
  title         = {A mixed virtual element method for a nonlinear {Brinkman}
                  model of porous media flow},
  volume        = {55},
  number        = {21},
  year          = {2018},
  doi           = {10.1007/s10092-018-0262-7}
}

@Misc{            Guennebaud.Jacob.ea:10,
  author        = {G. Guennebaud and B. Jacob and others},
  title         = {Eigen v3},
  howpublished  = {http://eigen.tuxfamily.org},
  year          = {2010}
}

@Article{         Juntunen.Stenberg:10,
  author        = {Juntunen, M. and Stenberg, R.},
  title         = {Analysis of finite element methods for the {B}rinkman
                  problem},
  journal       = {Calcolo},
  volume        = {47},
  year          = {2010},
  number        = {3},
  pages         = {129--147},
  doi           = {10.1007/s10092-009-0017-6}
}

@Article{         Konno.Stenberg:11,
  author        = {K\"onn\"o, J. and Stenberg, R.},
  title         = {{$H ({\rm div})$}-conforming finite elements for the
                  {B}rinkman problem},
  journal       = {Math. Models Methods Appl. Sci.},
  volume        = {21},
  year          = {2011},
  number        = {11},
  pages         = {2227--2248},
  doi           = {10.1142/S0218202511005726}
}

@Article{         Mardal.Tai.ea:02,
  author        = {Mardal, K. A. and Tai, X.-C. and Winther, R.},
  title         = {A robust finite element method for {D}arcy-{S}tokes flow},
  journal       = {SIAM J. Numer. Anal.},
  volume        = {40},
  year          = {2002},
  number        = {5},
  pages         = {1605--1631},
  doi           = {10.1137/S0036142901383910}
}

@Article{         Mora.Reales.ea:22,
  author        = {Mora, D. and Reales, C. and Silgado, A.},
  title         = {A $C^1$-virtual element method of high order for the
                  {Brinkman} equations in stream function formulation with
                  pressure recovery},
  journal       = {IMA J. Numer. Anal.},
  volume        = {42},
  pages         = {3632--3674},
  year          = {2022},
  doi           = {10.1093/imanum/drab078}
}

@Article{         Schenk.Gartner.ea:01,
  author        = {Schenk, O. and G\"{a}rtner,K. and Fichtner, W.~ and
                  Stricker, A.},
  journal       = { Future Gener. Comput. Syst.},
  pages         = {69--78},
  title         = {{Pardiso: A high-performance serial and parallel sparse
                  linear solver in semiconductor device simulation}},
  volume        = {18},
  number        = {1},
  year          = {2001},
  doi           = {10.1016/S0167-739X(00)00076-5}
}

@Article{         Vacca:18,
  author        = {Vacca, G.},
  title         = {An {$H^1$}-conforming virtual element for {D}arcy and
                  {B}rinkman equations},
  journal       = {Math. Models Methods Appl. Sci.},
  volume        = {28},
  year          = {2018},
  number        = {1},
  pages         = {159--194},
  doi           = {10.1142/S0218202518500057}
}

@Article{         Zhao.Chung.ea:20,
  author        = {Zhao, Lina and Chung, Eric and Lam, Ming Fai},
  title         = {A new staggered {DG} method for the {Brinkman} problem
                  robust in the {Darcy} and {Stokes} limits},
  journal       = {Comput. Meth. Appl. Mech. Engrg.},
  volume        = {364},
  number        = {112986},
  year          = {2020},
  doi           = {10.1016/j.cma.2020.112986}
}

\end{document}